\UseRawInputEncoding 
\documentclass[letterpaper,10pt]{amsart}
\usepackage{amsmath,amssymb,amsxtra,amsthm, amstext, amscd,amsfonts,fancyhdr, hyperref,enumerate,tikz}
\usepackage[all,cmtip]{xy}
\usepackage[OT2,T1]{fontenc}
\DeclareSymbolFont{cyrletters}{OT2}{wncyr}{m}{n}
\DeclareMathSymbol{\Sha}{\mathalpha}{cyrletters}{"58}

\newcommand{\bC}{{\mathbb{C}}}

\newcommand{\bN}{{\mathbb{N}}}

\newcommand{\bP}{{\mathbb{P}}}
\newcommand{\bQ}{{\mathbb{Q}}}
\newcommand{\bR}{{\mathbb{R}}}

\newcommand{\bT}{{\mathbb{T}}}

\newcommand{\bX}{{\mathbb{X}}}

\newcommand{\bZ}{{\mathbb{Z}}}

\newcommand{\Ba}{{\mathbf{a}}}
\newcommand{\Bb}{{\mathbf{b}}}
\newcommand{\Bc}{{\mathbf{c}}}
\newcommand{\Bd}{{\mathbf{d}}}
\newcommand{\Be}{{\mathbf{e}}}

\newcommand{\Bi}{{\mathbf{i}}}

\newcommand{\Bm}{{\mathbf{m}}}

\newcommand{\Bp}{{\mathbf{p}}}

\newcommand{\Bv}{{\mathbf{v}}}

\newcommand{\Bx}{{\mathbf{x}}}
\newcommand{\By}{{\mathbf{y}}}

  \newcommand{\C}{{\mathcal{C}}}
  \newcommand{\D}{{\mathcal{D}}}
  \newcommand{\E}{{\mathcal{E}}}

\renewcommand{\L}{{\mathcal{L}}}
  
  \newcommand{\N}{{\mathcal{N}}}
\renewcommand{\O}{{\mathcal{O}}}

  \newcommand{\Q}{{\mathcal{Q}}}
  
\renewcommand{\S}{{\mathcal{S}}}
  \newcommand{\T}{{\mathcal{T}}}

  \newcommand{\Y}{{\mathcal{Y}}}
  
  \newcommand{\HH}{\mathcal{H}}

\newcommand{\fp}{\mathfrak{p}}

\newcommand{\fX}{\mathfrak{X}}

\newcommand{\Pic}{\operatorname{Pic}}

\newcommand{\BA}{\mathbf{A}}

\newcommand{\ord}{\operatorname{ord}}

\newcommand{\PGL}{\operatorname{PGL}}

\newcommand{\sqf}{\operatorname{sqf}}

\newcommand{\ep}{\varepsilon}

\newcommand{\Jac}{\operatorname{Jac}}

\newcommand{\ol}{\overline}
\newcommand{\ul}{\underline}

\newcommand{\upchi}{{\raise.35ex\hbox{$\chi$}}}

\newcommand{\rad}{\operatorname{rad}}

\newcommand{\ra}{\rightarrow}

\makeatletter
\@namedef{subjclassname@2010}{%
  \textup{2010} Mathematics Subject Classification}
\makeatother

\newtheorem{theorem}{Theorem}[section]
\newtheorem{corollary}[theorem]{Corollary}
\newtheorem{proposition}[theorem]{Proposition}
\newtheorem{lemma}[theorem]{Lemma}
\newtheorem{conjecture}[theorem]{Conjecture}

\theoremstyle{definition}
\newtheorem{definition}[theorem]{Definition}
\newtheorem{question}[theorem]{Question}
\newtheorem{example}[theorem]{Example}
\newtheorem{remark}[theorem]{Remark}

\usepackage{color}
\newif\ifhascomments \hascommentstrue
\ifhascomments
  
  \newcommand{\brett}[1]{{\color{blue}[[\ensuremath{\heartsuit\heartsuit\heartsuit} #1]]}}
\else
  \newcommand{\matt}[1]{}
  \newcommand{\brett}[1]{}
\fi

\numberwithin{equation}{section}

\begin{document}

\title{Heights and quantitative arithmetic on stacky curves}

\author{Brett Nasserden}
\address{Department of Pure Mathematics \\
University of Waterloo }
\email{bnasserd@uwaterloo.ca}
\indent

\author{Stanley Yao Xiao}
\address{Department of Mathematics \\
University of Toronto \\
Bahen Centre \\
40 St. George Street, Room 6290 \\
Toronto, Ontario, Canada \\  M5S 2E4 }
\email{syxiao@math.toronto.edu}
\indent


\begin{abstract} In this paper we investigate a family of algebraic stacks, the so-called \emph{stacky curves}, in the context of the general theory of heights on algebraic stacks due to Ellenberg, Satriano, and Zureick-Brown. We first give an elementary construction of a height which is seen to be dual to theirs. Next we count rational points having bounded E-S-ZB height on a particular stacky curve, answering a question of Ellenberg, Satriano, and Zureick-Brown. We then show that when the Euler characteristic of stacky curves is non-positive, that the E-S-ZB height coming from the anti-canonical divisor class fails to have the Northcott property. Next we prove a generalized version of a conjecture of Vojta, applied to stacky curves with negative Euler characteristic and coarse space $\bP^1$, is equivalent to the $abc$-conjecture. Finally, we prove that in the negative characteristic case the purely ``stacky" part of the E-S-ZB height exhibits the Northcott property.
\end{abstract}

\maketitle

\section{Introduction}

Two of the outstanding conjectures in number theory are the so-called Manin-Batyrev conjecture \cite{BatmanConj} for the density of rational points on open subschemes of Fano varieties with respect to a Weil height, and Malle's conjecture \cite{MalleConj} on the number of number fields having bounded discriminant, fixed degree, and fixed Galois group. Both conjectures assert, roughly, that the number of objects to be counted with an appropriate height at most $X$ satisfy an asymptotic formula of the form 
\[C \cdot X^\alpha (\log X)^\beta, \]
where $C, \alpha, \beta$ are non-negative numbers with $C, \alpha > 0$, and that $C, \alpha, \beta$ can be computed explicitly within their respective geometric and arithmetic frameworks. \\

In a recent article J.~Ellenberg, M.~Satriano, and D.~Zureick-Brown formulate a bold conjecture that encompasses both the Manin-Batyrev and Malle conjetures as special cases  \cite[Main Conjecture]{ESZ-B}. Their conjecture concerns counting rational points with respect to a new theory of heights applying broadly to \emph{algebraic stacks} which extends Weil's theory of heights on varieties. While the Manin and Malle conjectures are well studied, comparatively little is known about the behaviour of rational points on algebraic stacks.\\

In this article, we study heights on the stacky analogue of a smooth projective algebraic curve defined over $\mathbb{Q}$, that is a \emph{stacky curve}. Just as working with algebraic curves greatly simplifies the theory of algebraic geometry, so does working with stacky curves greatly simplify the theory of heights on algebraic stacks. The theory of heights on our stacky curves has the benefit of being completely explicit, and can be understood in an elementary (but not easy) manner. Indeed, we will see that natural questions involving stacky curves leads to an equivalent formulation of the $abc$-conjecture.  \\

 Let $\fX$ be a ``nice" algebraic stack defined over a number field $K$. Ellenberg, Satriano, and Zureick-Brown show that there is a function

\[\textnormal{H}_\fX\colon\textnormal{Vec}(\fX)\ra \{\textnormal{functions }\fX(K)\ra \bR\},\]
where $\textnormal{Vec}(\fX)$ is the collection of isomorphism classes of finite rank \emph{vector bundles} on $\fX$, with the property that whenever $\fX=X$ is a projective algebraic variety then $\textnormal{H}_X$ has the following properties. 

\begin{enumerate}
    \item The restriction $\textnormal{H}_X\mid_{\Pic(X)}$ is Weil's height machine in the usual sense.
    \item If $\E$ is a vector bundle on $X$ then $\textnormal{H}_X(\E)=\textnormal{H}_X(\det E)$.
\end{enumerate}

We call the above construction the \emph{Ellenberg-Satriano-Zureick-Brown Height machine}\label{thm:HeightMachine}. The first point says the E-S-ZB Height machine recovers Weil's height machine when applied to an algebraic variety, and the second says that no new height functions are obtained for algebraic varieties.  \\

The fact that heights are assigned to all finite rank vector bundles and not just line bundles is a crucial feature; indeed the height function that recovers Malle's conjecture in this framework comes from a vector bundle and not a line bundle. This demonstrates that this level of generality is necessary. \\

We now state the following version of the main conjecture of Ellenberg-Satriano-Zureick-Brown, which recovers (weak forms of) both the Batyrev-Manin conjecture and Malle's conjecture: 

\begin{conjecture}[Stacky Batyrev-Manin-Malle Conjecture of Ellenberg-Satriano-Zureick-Brown]\label{conj:stackybatman}
Let $\fX$ be a ``nice" algebraic stack defined over a number field and let $\E$ be a "nice" vector bundle on $\fX$. Then there is an open dense substack $\mathcal{U}$ of $\fX$ such that for all $\ep>0$ one has

\[\#\{P\in \mathcal{U}(K):H_{\fX,\E}(P)\leq B\}=N_{\mathcal{U},\E,K,\ep}(B)=O_\ep(B^{a(\E)+\ep})\]
where $a(\E)$ is a number depending at most on $\E$.

\end{conjecture}

Here "nice" means it satisfies a Northcott type property that allows one to count points in an open substack. \\

In the framework of \cite{ESZ-B} the Malle and Manin-Batyrev conjectures represent two extremes of their theory of heights. The Manin conjecture involves counting points on a projective variety with respect to a Weil height and no theory of algebraic stacks is required. On the other hand, Malle's conjecture cannot be interpreted as a problem about counting rational points on a scheme. Indeed, as shown in \cite{ESZ-B} Malle's conjecture involves counting rational points on the classifying stack $BG$ where $G$ is a finite group. The theory of algebraic stacks is essential for this interpretation of the Malle conjecture, and the standard theory of heights on projective algebraic varieties is insufficient for this purpose.\\

In this article we study the Ellenberg-Satriano-Zureick-Brown's theory of heights on algebraic stacks in the specific case of \emph{stacky curves}. The family of stacky curves we are interested in lies between the two extremes described above. The stacky curves we consider may be thought of as a smooth projective curve $C$ defined over a number field, along the the data of finitely many points $P$ each endowed with a "stabilizer" group of the form $\bZ/m_P\bZ$ with $m_P>1$. The curve $C$ is called the \emph{coarse space} of the stacky curve, and is the "best approximation" of the stacky curve by an algebraic variety. We think of a stacky curve as being built out of the coarse space by adding multiplicities to finitely many points.\\

We now give a description of our main results.  First, we give a ground-up construction of a type of height function on stacky curves, that recovers the E-S-ZB anti-canonical height on a stacky curve with coarse space $\bP^1$. This construction is detailed in Section \ref{sec:MCurves}. We denote this height function by $\HH_{-K_{\fX}}$. Due to the technical nature of this aspect of our work, we defer stating the results until the next section. \\ 

The theory developed in \cite{ESZ-B} leads to two main conjectures, that depend on the positivity properties of a certain function $\textnormal{edd}\colon \fX(K)\rightarrow \mathbb{R}$. In the context of stacky curves \cite[Proposition 4.9]{ESZ-B} says that $\textnormal{edd}(x)=\log H_{-K_\fX}(x)$. In other words is the same as the logarithmic anti-canonical height. When $\textnormal{edd}(x)$ is \emph{positive} we think of $\fX$ as being "Fanoish" (see \cite[Page 39]{ESZ-B} for a provisional definition of this term) and \cite[Main Conjecture]{ESZ-B} should apply to the anti-canonical height. For stacky curves positivity is equivalent to the anti-canonical height having the Northcott property;there are finitely many points of bounded height. One of our main motivations for this paper is to answer a question of Ellenberg, Satriano, and Zureick-Brown about a stacky curve with "positive" anti-canonical height. Indeed, they asked whether one can obtain Conjecture \ref{conj:stackybatman} for the stacky curve, which has coarse space $\bP^1$ and three half-points. Our first theorem answers their question in the affirmative: 

\begin{theorem} \label{MT3}
Let $\fX$ be the stacky curve obtained by adding three half points to $\bP^1$. Then the Stacky Batyrev-Manin-Malle Conjecture is true for $H_{-K_\fX}$. 
\end{theorem}

In fact we prove more in Theorem \ref{MT3}, namely we show that the points of bounded E-S-ZB anti-canonical height satisfies an exact order of magnitude. \\

 We note that P.~Le Boudec has independently verified this case in a private communication. \\

A natural question which arises is when in the anti-canonical height "positive" in the sense described above. Equivalently, when does \cite[Main Conjecture]{ESZ-B} apply to the anti-canonical height of a stacky curve. Precisely, when does the  E-S-ZB anti-canonical height function produce a genuine height function, in the sense that it satisfies \emph{Northcott's property}. We show that the positivity of the anti-canonical height is equivalent to the positivity of the \emph{Euler characteristic} of the stacky curve. We can interpret our family of stacky curves as a mixture of geometry (the coarse space $\bP^1$) and arithmetic (the addition of various stacky points). In the purely geometric case of algebraic curves, the geometry determines the arithmetic in the sense that the anti-canonical height has Northcott's property if and only if the curve has positive \emph{Euler characteristic}. It turns out in the stacky case the Euler characteristic measures not only the geometry of the coarse space (base scheme) but also the additional structure imposed by the stacky points. This leads to the following theorem, which is an exact analogue of the above statement regarding algebraic curves: 

\begin{theorem} \label{intro:Northfail} Let $\fX$ be a proper smooth stacky curve defined over $\bQ$ that has coarse space $\bP^1$ or is isomorphic to a smooth projective curve. Then the anti-canonical height $\HH_{-K_\fX}$  has the Northcott property if and only if $\chi(\fX)>0$, where $\chi(\fX)$ is the Euler characteristic of $\fX$. 
\end{theorem}


The above result tells us that the E-S-ZB height machine when applied to the tangent bundle recovers the behavior of the Weil height machine when applied to the tangent bundle, providing additional evidence that the E-S-ZB theory of heights if the correct generalization of the classical theory. In particular, Theorem \ref{intro:Northfail} demonstrates that Conjecture \ref{conj:stackybatman} is a direct general ization of the classical Batyrev-Manin conjecture for Fano varieties. \\

When $\textnormal{edd}(x)$ is not positive, Ellenberg, Satriano, and Zureick-Brown have proposed a generalized \emph{Vojta's conjecture} applicable to the case of algebraic stacks \cite[Conjecture 4.23]{ESZ-B}. In the case of stacky curves, the stacky Vojta conjecture can be phrased in a particularly simple manner. The failure of $\textnormal{edd}(x)$ to be positive is equivalent to the failure of $\HH_{-K_X}$ to have the Northcott property. Given that the anti-canonical height $\HH_{-K_X}$  fails to have the Northcott property when the Euler characteristic is negative, a natural question is how far the function is from having Northcott's property. P. Vojta asked similar questions in his thesis (and see \cite{Voj} for an update) and indeed gave a far reaching series of conjectures relating the geometry and arithmetic of algebraic varieties. \\

In \cite[4.7]{ESZ-B} it is speculated that \cite[Conjecture 4.23]{ESZ-B} for stacky curves should follow from some version of the $abc$-conjecture. In the case of algebraic curves, Vojta's conjecture is known to be equivalent to the $abc$-conjecture.  We show that, much like the case of algebraic curves, the stacky analogue of Vojta's conjecture in the curve case is equivalent to the $abc$-conjecture. We formulate this as follows:

\begin{theorem} \label{MT4} Let $\fX$ be a proper smooth stacky curve defined over $\bQ$ that has coarse space $\bP^1$ or is isomorphic to a smooth projective curve. Further suppose that $\fX$ has negative Euler characteristic. Then the following statements are equivalent:
\begin{enumerate}
    \item The $abc$-conjecture holds; and
    \item For all $\fX$ satisfying the hypotheses of the theorem and for all $\delta > 0$ the function $\HH_{-K_{\fX}} \cdot H^\delta$ has Northcott's property, where $H([x,y]) = \max\{|x|, |y|\}$ is the usual height function on $\bP^1(\bQ)$.
\end{enumerate}

\end{theorem}

Theorem \ref{MT4} shows that Conjecture 4.23 in \cite{ESZ-B} is equivalent to the $abc$-conjecture, answering a question of Ellenberg, Satriano, and Zureick-Brown. \\

In \cite{ESZ-B} the authors wonder if the stacky Vojta conjecture might be more "in reach" for algebraic stacks obtained by rooting along a divisor $D$ on an scheme $X$. The proof of Theorem \ref{MT4} shows that if there is some $m\geq 4$ such that item (2) in Theorem \ref{MT4} holds for $\fX_m = \fX(\bP^1 : ((0,1,\infty) : (m,m,m))$ then a weak variant of the $abc$-conjecture can be derived. Specifically, there exists a positive number $c_m \geq 1$ such that for any co-prime $a,b,c \in \bZ$ with $a + b = c$ and $\ep > 0$ that 
\[\max\{|a|,|b|,|c|\} \ll_{\ep, m} \rad(abc)^{c_m + \ep}.\]
In particular, any progress on the stacky Vojta conjecture for curves would lead to substantial progress on the $abc$-conjecture. \\

\section{A further elaboration of our ideas} 
\label{intro2}

In this section we motivate and describe our main results in more detail, as well as describe our grounds-up height construction.

\subsection{An elementary height machine on stacky curves}

Our point of view of algebraic stacks is to adopt a bottom-up perspective. In other words, we define our algebraic stacks in terms of a base variety along with some extra data which is enough to uniquely construct an algebraic stack. As we are interested in a well behaved family of stacky curves this description will be particularly simple. The bottom up point of view allows us to discuss the objects we are interested in a concrete way that avoids technicalities and emphasizes the data most important for our purposes. The interested reader may consult \cite{BottomUp} for general results involving the bottom up perspective on algebraic stacks and \cite[Lemma 5.3.10]{CanRing} for the case of stacky curves.\\

Using the bottom up perspective, a stacky curve defined over a number field $K$ is \emph{uniquely} determined by the following data: A smooth $X$ defined over $K$. A finite number of \emph{stacky} points $P_1,\dots ,P_r$ along with integer multiplicities $m_{P_i}=m_i>1$ attached to each point $P_i$.  
 We use the notation.
	\[\fX=(X:(P_1,m_1), \cdots, (P_r,m_r))\] to denote the stacky curve with multiplicities $m_{P_i}=m_i$ at the points $P_i$. We will write $\fX(X:(\Ba, \Bm))$ as an abbreviation. We identify the rational points of the stack $\fX$ with those of the coarse space. That is, we require
	\[\fX(K)=X(K).\] 
	
In \cite{ESZ-B} some care is taken to work with the locus of stacky points on an algebraic stack. In particular, one must contend with the accumulation of \emph{infinitely new} stacky points. We ignore such difficulties since in the context of Conjecture \ref{conj:stackybatman} one may always remove such an accumulating collection of stacky points. Therefore from the point of view of counting points asymptotically on stacky curves this presents no issues.\\

We now will work exclusively with (stacky) curves of this shape, which is the algebraic stack obtained by replacing $\{a_1, \cdots, a_n\}\subseteq {\bP^1}$ with $B(\bZ/m_j \bZ)$ for $j = 1, \cdots, n$. For a precise definition of this object, see \cite{BottomUp}. In \cite{ESZ-B} these are treated in section 4.6.  The points $a_1, \cdots, a_n$ are "stacky points" which have an attached stabilizer group $\mathbb{Z}/m_j \mathbb{Z}$ and account for $\fX(\mathbb{P}^1_\bQ; (\Ba, \Bm))$ not being a projective variety. The other points behave like the points on the quasi-projective variety $\bP^1-\{a_1, \cdots, a_n\}$.  This construction is an example of a ``bottom-up" description of an algebraic stack; the algebraic stack $\fX(\mathbb{P}^1_\bQ;(a_1,m_1), \cdots, (a_n,m_n))$ is constructed from the data of $\mathbb{P}^1_\bQ$ and the points with their associated multiplicities. In fact a large class of algebraic stacks can be constructed this manner, see \cite{BottomUp} for details and for further references and recent appearances of these objects see \cite{CanRing}, \cite{LGStackyCurve},\cite{StackyCurves}, \cite[4.6]{ESZ-B}. \\

To apply the E-S-ZB height machine to a stacky curve $\fX(X:(\Ba,\Bm))$ we must choose a vector bundle on $\fX(X:(\Ba,\Bm))$. We will be primarily interested in the stacky curves $\fX(\bP^1_\bQ: (a_1,m_1),\dots,(a_r,m_r))$ and the vector bundle a line bundle. As in the classical case one may work with divisors rather then line bundles. It suffices to consider divisors of the form

\[\mathcal{D}=D+\sum_{i=1}^sc_i [a_i]\]

where $D$ is a divisor on $\bP^1$ and $0\leq c_i\leq m_i-1$. To associate a height on such a divisor, we associate a height to each $b_i[a_i]$ and extend linearly. Motivated by \cite{ESZ-B} we have the following construction. For each stacky point $a_i = [\alpha_i : \beta_i]$ with $\alpha_i,\beta_i$ coprime integers we associate to it the linear form $\ell_i(x,y) = \alpha_i y- \beta_i x$. For each $m_i$ define $\phi_{m_i}(n)$ is defined to be the smallest positive integer such that $n \phi_{m_i}(n)$ is a perfect $m_i$-th power. The height function associated to $c_i[a_i]$ is then

\[H_{c_i[a_i]}(x,y)=\phi_{m_i}(\ell_i(x,y)^{c_i})^{\frac{1}{m_i}}.\]

The linear form $\ell_i$ takes into account the point $a_i$ and $\phi_{m_i}$ accounts for the multiplicity of $a_i$, while the power $c_i$ accounts for the multiple of $[a_i]$. The introduction of the functions $\phi_{m_i}$ is due to \cite{ESZ-B} and working with these functions is a key feature of stacky curves with coarse space $\bP^1$. As any divisor $D$ on $\bP^1$ is linearly equivalent to $\dfrac{\deg D}{2}(-K_{\bP^1})$ we can define a height function for any divisor $\mathcal{D}=D+\sum_{i=1}^sc_i [a_i]$ on $\fX$ as

\begin{equation}\label{intro:Stacky Curve Height Machine}
H_{\mathcal{D}}(x,y)=\max\{\vert x\vert,\vert y\vert\}^{\mathcal{D}}\cdot \prod_{i=1}^r\phi_{m_i}(\ell(x,y)^{c_i})^{\frac{1}{m_i}}
\end{equation}

whenever $x,y$ are coprime integers. We call the above construction the \emph{Stacky Curve Height Machine}. The benefit of this construction is that it is completely explicit and requires no machinery from the theory of algebraic stacks to compute. For example, in \cite{ESZ-B} one computes the height using a stacky version of Arakeolv theory, and an object known as a tuning stack. Such difficulties can be avoided in our context.  In fact, this construction detects the \emph{$S$-integral points} of the stacky curve $\fX$ in a suitable sense. In particular a point $(x:y)$ is an integral point of $\fX(\bP^1_\bQ,(a_1,m_1),\dots,(a_r,m_r))$ if and only if

\[H_{[a_i]}=\phi_{m_i}(\ell_i(x,y))^{\frac{1}{m_i}}=1\textnormal{ for }1\leq i\leq r.\]

If we take $\mathcal{D}=K_{-\fX}$ then we obtain

\begin{equation}\label{Hdefgen}
H_{-K_\fX}(x,y)=\max\{\vert x\vert,\vert y\vert \}^{\chi(\fX)}\prod_{i=1}^r\phi_{m_i}(\ell_i(x,y))^{\frac{1}{m_i}}
\end{equation}
where

\begin{equation} \label{EulerChar} \chi(\fX)=\deg (-K_\fX)=\delta(\Bm) = 2 - \sum_{i=1}^n \left(1 - \frac{1}{m_i}\right),
\end{equation}

is the Euler characteristic of the stacky curve, which is defined as the degree of the anti-canonical divisor.

When we wish to emphasize the dependence of $(\Ba,\Bm)$ we write $H_{-K_\fX}([x:y])=\HH_{(\Ba, \Bm)}(x,y)$. 

\subsection{Properties of the anti-canonical E-S-ZB height $H_{-K_{\fX}}$}

The Northcott property of the naive height on $\bP^1$ implies that the E-S-ZB anti-canonical height $H_{-K_\fX}$ has the Northcott property whenever $\chi(\fX)=\delta(\Bm)>0$. On the other hand, if $\chi(\fX)\leq 0$ then it is not at all obvious whether $H_{-K_\fX}$ should have the Northcott property. The following question is fundamental: Let $\L$ be a line bundle on $\fX(X: (\Ba,\Bm))$ and let $H_{\L}$ be the associated E-S-ZB height. \emph{When does $H_\L$ have the Northcott property?} We will tackle this question when $\L=-K_\fX$ and $X=\bP^1$ leaving the general case for future study. In this situation, this amounts to studying the positivty properties of the function $\textnormal{edd}(x)$ defined in \cite{ESZ-B}. This is important for the following reason, when $\textnormal{edd}(x)$ is sufficiently positive then \cite[Main Conjecture]{ESZ-B} should apply. While when $\textnormal{edd}(x)$ is negative the stacky Vojta conjecture of \cite{ESZ-B} should be applicable.

\begin{theorem} \label{Northfail} Let $\{a_1, \cdots, a_n\} \subset \bP_\bQ^1$ and $\Bm = (m_1, \cdots, m_n)$ be a vector of multiplicities. Then whenever 
\[\chi(\fX(\bP^1 : (\Ba, \Bm)) = \delta(\Bm) \leq 0,\] 
the anti-canonical height $H_{-K_\fX}$ given by (\ref{Hdefgen}) on the curve $\fX(\bP^1 : (a_1, m_1), \cdots, (a_n, m_n))$ fails to have the Northcott property. 
\end{theorem}

This may seem surprising for the following reason: in the case when $\Bm = \underbrace{(2, \cdots, 2)}_n$ and $n \geq 5$, the height (\ref{Hdefgen}) should have the Northcott property, based on the heuristic that although the exponent of the naive height is negative, that the complexity of the equations should compensate for this. For example consider the following height

\[\HH([x:y])^2=\dfrac{\sqf(\mid x\mid)\sqf(\mid y\mid)\sqf(\mid x+y\mid)\sqf\mid x-y\mid)\sqf(\mid 2x+1 \mid )}{\max(\mid x\mid,\mid y\mid)} \]
where $x,y$ are coprime integers and $\sqf(\mid x\mid)$ is the square free part of $x$. The Northcott property asks if there are infinitely many $[x:y]\in \bP^1(\bZ)$ with $\HH([x:y])^2\leq T$ for some number $T$. A priori it may seem difficult to find infinitely many coprime pairs of integers $(x,y)$ such that $\sqf(\mid x\mid)\sqf(\mid y\mid)\sqf(\mid x+y\mid)\sqf\mid x-y\mid)\sqf(\mid 2x+1 \mid )$ is small compared to $\max(\mid x\mid,\mid y\mid)$.   \\

If we assume that the E-S-ZB theory should behave roughly like its classical counterpart we can argue the converse: when $\chi(\fX)\leq 0$ one should have that $H_{-K_\fX}$ should fail to have the Northcott property.\\

\emph{Arithmetic intuition from the theory of curves:}
For an algebraic variety $X$ defined over a number field $K$, we say that a function $H_X$ has the \emph{strong Northcott property} if for any positive numbers $T, D$ the set
\begin{equation} \N(X, K; D, T) = \{P \in X(\ol{\bQ}) : H_X(P) \leq T, [K(P) : K] \leq D\} 
\end{equation}
is finite. Here $K(P)$ is the finite extension of $K$ obtained by adjoining the coordinates of $P$. We introduce some new notions to make our statements as precise as possible.

\begin{definition}[\cite{arr}]
Let $C$ be a smooth and geometrically integral projective curve defined over a number field $K$. Let $e\geq 1$ be an integer. Define

\[C_e=\{P\in C(\bar{K})\colon [K(P):K]\leq e\}.\]
The \emph{arithmetic degree of irrationality} of $C$ is defined to be

\[\textnormal{a.rr}_K(C)=\min\{e\in \bZ_{\geq 1}\colon C_e \textnormal{ is infinite}.\}\]
\end{definition}
The arithmetic degree of irrationality is always finite since the gonality of a curve provides an upper bound (see \cite{arr} for an introduction and further references). Indeed, if $C$ is a smooth projective and geometrically integral curve defined over a number field $K$ with gonality $d$ then there is a degree $d$ finite surjective morphism $\phi\colon C
\ra \bP^1_K$, which provides infinitely many points on $C$ of degree at most $d$. The notion of arithmetic irregularity can be related to the failure of the Northcott property for the anti-canonical height.

\begin{proposition}\label{prop:arrprop}
Let $C$ be a smooth projective and geometrically integral curve defined over a number field $K$ with $\chi(C)\leq 0$. Let $-h_C$ be the logarithmic anti-canonical height of $C$. Then \begin{align*}
    \textnormal{a.rr}_K(C)&=\min\{D\in \bZ_{\geq 1}\colon \#\N(C, K; D, T)=\infty\textnormal{ for some real number }T\}\\ &=\min \{D\in \bZ_{\geq 1}\colon \#\{P\in C(\ol{K}) : h_C(P) \leq T, [K(P) : K] \leq D\}=\infty\textnormal{ for some real number }T\}
\end{align*}
\end{proposition}
\begin{proof}
First suppose that $\chi(C)=0$. Then $g(C)=1$ and $C$ is a genus 1 curve. Therefore the canonical divisor is trivial. Then we have that $\textnormal{a.rr}_K(C)\in \{1,2\}$. If $\textnormal{a.rr}_K(C)=1$ there $C$ has infinitely many $K$ points and thus $-M_1\leq -h_C(P)\leq M_1$ for some constant $M_1$ and so $\min\{D\in \bZ_{\geq 1}\colon \#\N(C, K; D, M_1)=\infty\}=1.$ Conversely if $\min\{D\in \bZ_{\geq 1}\colon \#\N(C, K; D, M_1)=\infty\}=1$ then there are infinitely $K$ points in $C$ and we have $\textnormal{a.rr}_K(C)=1$. If $\textnormal{a.rr}_K(C)=2$ then $C(K)$ is finite. Thus $\min\{D\in \bZ_{\geq 1}\colon \#\N(C, K; D, M_1)=\infty\}=2$. Conversely if $\min\{D\in \bZ_{\geq 1}\colon \#\N(C, K; D, M_1)=\infty\}=2$ then $C(K)$ is finite. We may now assume that $\chi(C)<0$ which means $g(C)>1$. Therefore $C$ is of general type, and $h_C$ is a height associated to the canonical divisor. As $g>1$ the canonical divisor is ample, therefore there is some constant $T$ with $h_C(P)\geq T$ for all $P\in C(\bar{K})$. Therefore, $-h_C(P)\leq -T$ for all $P\in C(\bar{K})$. Thus 
\[\#\N(C, K; D, -T)=\infty\text{ if and only if} \#C_D=\infty\] and we obtain the desired inequality.
\end{proof}

The anti-canonical height of an algebraic curve of non-positive Euler characteristic defined over a number field \emph{always fails to have the strong Northcott property}. In the case of stacky curves of the form $\fX(\bP^1: (\Ba,\Bm))$ we have that the arithmetic irregularity and gonality are both equal to one. It is clear that there are infinitely many rational points and the coarse space mapping gives a cover $\fX\ra \bP^1$. Therefore if the E-S-ZB theory behaves analogously to the theory of Weil heights we expect that the E-S-ZB anti-canonical height on stacky curves with non-positive Euler characteristic to also fail to have Northcott's property.\\

\emph{Geometric intuition from the interplay of positivity properties and the theory of heights:} Let $X$ be a smooth projective variety defined over a number field $K$. If $\L$ is a line bundle on $X$ there is a height function $H_\L$ associated to $\L$. If $\L$ is ample then $H_\L$ has the Northcott property, while if $\L$ is not ample we expect $H_\L$ to fail to have the Northcott property. Thus for a stacky curve $\fX$ to determine when $H_\L$ has the Northcott property we should determine when $\L$ is an "ample line bundle". There does not seem to be a general theory of ample line bundles on algebraic stacks, however for stacky curves there is a natural notion extending the classical definition. Let $\fX$ be a nice stacky curve. Then there is a natural notion of degree for line bundles, extending the usual notion of the degree of a line bundle on a curve. It is natural the to call a line bundle $\L$ on $
\fX$ ample if $\deg_\fX \L>0$. Thus one might expect that $H_\L$ has the Northcott property if and only if $\deg_\fX \L>0$. Now take $\L=\T_\fX$. Then one has  
\[\deg_\fX \T_\fX=\chi(\fX)\]
and consequently we should expect that the anti-canonical height $H_{\T_\fX}$ has the Northcott property if and only if $\deg_\fX\T_\fX=\chi(\fX)>0$.\\

In summary, Theorem \ref{Northfail} shows that arithmetic and geometric intuition described above prove to be correct when $\fX$ is birational to $\bP^1$. In other words, the E-S-ZB theory behaves as predicted by the classical situation of curves and the classical theory of heights associated to ample line bundles. This answers a question posed by Ellenberg. It also suggests that the theory of heights may be of help in developing a theory of positivity in the setting of algebraic stacks. Indeed, one potential way of proceeding in certain cases is the following. It is relatively easy to show that for any line bundle $\L$ on a stacky curve $\fX=\fX(X: (\Ba,\Bm))$ one can find an integer $N_\fX$ such that $N_\fX\L=\pi_\fX^*L$ for some line bundle $L$ on $X$, here $\pi_\fX$ is the coarse space mapping. Since there is a good theory of ampleness on $X$ one might try and say that $\L$ is ample if and only if $N_\fX\L$ is ample. This is equivalent to demanding that $\deg_\fX \L>0$ as \[\deg_\fX \L>0\iff \deg_\fX N_\fX\L>0\iff \deg_\fX \pi_\fX^*L>0\iff \deg_X L>0. \]

As the second notion may be extended to certain other classes of algebraic stacks, our results suggests that this may be an appropriate generalization of ampleness in arithmetic settings. On the other hand, this notion of ampleness only depends on coarse space of the curve, which makes it potentially unsatisfactory from a stacky point of view. Getting to the bottom of this tension seems to very interesting and deserves further study.\\

The proof of Theorem \ref{Northfail} uses the following theorem about elliptic curves:

\begin{theorem} \label{g1 thm} Let $F \in \bZ[x,y]$ be a  non-singular binary quartic form. Then there exists square-free $d \in \bZ$ such that the the curve
\[dz^2 = F(u,v)\]
has a rational point and such that its Jacobian has positive rank as an elliptic curve defined over $\bQ$.
\end{theorem}

The proof of Theorem \ref{g1 thm} is provided to us by A.~Shnidman in \cite{ShnidMO}, and we graciously acknowledge his assistance. \\ 

Combining these results gives the following uniform statement.

\begin{corollary} \label{northfailcor} 
Let $\fX$ be a smooth proper stacky curve defined over $\bQ$ such that $\fX$ has coarse space $\bP^1_\bQ$ or $\fX$ is a projective algebraic curve. Let $H_{\fX}$ be the height associated to the anti-canonical divisor $-K_\fX$. Then $\chi(\fX)>0$ if and only if $H_\fX$ has the strong Northcott property.

\end{corollary}
\begin{proof}
If $\fX$ is an algebraic stack and not an algebraic curve then it is of the form $\fX=\fX(\bP^1_\bQ;(\Ba,\Bm))$ and (\ref{Northfail}) gives the desired result. On the other hand, if $\fX$ is a smooth projective and geometrically integral curve then $\chi(\fX)\leq 0$ implies that $-H_{C}$ does not have the Northcott property by the fact that $\textnormal{a.rr}_K(C)<\infty$.
\end{proof}

Corollary \ref{northfailcor} compels the following question: how should one define the arithmetic irregularity of a stacky curve $\fX=\fX(X:(\Ba,\Bm))$? There are two natural definitions that come to mind when $\chi(
\fX)\leq 0$. We might define $\textnormal{a.rr}_K(\fX(X:(\Ba,\Bm)))=\textnormal{a.rr}_K(X)$. That is, the arithmetic irregularity of the stack is the arithmetic irregularity of the underlying base curve/coarse space. On the other hand using one might try to use (\ref{prop:arrprop}) to define the arithmetic irregularity of $\fX$ when $\chi(\fX)<0$ as \[\min\{D\in \bZ_{\geq 1}\colon \#\{P\in X(\bar{\bQ}):H_{T_\fX}(P)\leq T,[\bQ(P):\bQ]\leq D,\textnormal{ for some real number }T \}=\infty\}.\]

Our results show that these definitions are the same when $\chi(\fX)\leq 0$, and $\fX=\fX(\bP^1_\bQ,(\Ba,\Bm)),$ precisely as what occurs in the case of algebraic curves.  \\

In cases where we can prove that the Northcott property fails, according to \cite{ESZ-B} there should be a stacky Vojta conjecture. This can be explained in the following elementary manner. We would like to know whether the height $\HH_{(\Ba, \Bm)}$ can be modified to recover the Northcott property. This has the following motivation from the classical setting. Let $L$ be a line bundle on a smooth projective variety $X$. Let $M$ be a chosen ample line bundle on $X$. Then one considers

\[\inf\{t\in \bR_{\geq 0}\colon L+tM\textnormal{ is ample}\}\]
as a measure of how far $L$ is being from an ample line bundle. We would like to ask similar questions for our height functions. Difficulties arise because the E-S-ZB height machine is not functorial in the usual sense; in the setting of algebraic varieties one can work with a linear spaces of divisors and then apply the height machine which by functoriality will respect the linear structure. Such methods are not immediately available to us. Instead we will apply the height machine, and then apply linear operations. Our motivating question is as follows, if $\fX=\fX(\bP^1:(\Ba,\Bm))$ with $\chi(\fX)\leq 0$, what can be said about

\[\inf\{t\in \bR_{\geq 0}\colon H_{\bP^1}^tH_{-K_\fX}\textnormal{ has the Northcott property}\}.\]

Clearly, if we change the exponent in the classical part of the height so that it is positive, then we will recover the Northcott property. In fact we expect that something far less drastic suffices. \\

For a real number $\delta$ and the curve $\fX(\bP^1 : (\Ba, \Bm))$, define the height
\begin{equation} \label{ESZBmod} \HH_{(\Ba, \Bm)}^\delta(x,y) = \prod_{i=1}^n \phi_{m_i}(\ell_i(x,y))^{1/m_i} \max\{|x|, |y|\}^{\delta}.
\end{equation}
We then see that $\HH_{(\Ba, \Bm)} = \HH_{(\Ba, \Bm)}^{\delta(\Bm)}$. Next put
\begin{equation} \gamma(\fX) = \inf \{\delta \in \bR: \HH_{(\Ba, \Bm)}^\delta \text{ has the Northcott property on }  \fX \}.
\end{equation}
In fact $\gamma(\fX)$ depends only on $\Bm$, so we may also write it as $\gamma(\Bm)$. We make the following conjecture: 

\begin{conjecture}[Northcott conjecture for stacky curves with coarse space $\bP^1$]\label{conj:NorthcottConj} For all $(\Ba, \Bm)$, we have $\gamma(\fX) = \min\{\delta(\Bm), 0\}$.
\end{conjecture}

Conjecture \ref{conj:NorthcottConj} is in fact a version of Vojta's conjecture for stacky curves, and agrees with a conjecture of Ellenberg, Satriano, and Zureick-Brown in \cite{ESZ-B}. Towards this conjecture, we have the following:

\begin{theorem} \label{Northsuccess} We have $\gamma(\fX) = 0$ if $\chi(\fX) = \delta(\Bm) \geq 0$. Moreover, the height $\HH_{(\Ba, \Bm)}^0$ has the Northcott property if and only if $\delta(\Bm) < 0$.
\end{theorem}

Combined with Theorem \ref{Northfail} the conjecture predicts that the set of $\delta\in \bR$ such that $\HH^\delta_{\Ba,\Bm}$ has the Northcott property is an interval of the form $(\delta(\Bm),\infty)$ when $\delta(\Bm)<0$ and $(0,\infty)$ when $\delta(\Bm)\geq 0$. Therefore, while Theorem \ref{Northfail} tells us we cannot count points with $\HH_{\Ba,\Bm}$, Conjecture \ref{conj:NorthcottConj} predicts that we can count points using $\HH^{\delta(\Bm)+\ep}_{\Ba,\Bm}$ for any $\ep >0$.  \\

Turning back to our motivational question in this area, we find that \ref{conj:NorthcottConj} predicts that \[\inf\{t\in \bR_{\geq 0}\colon H_{\bP^1}^tH_{\T_\fX}\textnormal{ has the Northcott property}\}=0\] 
when $\chi(\fX)\leq 0$. Furthermore, the infimum is a limit point that is not achieved. From this perspective it predicts that the anti-canonical height $H_{-K_\fX}$ lies on the boundary of those line bundles with the Northcott property. Furthermore Conjecture \ref{conj:NorthcottConj} suggests that the Euler characteristic of the curve can be recovered in a natural way from the E-S-ZB theory of heights, as the smallest possible exponent for the \emph{classical part of the height} such that the resulting function fails to have the Northcott property. \\

We proceed to prove that Conjecture \ref{conj:NorthcottConj} is a consequence of the $abc$-conjecture. However, it seems that we are very far from being able to prove such a result as strong as Conjecture \ref{conj:NorthcottConj} unconditionally. 

\begin{theorem} \label{nearnorth} Suppose that the $abc$-conjecture holds. Then for any $\delta > \delta(\Bm)$ the function $\HH_{(\Ba, \Bm)}^{\delta} (x,y)$ on $\fX(\bP^1: (\Ba, \Bm))$ has Northcott's property. 
\end{theorem}

In fact Conjecture \ref{conj:NorthcottConj} is \emph{equivalent} to the $abc$-conjecture; see Theorem \ref{MT4}. The proof of the converse is quite different and so we give it in a separate subsection. \\

\subsection{Quantitative arithmetic on stacky curves} 

In the positive Euler characteristic case, we consider a particular family of stacky curves, which includes an important example suggested by J.~Ellenberg\footnote[1]{"What?s up in arithmetic statistics?" Number Theory Web Seminar, July 23, 2020} and show that our theory of heights matches \cite{ESZ-B} in this instance. Finally, we verify a specific instance of the main conjecutre in \cite{ESZ-B} given by Ellenberg-Satriano-Zureick-Brown \footnotemark[1]  using analytical methods. We remark that P.~Le Boudec had obtained the same result as us in independent work (private communication). \\ 

We study the expression (\ref{EulerChar}) a bit more carefully. It is easy to deduce that $\delta(\Bm) \geq 0$ if and only if $n \leq 4$, and $\delta(\Bm) > 0$ only if $n \leq 3$. We will not consider the case $n \leq 2$ in this paper. \\ 

If we assume $m_1 \leq m_2 \leq m_3$ then the only cases when we have positive Euler characteristic are when $m_1 = m_2 = 2$, $m_1 = 2, m_2 = m_3 = 3$ or $m_1 = 2, m_2 = 3, m_3 = 4$. In each of these three cases the Northcott property for $\HH_{(\Ba, \Bm)}$ holds trivially. \\

We now focus on the simplest cases, where $m_1 = m_2 = 2$ and $m_3 = m, m \geq 2$. Using that $\PGL_2$ acts 3-transitively on $\bP^1$, we reduce to the case $\{a_1, a_2, a_3\} = \{0,-1,\infty\}$. For $[x,y] \in \bP^1$ we may then set
\[x = x_1 x_2^2, y = y_1 y_2^2\]
with $x_1, y_1$ square-free, and
\[x + y = z_1 z_2^2 \cdots z_{m-1}^{m-1} z_m^m,\]
with $z_1, \cdots, z_{m-1}$ square-free. In this notation, the E-S-ZB height is given by 
\begin{equation}\label{eq:HeightDef}
    \HH_{\Bm}(x,y) = |x_1|^{1/2} |y_1|^{1/2} |z_1^{m-1} \cdots z_{m-1}|^{1/m} \max\{|x_1 x_2^2|,|y_1 y_2^2|\}^{1/m}.
\end{equation}
We normalize the height so that the exponent of the ``classical part" is equal to one, to obtain the normalized height
\begin{equation} \label{normht} H_m(x,y) = |x_1|^{m/2} |y_1|^{m/2} |z_1^{m-1} \cdots z_{m-1}| \max\{|x_1 x_2^2|, |y_1 y_2^2|\}.
\end{equation}
We now put
\begin{equation} \label{NT} 
N_m(T) = \# \{(x_1, x_2), (y_1, y_2), (z_1, \cdots, z_m) : \gcd(x_1 y_1, x_2 y_2) = 1, x_1, y_1, z_1, \cdots, z_{m-1} \text{ square-free and pairwise co-prime },
\end{equation}
\[x_1 x_2^2 + y_1 y_2^2 = z_1 z_2^2 \cdots z_{m-1}^{m-1} z_m^m, |x_1|^{m/2} |y_1|^{m/2} |z_1^{m-1} \cdots z_{m-1}| \max\{|x_1 x_2^2|, |y_1 y_2^2|\} \leq T\}.\]
We prove the following theorem, which gives a crude upper bound for $N_m(T)$:
\begin{theorem} \label{MT0} Let $\fX = \fX(\bP^1 : (0,2), (\infty, 2), (-1, m))$ and let $H_m$ be the height function on $\fX$ defined by (\ref{normht}). Then for any $\ep > 0$ we have
\[T^{1/m} \ll N_m(T) \ll_\ep T^{2/m + \ep}.\]
\end{theorem}
When $m = 2$ the upper bound of Theorem \ref{MT0} is essentially the trivial bound, but it is non-trivial as soon as $m > 2$. In general we expect the exponent in Theorem \ref{MT0} to be equal to the lower bound. Indeed, this can be verified when $m = 2$. Even more, we can give an exact order of magnitude for $N_2(T)$: 
 
\begin{theorem} \label{MT} There exist positive numbers $c_1, c_2, c_3$ such that 
\[c_1 T^{1/2} (\log T)^3 < N_2(T) < c_2 T^{1/2} (\log T)^3 \]
for all $T > c_3$. 
\end{theorem}

In particular, we confirm the stacky Manin-Batyrev  conjecture \cite[Main conjecture]{ESZ-B} for $\fX(\bP^1_\bQ,(a,2),(b,2),(c,2))$. For this stacky curve, \cite[Main conjecture]{ESZ-B} predicts that $N_2(T) = O_\ep \left(T^{1/2 + \ep}\right)$.\footnotemark[1] Our theorem gives an exact order of magnitude for $N_2(T)$. We remark, once again, that P.~Le Boudec had obtained the same result. Further, our counting arguments are similar to those obtained by Le Boudec in \cite{LB} which studies the equation (\ref{variety}). \\

The other cases with positive Euler characteristic do not yield to the simple analytic counting arguments used to prove Theorem \ref{MT}, though in principle counting rational points by height is a well-posed problem. We plan on returning to this issue in the future. \\ 

We illustrate how the stacky curve height machine (equation \ref{intro:Stacky Curve Height Machine}) allows one to detect integral points on $M$-curves. In this case the standard height is given by $H_s(a,b) = \max\{|a|,|b|\}$ and the stacky height given by (\ref{eq:HeightDef}). They are equal precisely when
\[|\sqf(a) \sqf(b) \sqf(a+b)| = 1, \]
or in the notation of (\ref{variety}), that $|x_1| = |x_2| = |x_3| = 1$. (\ref{variety}) then turns into 
\[\pm y_1^2 \pm y_2^2 = \pm y_3^2, \]
and up to rearranging we are essentially counting points on the conic 
\begin{equation} \label{intcon} y_1^2 + y_2^2 = y_3^2.\end{equation}
Therefore if we denote by $\N(T)$ the number of integral points (in the sense of Definition \ref{DarInt}) on $\bP_{2,2,2}^1$ then:

\begin{corollary} \label{intct}  There exist positive numbers $c_1, c_2, c_3$ such that for all $T > c_3$ we have \[c_1 T^{1/2} < \N(T) <c_2 T^{1/2}.\]
\end{corollary}

The proof is elementary, since the curve can be explicitly parametrized by 
\[y_1 = u^2 - v^2, y_2 = 2uv, y_3 = u^2 + v^2.\]
The condition $\max\{|y_1|, |y_2\} \leq T^{1/2}$ is subsumed by $u^2 + v^2 \leq 4T^{1/2}$ say, so number of possible $u,v$'s is $\asymp T^{1/2}$ as desired. \\

Theorem \ref{MT} and Corollary \ref{intct} imply that asymptotically $0$-percent of the rational points on $\fX(\bP^1 :(0,2),,(-1,2),(\infty,2))(\bQ)$ are integral, in the sense of Darmon (Definition \ref{DarInt}). 

To close off this subsection We note that in \cite{LGStackyCurve} Bhargava and Poonen study situations where the rational and integral points of a stacky curve satisfy the Hasse Principle. Motivated by this work we prove that the integral points on $\fX(\bP^1 : (a_1, 2), \cdots, (a_n, 2))$ \emph{satisfy Hasse's principle}. 

\begin{theorem} Let 
\[\fX = \fX(\bP_\bQ^1 : ([a_1 : -b_1], 2), ([a_2 : -b_2], 2), ([a_3, -b_3] : 2) ).\]
Then $\fX$ has integral points if and only if the ternary quadratic form 
\[Q_\Ba(u,v,w) = \det \begin{pmatrix} u^2 & v^2 & w^2 \\ a_1 & a_2 & a_3 \\ b_1 & b_2 & b_3 \end{pmatrix} \]
defines a conic with a rational point. 
\end{theorem}

\subsection*{Notation} We denote by $d_k(n)$ for the number of ways of writing $n$ as a product of $k$ (not necessarily distinct) positive integers, and write $d(n) = d_2(n)$ for the usual divisor function. We will also use the big-$O$ notation as well as Landau's notation. In particular, we will denote in the subscripts any dependencies; if there are no subscripts, then the implied constants are absolute. 

\section{(Stacky) Heights on $M$-curves}\label{sec:MCurves}

In this section we give an alternative construction of the height functions constructed in \cite{ESZ-B} in a special case: we construct the ESZ-B heights associated to line bundles on stacky curves with coarse space $\bP^1_\bQ$. This construction is elementary as it avoids the language of algebraic stacks.\\ 

For our construction, we eschew the theory of algebraic stacks in favour of Darmon's $M$-curves, which is an essentially equivalent theory that emphasizes the bottom up perspective to algebraic stacks. In other words, we only keep track of the minimum amount of data needed to construct the stack. This is analogous to only keeping track of a particular Weierstrass equation of an elliptic curve.

\begin{definition}[\cite{MCurves}]\label{def:MCurve} Let $K$ be a number field.
 An \textbf{$M$-curve over $K$} consists of the following data: 
 
 \begin{itemize}
     \item A smooth projective curve $X$ defined over a number field $K$, and 
     \item For each $P\in X(K)$ a multiplicity $m_P\in \mathbb{Z}_{\geq 1}\cup\{\infty\}$ with $m_P=1$ for all but finitely many $P$.
     \item The $K$-rational points of $\fX$ are defined to be $X(K)$. 
 \end{itemize} 
 \end{definition}

 Before continuing let us fix some notation. We assume that $\fX=(X:(P_1,m_1), \cdots, (P_r,m_r))$ is an $M$-curve as in Definition \ref{def:MCurve}. We furthermore assume that $1<m_{P_i}<\infty$ for each $i$. 


\subsection{Translation between $M$-curves and stacky curves} We now explain the connection between stacky curves and $M$-curves. Recall that given a nice stacky curve $\fX$ over a number field $K$ there is a morphism $\pi\colon \fX\rightarrow X$ to a smooth projective curve $X$ called the coarse space morphism, which is the universal morphism from $\fX$ to a scheme. Often it is more practical to construct $\fX$ from $X$ by specifying a collection of points $P_1,...,P_r$ in the coarse space $X$ and attaching cyclic stabilizer groups $\mu_{P_i}$ to each $P_i$. This is the \emph{bottom up} approach of constructing an algebraic stack. One can think of this as specifying the ramification data of the coarse space morphism $\pi\colon \fX\rightarrow X$. We often think of these points as ``fractional points", because in the divisor class group of the associated curve we have added the point $\frac{1}{\#\mu_{P_i}}P_i$. In other words, we think of a stacky curve as a smooth curve $X$ with a choice of points $P_1,...,P_r$ with stabilizer groups $\mu_{P_i}$ attached to each $P_i$.  This data defines an $M$-curve $\fX=(X:(P_1,\#\mu_{P_1}),\cdots (P_r,\#\mu_{P_r}))$. Conversely, given an $M$-curve $\fX=(X:(P_1,m_1),\cdots (P_r,m_r))$ with each $1<m_{P_i}<\infty$ we consider the stacky curve with stacky points $P_i$ having the stabilizer group $\mu_{m_{P_i}}$. In this way one may establish a bijection between smooth proper geometrically connected Deligne-Mumford stacks of dimension 1 over $K$ with stacky points defined over $K$ that contain an open dense subscheme and possess a projective coarse moduli space and $M$-curves over $K$ with finite multiplicities.

\subsection{Construction of heights}\label{subsec:HeightsMCurves}

We will give an alternative construction of heights on a stacky curve associated to line bundles. Our construction only depends on the definition of an $M$-curve and basic arithmetic. We then show that our height construction corresponds to the heights associated to line bundles in \cite{ESZ-B} when the coarse space is $\bP^1_\bQ$. As in the classical setting, we will work with height functions up to some bounded function. Line bundles on a stacky curve can be described as follows.

\begin{lemma}[\cite{ToricDM} Section 1.3]\label{lem:bundlerep}
	Let $\fX=\fX(X:(P_1,...,P_r),(m_1,...m_r))$. Let $\O_{X}(P)$ be the line bundle associated to the divisor $P$ on $X$. Then there are line bundles $\L_{P_i}$ on $\fX$ such that
	
	\begin{equation}
		\L_{P_i}^{\otimes m_i}\cong \pi_\fX^*\O_X(P)
	\end{equation}
where $\pi_\fX\colon \fX\ra X$ is the coarse space map. Moreover, we have that any line bundle $\L$ on $\fX$ can be uniquely written as
	
	\begin{equation}\label{eq:linebundledecomp}
		\L\cong \pi_\fX^*M\otimes \prod_{i=1}^r\L_{P_i}^{\otimes d_i}
	\end{equation}
	where $0\leq d_i<m_i$ and $M$ is a line bundle on $X$. 
\end{lemma}

The lemma tells us that every line bundle on $\fX$ decomposes canonically as a line bundle on the coarse space $X$ and some stacky line bundles $\L_{P_i}$ that do not arise as line bundles on $X$. However, $\L_{P_i}^{\otimes m_i}$ does arise from a line bundle on $X$. On the other hand, in \cite[Definition 2.21]{ESZ-B} a height function $H_\mathcal{V}$ decomposes as

\begin{equation}\label{eq:ESZ-Bdecomp}
H_{\mathcal{V}}^{\textnormal{st}}\cdot \prod_\nu\delta_{\mathcal{V};\nu}.
\end{equation}

This decomposition highlights the new features of heights on algebraic stacks. The height $H_{\mathcal{V}}$ is no longer additive in the sense that $H_{\mathcal{V}^{\otimes N}}\neq H_{\mathcal{V}}^N$. Furthermore, the height is no longer stable under field extensions. However, $H_{\mathcal{V}}^{\textnormal{st}}$ is stable under field extensions and additive. On the other hand, $H_{\mathcal{V}}$ does not possess a canonical decomposition into local factors, yet $\prod_\nu\delta_{\mathcal{V};\nu}$ does admit such a decomposition.\\

Given a line bundle $\L$ on $\fX$ we write \[\L\cong \pi_\fX^*M\otimes \prod_{i=1}^r\L_{P_i}^{\otimes d_i}\] using Lemma \ref{lem:bundlerep}. Taking inspiration from the decomposition (\ref{eq:ESZ-Bdecomp}) we will build the stable height $H_{\L}^{\textnormal{st}}$ from  Weil height functions $H_{M}$ and $H_{P_i}$ on $X(K)$. The local factors $\prod_\nu\delta_{\mathcal{V};\nu}$ will depend on the purely stacky data of the pairs $(P_i,m_i)$. We will construct the local factors $\prod_\nu\delta_{\mathcal{V};\nu}$ from heights $H_{\L_{P_i}}$ associated to the stacky line bundles $\L_{P_i}$. \\

When working with a stacky curve $\fX=\fX(\bP^1_\bQ,(\Ba,\Bm))$ this becomes particularly explicit. The stable part of the height $H_{\mathcal{V}}^{\textnormal{st}}$ will be a Weil height function on $\bP^1_\bQ$, in other words a function of the form $H_{\mathcal{V}}^{\textnormal{st}}(x:y)=\max\{\vert x\vert,\vert y\vert\}^d$ when $x,y$ are coprime integers. To describe the local terms we will construct heights $H_{\L_{P_i}^{\otimes d_i}}$ which each possess a decomposition

\[H_{\L_{P_i}^{\otimes d_i}}=\prod_{\nu}\delta_{\fX,P_i,d_i\nu}\]
and recover the local decomposition as

\[\prod_{i=1}^r\prod_\nu \delta_{\fX,P_i,d_i\nu}. \]

We first construct the stable or \emph{classical} part of the height. For this we require the notion of degree of a line bundle on an $M$-curve.

\begin{definition}
	Let $\fX=\fX(X:(P_1,...,P_r),(m_1,...m_r)))$. Let $\mathcal{L}=\pi_\fX^* D\otimes \prod_{i=1}^r\L_{P_i}^{\otimes d_i}$. Then we define
	
	\[\deg_\fX \mathcal{\L}=\deg_X D+\sum_{i=1}^r\dfrac{d_i}{m_i}\]
\end{definition}

We may now define the stable height, which will correspond to the part of the height which can be computed classically. 

\begin{definition}\label{def:stableheight}
	Let $\fX=\fX(X:(P_1,...,P_r),(m_1,...m_r)))$ and let $\L$ be a line bundle on $\fX$ with
	
	\[\L\cong \pi_\fX^*M\otimes \prod_{i=1}^r\L_{P_i}^{\otimes d_i}\]
	
	where $0\leq d_i<m_i$ and $M$ is a line bundle on $X$ with $\pi_
	\fX$ being the coarse space map. We define the \emph{stable height} associated to $\L$ as 
	
	\[H_{\L}^{\textnormal{st}}=H_M\cdot\prod_{i=1}^rH_{P_i}^{\frac{d_i}{m_i}}\]
\end{definition}

The stable height can be defined abstractly for all stacky height functions. Later in \ref{prop:stablecomparison} we prove that this definition matches that in \cite{ESZ-B}. The stable height should be thought of as the part of the height that can be computed using classical height functions. Our heights will decompose as a product

\begin{equation}
H_\L=H_\L^{\textnormal{st}}\cdot \prod_{i=1}^rH_{\fX,\L_{P_i}^{\otimes d_i}}. 	
\end{equation}	

 We now concentrate on constructing the height function $H_{\L_{P_i}^{\otimes d_i}}$ associated to the stacky line bundle $\L_{P_i}^{\otimes d_i}$. \\


Fix an $M$-curve $\fX=(X:(P_1,m_1), \cdots (P_r,m_r))$ defined over a number field $K$. We further choose a finite set of primes $S$ of $\O_K$ containing all the primes of bad reduction and all infinite places of $K$. We further choose a smooth and proper model $\underline{X}$ of $X$ over $\O_{K,S}$. Everything we do is \emph{relative to this choice of model} and the finite set of primes $S$.  \\ 

To define meaningful heights we require the basics of arithmetic intersection theory as described by Darmon. It is unsurprising that intersection arises in the construction of heights, as this already occurs in the classical setting via Arakelov's construction of heights through arithmetic intersection theory.  
\begin{definition}[Darmon, \cite{MCurves}] 
Fix an $M$-curve $\fX=(X:(P_1,m_1), \cdots (P_r,m_r))$ defined over a number field $K$. Choose a finite set of primes $S$ of $\O_K$ containing all the primes of bad reduction and all infinite places of $K$. We further choose a smooth and proper model $\underline{X}$ of $X$ over $\O_{K,S}$. Let $P,Q$ be distinct points in $X(K)$ and place $\nu$ a place in $K$ with $\nu\notin S$. Take $\fp_\nu \subset \O_K$ to be the prime ideal associated to $\nu$. We define the \textbf{intersection multiplicity} of $P$ and $Q$ at $\nu$ as follows.
	
	\[(P\cdot Q)_\nu:=\max\{m:\textnormal{ the images of }P,Q
	\textnormal{ in }\underline{X}(\O_{K,S}/\fp_\nu^m)\textnormal{ are equal.}\}\]
where the maximum over the empty set is defined to be 0.	 
\end{definition}


We now package all the intersection multiplicities together while taking into account the arithmetic of the field extension $K\mid \bQ$.


\begin{definition}\label{def:LambdaFactor}
Fix an $M$-curve $\fX=(X:(P_1,m_1), \cdots (P_r,m_r))$ defined over a number field $K$. Choose a finite set of primes $S$ of $\O_K$ containing all the primes of bad reduction and all infinite places of $K$. We further choose a smooth and proper model $\underline{X}$ of $X$ over $\O_{K,S}$. Given a prime $\fp_\nu\subseteq \O_{K}$ we let \[\mathfrak{f}_\nu=[\O_K/\fp_\nu\colon \bZ/(\fp_\nu\cap \bZ)].\]

Now fix $P\in X(K)$,$t\in X(K)\setminus\{P\}$ and $\nu\notin S$.  We define

\begin{equation}
(t\cdot P)_p=\sum_{\nu\notin S,\nu\mid p}\mathfrak{f}_\nu\cdot (t\cdot P)_\nu.
\end{equation}	

Now define

\begin{equation}\label{eq:LocalLambdaFactor}
	\lambda_{S,\underline{X},\nu}(P,t)=\lambda_{\nu}(P,t)=\textnormal{N}(\fp_\nu)^{(t\cdot P)_\nu}.    
\end{equation} 

We then package these together as 

\begin{equation}\label{eq:LambdaFactor}
	\lambda_{S,\ul{X}}(P,t) = \lambda(P,t)=\prod_{\nu\notin S}\lambda_\nu(P,t)=\prod_pp^{(t\cdot P)_p}.
\end{equation}
\end{definition}

The integer $\lambda(P,t)$ is an exponential version of the familiar looking intersection product

\begin{align}
\sum_{\nu\notin S}(t\cdot P)_\nu \log(p)=\sum_{p}\left(\sum_{\nu\notin S,\nu \mid p}\mathfrak{f}_\nu\cdot (t\cdot P)_\nu\log(p)\right)=\sum_p(t\cdot P)_p\log (p).
\end{align}

This is the first stage in the construction of the product of local terms in \ref{eq:ESZ-Bdecomp}. We now show how to compute $\lambda_{S, \ul{X}}$ in the following situations, which covers many cases of interest. 

\begin{example}\label{ex:P1example}
Let $\fX=(\bP^1_\bQ; (P_0,a),(P_1,b),(P_{\infty},c) )$ with $S=\{\nu_\infty\}$ and $P_1=[0:1],P_1=[1:-1],P_{\infty}=[1:0]$. Fix $[x:y]=t\in \bP^1(\bZ)$ with $x,y$ coprime integers. Let $\nu$ be the place associated to a prime $p_\nu$. Then $\lambda_\nu(P_0,t)=p_\nu^{(t\cdot P_0)_\nu}$ and $(t\cdot P_0)_\nu$ is the largest $m$ such that $[x:y]=[0:1] \mod p^m$ which is $\ord_\nu(x)$. Therefore we compute
\[\lambda(P_0,t)=\prod_{\nu\neq \nu_\infty}p_\nu^{(t\cdot P_0)_\nu}=\prod_{\nu\neq \nu_\infty}p_\nu^{\ord_\nu(x)}=\mid x\mid \]

Similarly we compute

\[\lambda(P_1,t)=\mid x+y\mid,\lambda(P_\infty,t)=\mid y\mid \]

\end{example}

\begin{example}
Let $E/\bQ$ be the elliptic curve given by affine equation $y^2=x^3+Ax+B$ with $A,B\in \bZ$. Let $S$ be the set of primes of bad redudction of $E$ and the infinite place. Let $P=(a,b)\in E(\bZ)$.  If $t=(x,y)$ is any other point in $E(\bZ)$ then moving to the projective model we have that for any $\nu\notin S$ with associated prime $p_\nu$  that $(t\cdot P)_\nu=\max_m ([x:y:1]=[a:b:1]\mod p^m)$ which is the same as  $\min(\ord_p(x-a),\ord_p(y-b))$. We then have

\[\lambda(P,t)=\prod_{\nu\notin S}p_\nu^{\min(\ord_\nu(x-a),\ord_\nu(y-b))}.\]

Let us be more explicit and take $E$ to be given by $y^2=x(x-1)(x-\lambda)$ with $\lambda$ an integer and $S$ the set of primes of dividing $\lambda$ along with the infinite prime. Take $P=(0,0)=O$. Now take $(x,y)
\in E(\bZ)$ with $\ord_p(x)=a>1$ and $\ord_p(y)=b>1$. Then we have $x=x^\prime p^a,y=y^\prime p^b$ and so

\[p^{2b}(y^\prime)^2=p^ax^\prime(p^ax^\prime-1)(p^ax^\prime-\lambda).\]

If $b>a$ then $2b-a>a$ and $p^{2b-a}\mid x^\prime(p^ax^\prime-1)(p^ax^\prime-\lambda)$. Since $\gcd(x^\prime(p^ax^\prime-1),p)=1$ we have that $p\mid (p^ax^\prime-\lambda)\Rightarrow p\mid \lambda$ and so $p$ is a prime of bad reduction. Thus $b\leq a$ when $\ord_p(x),\ord_p(y)>1$. Therefore we have

\[\lambda(O,(x,y))=
    \prod_{p\nmid \lambda\textnormal{ and }p\mid \gcd(x,y)}p^{\ord_p(y)}.
\]

Thus $\lambda(O,(x,y))$ measures the size of the $Y$ coordinate of a point which has a non-trivial greatest common divisor.

\end{example}

	
We will define \[H_{\L_{P_i}^{\otimes d_i}}=\sigma_{m_i,d_i}(\lambda(P_i,t))\] for some function $\sigma_{m_i,d_i}\colon \bZ_{\geq 0}\ra \bZ_{\geq 0}$.  Recall that
\[\lambda(P,t)=\prod_{\nu\notin S}\textnormal{N}(\fp_\nu)^{(t\cdot P)_\nu}=\prod_pp^{(t\cdot P)_p}.\]

The information contained in $\lambda(P,t)$ are the intersection multiplicities $(t\cdot P)_p$. As we are only interested in these exponents we will demand that

\[\sigma_{m_i,d_i}(\lambda(P_i,t))=\prod_{p}\sigma_{m_i,d_i}(p^{(t\cdot P)_p})=\prod_p p^{N_{m_i,d_i}((t\cdot P)_p))}\]

where $N_{m_i,d_i}\colon \bZ_{\geq 0}\ra \bZ_{\geq 0}$. Recall that $\L_{m_{P_i}}^{\otimes m_i}=\pi^*_{\fX}\O_X(P_i)$. In other words $\L_{m_{P_i}}$ is torsion in the group $\Pic(\fX)/\pi^*_\fX(\Pic(X))$ of order $m_i$. To reflect this we demand that $N_{m_i,d_i}(p^{(t\cdot P_i)_p})$ only depends on $(t\cdot P_i)_p\mod m_i$.\\

We remained purposely agnostic until this point as to which properties $N_{m_i,d_i}$ must have to allow for a wide variety of potential height functions in all cases. However, we suspect that we can precisely identify the relevant functions $N_{m_i,d_i}$. See question \ref{ques:eszbques}. We now specify the specific functions $N_{m_i,d_i}$ which we will use.

\begin{definition}\label{def:sizefunctions}
We now construct the functions $N_{m}\colon \bZ/m\bZ\ra \bZ_{\geq 0}$ that will be relevant to us. Let $[0],...,[m-1]$ be a set of representatives of the equivalence classes of $\bZ/m\bZ$.  

\begin{enumerate}
	\item The most canonical approach is to consider the remainder of any representative of $n$ modulo $m$ of a given equivalence class. We define the canonical size function as 
	
	\[N_{m,\textnormal{can}}([r])=r\]
	
	 for $0\leq r<m$.  
	\item 
	$\bZ/m\bZ$ also has a canonical involution that can be composed with the canonical function defined above. We define
	
	\[N_{m,-}([r])=\begin{cases} 0 & r=0 \\
		m-r & r\neq 0
	\end{cases} \]
	
	Alternatively we have
	
	\[N_{m,-}([r])=N_{m,\textnormal{can}}([-r]).\]

	\item More generally one could consider
	
	\begin{equation}\label{eq:mainsize}
	    N_{m,d}([r])=N_{m,\textnormal{can}}([-rd])
	\end{equation}
	
	for any $d\in \bZ$. With this notation, $N_{m,-}=N_{m,1}$.

	
\end{enumerate}
    
\end{definition}

We can now define the height functions associated to the stacky line bundles $\L_{P_i}^{\otimes d_i}$

\begin{definition}
Fix an $M$-curve $\fX=(X:(P_1,m_1), \cdots (P_r,m_r))$ defined over a number field $K$. Choose a finite set of primes $S$ of $\O_K$ containing all the primes of bad reduction and all infinite places of $K$. We further choose a smooth and proper model $\underline{X}$ of $X$ over $\O_{K,S}$.
Consider the line bundle $\L_{P_i}^{d_i}$ on $\fX$ given by \ref{lem:bundlerep}. The stacky height function associated $\L_{P_i}^{d_i}$ is a function 

\[H_{\L_{P_i}^{d_i}}\colon X(K)\setminus{P_i}\ra \bZ_{\geq 0}\]

defined by

\[H_{\L_{P_i}^{d_i}}(t)=\left(\prod_p p^{N_{m_i,d_i}((t\cdot P_i)_p)}\right)^{\frac{1}{m_i}}\]

where $N_{m_i,d_i}\colon \bZ_{\geq 0}\ra \bZ_{\geq 0}$ is the function defined by equation \ref{eq:mainsize}.

\end{definition}








Returning to the situation of a fixed $M$-curve $\fX=(X:(P_1,m_1), \cdots (P_r,m_r))$ we now may define a stacky height machine for line bundles on $\fX$.

\begin{definition}[The stacky height machine]\label{def:stackycurveheightmachine}
Fix an $M$-curve $\fX=(X:(P_1,m_1), \cdots (P_r,m_r))$ defined over a number field $K$. Choose a finite set of primes $S$ of $\O_K$ containing all the primes of bad reduction and all infinite places of $K$. We further choose a smooth and proper model $\underline{X}$ of $X$ over $\O_{K,S}$.
	Consider the line bundle
	
	\begin{equation*}
		\L\cong \pi_\fX^*M\otimes \prod_{i=1}^r\L_{P_i}^{\otimes d_i}
	\end{equation*} 

on $\fX$.  Let $\mod m_i\colon \bZ/m_i\bZ\ra \bZ_{\geq 0}$ be the function that identifies an integer with its remainder modulo $m_i$. The stacky height associated to $\L$ is defined to be

\begin{align}
H_\L(t)&=H_{\L}^{\textnormal{st}}\cdot \prod_{i=1}^rH_{\fX,\L_{P_i}^{\otimes d_i}}\\
&=H_M\cdot \prod_{i=1}^r\left(H_{P_i}^{d_i}\cdot \prod_{p}p^{N_{m_i,d_i}((t\cdot P_i)_p)}\right)^{\frac{1}{m_i}}\\
&=H_M\cdot \prod_{i=1}^r\left(H_{P_i}^{d_i}\cdot \prod_{p}p^{-d_i(t\cdot P_i)_p \mod m_i}\right)^{\frac{1}{m_i}}
\end{align}

We call \[H_{\L}^{\textnormal{st}}=H_M\cdot \prod_{i=1}^rH_{P_i}^{\frac{d_i}{m_i}}\] 

the \emph{classical part} of the height $H_\L$ and

\[H_{\L,\textnormal{stacky}}=\prod_{i=1}^rH_{\fX,\L_{P_i}^{\otimes d_i}}\]

the \emph{stacky} part of the height $H_\L$.
\end{definition}

\begin{corollary}
Fix an $M$-curve $\fX=(\bP^1:(P_1,m_1), \cdots (P_r,m_r))$ defined over a number field $K$. Choose a finite set of primes $S$ of $\O_K$ containing all the primes of bad reduction and all infinite places of $K$. We further choose a smooth and proper model $\underline{X}$ of $X$ over $\O_{K,S}$.
The canonical height function $H_{K_\fX}$ arises from choosing the \emph{canonical} size function $N_{m_i,\textnormal{can}}\colon \bZ/m_i\bZ\ra \bZ_{\geq 0}$ that identifies $\bZ/m_i\bZ$ with $0,...,m_i-1$. In other words,

\[H_{K_\fX}=H_{\bP^1}^{\deg K_\fX}\prod_{i=1}^r\left(\prod_pp^{(t\cdot P_i)_p \mod m_i}\right)^{\frac{1}{m_i}}.\]
On the other hand the anti-canonical height function $H_{-K_\fX}$ arises by choosing the \emph{anti-canonical} size function $N_{m_i,-}\colon \bZ/m_i\bZ\ra \bZ_{\geq 0}$  that sends a residue $r$ to $m_i-r$.  That is,

\[H_{-K_\fX}=H_{\bP^1}^{-\deg K_\fX}\prod_{i=1}^r\left(\prod_pp^{-(t\cdot P_i)_p \mod m_i}\right)^{\frac{1}{m_i}}.\]

\end{corollary}
\begin{proof}
This follows directly from the definition, \ref{cor:P^1stableheight}, and the fact that $K_\fX$ corresponds to the line bundle

\[\pi^*_\fX \O_{\bP^1}(K_{\bP^1})\otimes \prod_{i=1}^r\L_{P_i}^{\otimes m_i-1}.\]
\end{proof}

We now introduce two multiplicative functions $\phi_m$ and $r_m$ that depend on an integer $m\geq 1$. The function $\phi_m$ will be the multiplicative function associated to $N_{m,-}\colon \bZ/m\bZ\ra \bZ_{\geq 0}$ and $r_m$ will correspond to $N_{m,\textnormal{can}}\colon \bZ/m\bZ\ra \bZ_{\geq 0}$. These functions will be crucial to understanding heights on the $M$-curve $\fX(\bP^1_\bQ,(P_1,m_1),\dots,(m_r,P_r))$. The functions $\phi_m$ and $r_m$ are \emph{dual} to one another in a certain sense. This can be seen on the level of the functions $N_m$ as

\[N_{m,-}(-x)=N_{m,\textnormal{can}}(x).\]
This duality is key to understanding the non-linear aspects of heights on stacky curves. For example, the differences between the height functions $H_\L$ and $H_{\L^{-1}}$. To understand the duality between $\phi_m$ and $r_m$ we define a function $\textnormal{rad}_m(x)$. We have that $\phi_m$ and $r_m$ satisfy a functional equation

\[r_m\phi_m=\textnormal{rad}_m^m.\]

The function $\textnormal{rad}_m(x)$ provides a lower bound on the usual radical of integer, in other words

 \[\textnormal{rad}_m(x)\leq \textnormal{rad}(x)\]
 
These two properties allows us to relate the stacky Vojta conjecture of \cite{ESZ-B} to the $abc$-conjecture. 



\begin{definition}
	Let $x$ be a positive integer with prime factorization $x=\prod_{p}p^{\ord_p(x)}$. 
	
	\begin{enumerate}
		\item Using the division algorithm we define integers $q_{p,m}(x),r_{p,m}(x)$ by the equation $\ord_p(x)=q_{p,m}(x)m+r_{p,m}(x)$ where $0\leq r_{p,m}(x)<m$. 
		\item Define $q_m(x)=\prod_{p}p^{q_{p,m}(x)}$ and $r_m(x)=\prod_{p}p^{r_{m,p}(x)}$. 
		\item Set $\phi_m(x)$ to be the least positive integer such that $x\phi_m(x)$ is a $m^{th}$ power.
		\item We define the $m$-radical of $x$ to be the product of all prime divisors of $x$ whose order is not a factor of $m$. In other words, \[\rad_m(x)=\prod_{p\textnormal{ s.t. }\ord_p(x)\neq 0 \mod m}p\].

	\end{enumerate} 
	
\end{definition}

We think of $r_m$ as the function associated to $N_{m,\textnormal{can}}$ and $\phi_m$ the function associated to $N_{m,-}$ 

\begin{proposition}\label{prop:betterdescription}

Let $x\in \bZ_{\geq 0}$. Then we have

\begin{align*}
r_m(x)&=\prod_{p}p^{r_{p,m}(x)}=\prod_{p}p^{N_{m,\textnormal{can}}(\ord_p(x))},\\
\phi_m(x)&=\prod_{p\textnormal{ s.t. }\ord_p(x)\neq 0 \mod m}p^{m-r_{p,m}(x)}=\prod_{p\mid \rad_m(x)}p^{m-r_{p,m}(x)}=\prod_pp^{N_{-,m}(\ord_p(x))}.
\end{align*}

In particular we have

\begin{align*}
	H_{\fX,\L_{P_i}^{-1}}(t)&=r_{m_i}(\lambda(P_i,t))^{\frac{1}{m_i}},\\
	H_{\fX,P,\L_{P_i}^{\otimes d_i}}(t)&=\phi_m(\lambda(P_i,t)^{d_i})^{\frac{1}{m_i}}.
\end{align*}
\end{proposition}

 From the above formulas we obtain the following.

\begin{proposition}
Fix and integer $m>1$ and let $x\in \bZ$. Then
\begin{enumerate}
	\item Both $r_m$ and $\phi_m$ are multiplicative functions.
	\item $\phi_m(x)r_m(x)=\rad_m(x)^m$
	\item $\phi_m(x)=1\iff r_m(x)=1$.
	\item If $m=2$ then $r_m(x)=\phi_m(x)$. 
\end{enumerate}
\end{proposition}
The above gives an alternative relationship between these height functions> After taking $\frac{1}{m}$ powers the height functions decompose $\rad_m$ into a "canonical" part and "anti-canonical" part.  With these definitions in hand we define height functions on the $M$-curve $\fX$ to be a function built out  from the functions $\phi_{m_P}(\lambda(P,t))^{\frac{1}{m_P}},$ and the classical Weil heights on the coarse space $X$.

\begin{corollary}
Fix an $M$-curve $\fX=(X:(P_1,m_1), \cdots (P_r,m_r))$ defined over a number field $K$. Choose a finite set of primes $S$ of $\O_K$ containing all the primes of bad reduction and all infinite places of $K$ and a smooth and proper model $\underline{X}$ of $X$ over $\O_{K,S}$.
	Consider the line bundle
	
	\begin{equation*}
		\L\cong \pi_\fX^*M\otimes \prod_{i=1}^r\L_{P_i}^{\otimes d_i}
	\end{equation*}
with $0\leq d_i\leq m_i-1$. Then 

\begin{equation}
    H_\L=H_M\cdot \prod_{i=1}^rH_{P_i}^{\frac{d_i}{m_i}}\cdot \prod_{i=1}^r \phi_{m_i}(\lambda(P_i,t)^{d_i})^{\frac{1}{m_i}}
\end{equation}

\end{corollary}

In particular, when $X=\bP^1$ we have

\begin{align}
H_{K_\fX}(t)&=H_{\bP^1}^{-\chi(\fX)}(t)\cdot\prod_{i=1}^rr_{m_i}(\lambda(P_i,t))^{\frac{1}{m_i}},\\
H_{-K_\fX}(t)&=H_{\bP^1}^{\chi(\fX)}(t)\cdot\prod_{i=1}^r\phi_{m_i}(\lambda(P_i,t))^{\frac{1}{m_i}},
\end{align}
where $\chi(\fX)=-\deg K_{\fX} $ is the Euler characteristic of $\fX$.\\

One interesting feature of the heights given by (\ref{def:stackycurveheightmachine}) is that they differentiate between \emph{rational} and \emph{integral} points on $M$-curves; this is one of the main features of these heights. This provides additional evidence that Ellenberg, Satriano, and Zuerick-Brown's theory of heights on algebraic stacks is an appropriate one. 

The connection to \cite{ESZ-B} and our heights is the following which is proved in \ref{pf:proofofmaincomparison}.

\begin{theorem}\label{thm:maincomarison}
Fix an $M$-curve $\fX=(\bP^1_\bQ:(P_1,m_1), \cdots (P_r,m_r))$. Choose $S$ to be the set of all finite primes of $\bZ$ and let $\bP^1_\bZ$ be the canonical model of $\bP^1_\bQ$ over $\bZ$. Let $\L$ be a line bundle on $\fX$. Let $H^{\textnormal{ESZ-B}}_\L$ be the height constructed in \cite{ESZ-B} associated to $\L$. Then there is some constant $C>0$ with

\begin{equation}
    C^{-1}\cdot H^{\textnormal{ESZ-B}}_\L\leq H_\L\leq C\cdot H^{\textnormal{ESZ-B}}_\L.
\end{equation}

That is to say, up to a constant the stacky height machine of \ref{def:stackycurveheightmachine} agrees with the ESZ-B height machine of \cite{ESZ-B} when the coarse space is $\bP^1_\bQ$. 

\end{theorem}

We now explain how the functions $\phi_m$ and $r_m$ can be used to understand the difference between $H_{\L}$ and $H_{\L^{\otimes n}}$.

\begin{proposition}\label{prop:dualityprop}
	
Let $m\in \bZ_{\geq 1}$ and choose an integer $d\geq 0$.

\[\phi_m(x^{-d \mod m})=r_m(x^{d \mod m}) .\]
\end{proposition}
\begin{proof}
Since $\phi_m$ is multiplicative it suffices to prove the statement for $x=p^a$ where $p$ is some prime. Note that $\phi_m((p^a)^{n}=p^{m-na \mod m})$. Therefore $\phi_m((p^a)^{-d \mod m})=p^{m+da \mod m}$. On the other hand, $r_m((p^a)^{d \mod m}=p^{da \mod m}=p^{m+da\mod m}$ as needed. 	
	
\end{proof}

\begin{theorem}[Duality theorem]\label{thm:dualitythm}
Let $\fX=\fX(X:(P_1,m_1),...,(P_r,m_r))$ be an $M$-curve and let

	\[\L\cong \pi_\fX^*M\otimes \prod_{i=1}^r\L_{P_i}^{\otimes d_i}\]

with $0\leq d_i\leq m_i-1$. Fix an integer $n\neq 0$ and write $n_i=nd_i\mod m_i$. 

\begin{enumerate}
	\item Then we always have
	
	\begin{align*}
	H_{\L^{\otimes n}}&=(H_{\L}^{\textnormal{st}})^n\cdot \prod_{i=1}^r\phi_{m_i}(\lambda(P_i,t)^{n_i})^{\frac{1}{m_i}}\\
		&=(H_{\L}^{\textnormal{st}})^n\cdot \prod_{i=1}^r\phi_{m_i}(\lambda(P_i,t)^{nd_i\mod m_i})^{\frac{1}{m_i}}.
	\end{align*}
\item If $n>0$ then

\begin{align*}
	H_{\L^{\otimes -n}}=(H_{\L}^{\textnormal{st}})^{-n}\cdot \prod_{i=1}^rr_{m_i}(\lambda(P_i,t)^{nd_i \mod m_i})^{\frac{1}{m_i}}
\end{align*}

In particular,

	\[H_{\L^{-1}}=(H_{\L}^{\textnormal{st}})^{-1}\cdot \prod_{i=1}^rr_{m_i}(\lambda(P_i,t)^{d_i})^{\frac{1}{m_i}}.\]

\end{enumerate} 
\end{theorem} 
\begin{proof}
Write $nd_i=m_iq_i+r_i$ with $0\leq r_i<m_i$. Note that $\L_{P_i}^{\otimes m_i}=\pi^*_\fX \O_X(P_i)$. So we have that that  

\begin{equation}\label{eq:dualityeq1}
	\L^{\otimes n}=\pi^*_\fX(M^{\otimes n})\cdot\prod_{i=1}^r\pi_\fX^*(\O_X(q_iP_i))\cdot \prod_{i=1}^r\L_{P_i}^{\otimes r_i}.
\end{equation}	

Therefore we have

\begin{align*}
H_{\L^{\otimes n}}^{\textnormal{st}}&=H_m^n\cdot\prod_{i=1}^rH_{P_i}^{q_i}\cdot \prod H_{P_i}^{\frac{r_i}{m_i}}=H_m^n\cdot\prod_{i=1}^rH_{P_i}^{\frac{m_iq_i+r_i}{m_i}}
\\
&=H_m^n\cdot\prod_{i=1}^rH_{P_i}^{\frac{nd_i}{m_i}}=\left(H_M\cdot \prod_{i=1}^rH_{P_i}^{\frac{d_i}{m_i}}\right)^n=(H_{\L}^{\textnormal{st}})^n.
\end{align*}

Now we have that

\begin{align*}
	H_{\L^{\otimes n}}&=H_{\L^{\otimes n}}^{\textnormal{st}}\prod_{i=1}^r\phi_{m_i}(\lambda(P_i,t)^{r_i})^{\frac{1}{m_i}}\\
	&=H_{\L^{\otimes n}}^{\textnormal{st}}\prod_{i=1}^r\phi_{m_i}(\lambda(P_i,t)^{nd_i \mod m_i})^{\frac{1}{m_i}}
\end{align*}	

by the definition of the $r_i$. Now let $n>0$. By proposition \ref{prop:dualityprop} we have $\phi_m(\lambda(P_i,t)^{-dn \mod m})=r_m(\lambda(P_i,t)^{nd_i \mod m_i})$ giving the desired conclusion of (2). 

\end{proof}

\begin{corollary}

	Let $\fX=\fX(X:(P_1,m_1),...,(P_r,m_r))$ be an $M$-curve and let
	\[\L\cong \pi_\fX^*M\otimes \prod_{i=1}^r\L_{P_i}^{\otimes d_i}\]
	
	with $0\leq d_i\leq m_i-1$. Then

\[H_\L(t)\cdot H_{\L^{-1}}(t)=\prod_{i=1}^r\textnormal{rad}_{m_i}(\lambda(P_i,t)^{d_i})\]

The following special cases are of particular interest.

\begin{enumerate}
	\item 
	Let $\fX=(\bP^1_
	\bQ,([1:-1],m))$. Then if $x,y$ are coprime integers we have
	
	\[H_{K_{
			\fX}}([x:y])\cdot H_{-K_{\fX}}([x:y])=\textnormal{rad}_{m}(x+y)\leq \textnormal{rad}(x+y).\]
		
\item 	Let $\fX=(\bP^1_
\bQ,([0:1],m_1),([1:0],m_2),([1:-1],m_3))$. Let $t=[x:y]$ with $x,y,x+y$ coprime integers then we have 

\[H_{K_\fX^{}}(t)\cdot H_{-K_\fX}(t)=\textnormal{rad}_{m_1}(x)\textnormal{rad}_{m_2}(y)\textnormal{rad}_{m_3}(x+y)\leq \textnormal{rad}(xy(x+y)).\]
\end{enumerate} 

\end{corollary}
\begin{proof}
We have
\[H_\L=H_\L^{\textnormal{st}}\cdot \prod_{i=1}^m\phi_{m_i}(\lambda(P,t)^{d_i})^\frac{1}{m_i}.\]

By \ref{thm:dualitythm} we have

\[H_{\L^{-1}}=(H_\L^{\textnormal{st}})^{-1}\cdot \prod_{i=1}^mr_{m_i}(\lambda(P,t)^{d_i})^\frac{1}{m_i}.\]

Therefore

\[H_\L(t)\cdot H_{\L^{-1}}(t)=\prod_{i=1}^r\phi_{m_i}(\lambda(P_i,t)^{d_i})^{\frac{1}{m_i}}\cdot r_{m_i}(\lambda(P_i,t)^{d_i})^{\frac{1}{m_i}}. \]

Since $\phi_{m_i}(x)\cdot r_{m_i}(x)=\textnormal{rad}_{m_i}(x)^{m_i}$ we have

\begin{align*}
	\prod_{i=1}^r\phi_{m_i}(\lambda(P_i,t)^{d_i})^{\frac{1}{m_i}}\cdot r_{m_i}(\lambda(P_i,t)^{d_i})^{\frac{1}{m_i}}
	&=\prod_{i=1}^r\left(\phi_{m_i}(\lambda(P_i,t)^{d_i})\cdot r_{m_i}(\lambda(P_i,t)^{d_i})\right)^{\frac{1}{m_i}}\\
	&=\prod_{i=1}^r\left(\textnormal{rad}_{m_i}(\lambda(P_i,t)^{d_i})^{m_i}\right)^{\frac{1}{m_i}}=\prod_{i=1}^r\textnormal{rad}_{m_i}(\lambda(P_i,t)^{d_i})
\end{align*}	

as needed. The remaining results follows from \ref{ex:P1example}.

\end{proof}
This suggests that the product of a stacky height with the stacky height of its dual may be of arithmetic interest. Indeed, for a single stacky point point we see that the product of a stacky canonical height and anti-canonical height gives a lower bound

\[H_{K_{\fX}}([x:y])\cdot H_{-K_{\fX}}([x:y])\leq\textnormal{rad}(x+y)).\]

Moving to three stacky points gives a lower bound 

\[H_{K_{\fX}}([x:y])\cdot H_{-K_{\fX}}([x:y])\leq\textnormal{rad}(xy(x+y)).\]

It is now apparent that stacky heights can be related to the $abc$-conjecture. In fact the $abc$-conjecture would now follow from the following statement about stacky heights: Let $\fX=(\bP^1_\bQ,([0:1],m_1),([1:0],m_2),([1:-1],m_3))$ and a constant $C_{\epsilon,\fX}>0$ such that

\[C_{\epsilon,\fX}\cdot \lambda(P_i,t)\leq \left(H_{K_\fX}(t)\cdot H_{-K_\fX}(t)\right)^{1+\epsilon}\]

where $P_1=[0:1],P_2=[1:0]$ and $P_3=[1:-1]$.

\begin{question}
Let $\fX=(X:(P_1,m_1),\dots,(m_r,P_r))$ be an $M$-curve. Then can one always find constants $\epsilon(\fX)>0$ and $C_{\fX}>0$ such that

\[C_{
\fX}\cdot \lambda(P_i,t)\leq \left(H_{K_\fX}(t)\cdot H_{-K_\fX}(t)\right)^{1+\epsilon(\fX)}.\]

\end{question}

Any positive answer to such a question would be extremely interesting, as taking $\fX=(\bP^1_\bQ,([0:1],m_1),([1:0],m_2),([1:-1],m_3))$ would lead to a bound

\[C_\fX\cdot \vert x+y\vert=C_\fX\cdot\lambda([1:-1],[x:y]))\leq \left(H_{K_\fX}([x:y])\cdot H_{-K_\fX}([x:y])\right)^{1+\epsilon(\fX)}\leq \textnormal{rad}(xy(x+y))^{1+\epsilon(\fX)}\]

when $x,y,x+y$ are coprime integers.

\subsection{Integral Points on $M$-Curves}

Here we show that the the height functions defined in  (\ref{def:stackycurveheightmachine}) can be used to obtain information about integral points on $\fX$. Fix a stacky curve $\fX=\fX(X:(\Bp,\Bm))$. \\

Let $H_{\fX,P_i}$ be the height function associated to the the line bundle $\L_{P_i}$. We have $H_{\fX,P_i}(t)=\phi_{m_i}(\lambda(P_i,t))^{1/m_i}$. We will be interested in the height function
\[H_{\D_\fX}(t)=\prod_{i=1}^r\phi_{m_i}(\lambda(P_i,t))^{1/m_i}.\]

associated to the line bundle
\[\D_\fX=\prod_{i=1}^r\L_{P_i}.\]

This is the \emph{stacky part} of the anti-canonical height. We could also work with the stacky part of the canonical height

\[H_{-\D_\fX}=\prod_{i=1}^rr_{m_i}(\lambda(P,t))^{\frac{1}{m_i}}\]We will prove that the set of integral points is contained in the set of points where $H_{\D_\fX}(t)=1$. When we take $K=\bQ$ we see that this condition is sufficient. In other words the $S$-integral points are those where the stacky part of the height is trivial. Following Darmon \cite{MCurves} we have the following notion of integral points on an $M$-curve:
\begin{definition}[Darmon] \label{DarInt} 
	
	Let $\fX=(X:(P_1,m_1), \cdots (P_r,m_r))$ be a $M$-curve over a number field $K$, $S$ a finite set of places of $K$ containing all primes of bad reduction. Let $\underline{X}$ be a smooth proper model for $X$ over $\O_{K,S}$. The $(\underline{X},S)$-integral points of $\fX$ (usually abbreviated to $S$-integral points of $\fX$) are the points $t\in X(K)$ such that
	
	\begin{equation}\label{eq:SInt}
		(t\cdot P)_\nu\equiv 0\mod m_P    
	\end{equation}
	for all $P\in X(K)$ and  $\nu\notin S$. 
	
\end{definition}

We shall prove the following theorem.

\begin{theorem} \label{height MT}
	Let $\fX=(X:(P_1,m_1), \cdots ,(P_r,m_r))$ be an $M$-curve over $K$ satisfying our assumptions and choose $S$ and a model $\underline{X}$ as we have specified. Then we have the following conclusions.
	
	\begin{enumerate}
		\item \[\fX(\O_{K,S,\underline{X}})\subseteq \bigcap_{i=1}^r\fX(P_i;K)\]
		where  $\fX(P_i;K)=\{t\in X(K)\colon H_{\fX,P_i}(t)=1\}$.
		
		\item If $K=\bQ$ then
		\[\fX(\O_{K,S,\underline{X}})= \bigcap_{i=1}^r\fX(P_i;K).\]
		In particular, the set of $S$-integral points of $\fX$ is precisely the set of points where $H_{\fX,P_i}(t)=1$ for all $i=1,...,r$. 
		
	\end{enumerate}

\end{theorem} 
Fix a prime $\nu\notin S$ and write $(t\cdot P)_\nu=m_{P}^{e_{\nu,P}(t)}\cdot q_{\nu,P}(t)$ where $e_{\nu,P}(t)\geq 0$ and $q_{\nu,P}(t)\geq 0$ is not divisible by $m_P$. In other words $q_{\nu,P}(t)$ is the $m_P$-free part of $(t\cdot P)_\nu$. Set $\textnormal{N}(\fp_\nu)=p_\nu^{f(\nu)}$. Then \begin{equation}\label{eq:locterms}
	\lambda_\nu(P,t)=p_\nu^{f(\nu)(t\cdot P)_\nu}=p_\nu^{m_{P}^{e_{\nu,P}(t)}\cdot q_{\nu,P}(t)\cdot f(\nu)}
\end{equation}
and
\begin{equation}\label{eq:expanded}
	\lambda(P_i,t)=\prod_{\nu\notin S}p_\nu^{m_{P_i}^{e_{\nu,P_i}(t)}\cdot q_{\nu,P_i}(t)\cdot f(\nu)}
\end{equation}
Using the functions $\lambda(P,t)$ we can find subsets of the rational points that contain all integral points. 
\begin{proposition}\label{prop:cuttingout}
	Suppose that $m_P>1$. Define $\fX(P;K)=\{t\in X(K)\colon H_{\fX,P}(t)=1\}$. Then \[\fX(\O_{K,S,\underline{X}})\subseteq \fX(P;K).\]
\end{proposition}
\begin{proof}
Note that $H_{\fX,P}(t)=1\iff  \phi_{m_P}(\lambda(P,t))=1\iff r_{m_P}(\lambda(P,t))=1$.	Suppose now that $t$ is an $S$-integral point. Then $(t\cdot P)_\nu\equiv 0\mod m_{P}\Rightarrow e_{\nu,P}(t)>0$ for all $\nu\notin S$. Thus 
	
	\[\lambda(P,t)=\prod_{\nu\notin S}p_\nu^{m_{P_i}^{e_{\nu,P_i}(t)}\cdot q_{\nu,P_i}(t)\cdot f(\nu)}=\left(\prod_{v\notin S}p_\nu^{m_{P_i}^{e_{\nu,P_i}(t)-1}\cdot q_{\nu,P_i}(t)\cdot f(\nu)}\right)^{m_P},\]
	whence $H_{\fX,P}(t)=1$ as $\lambda(P,t)$ is an $m_P$-power.
\end{proof}
We see that each point $P$ with multiplicity $m_P>1$ imposes a height dropping condition on the set of integral points. Thus to study integral points it suffices to study

\[\fX(\O_{K,S,\underline{X}})\subseteq \bigcap_{m_P>1}\fX(P;K).\]

We obtain the following, that the stacky part of the anti-canonical height cuts out the integral points.

\begin{corollary}
Let $\fX=\fX(X:(\Bp,\Bm))$ be a stacky curve and $\D_\fX=\prod_{i=1}^r\L_{P_i}$.  Then
	
	\[\fX(\O_{K,S})\subseteq \{t\in X(K)\colon H_{\D_\fX}(t)=1\}.\]
In other words, the integral points are precisely those points where the stacky part of the anti-canonical height vanishes. 
\end{corollary}

\begin{example}
Let $X=C$ be an elliptic curve. Then $H_{\D_\fX}=H_{-K_\fX}$. In other words the integral points of the stacky elliptic curve are precisely the points where the anti-canonical height vanishes. Since $\phi_m(x)=1\iff r_m(x)=1$ we have that the $S$-integral points of a stacky elliptic curve are also precisely the points where the stacky canonical height vanishes. If $\fX$ is a scheme and $X$ is an elliptic curve then this is certainly true, as the canonical height is trivial and every rational point is integral, and vice versa. 
\end{example}


\begin{proof}[Proof of Theorem \ref{height MT}]
	We have already shown part $(1)$ of Theorem \ref{height MT} in (\ref{prop:cuttingout}). We turn to part $(2)$ and assume that $K=\bQ$. We know that that $\bigcap_{m_P>1}\fX(P;\bQ)\subseteq \fX(\O_{\bQ,S,\underline{X}})$ by (\ref{prop:cuttingout}). We now show the reverse inclusion. Let $t\in X(\bQ)$ with $H_{\fX,P}(t)=1$ for all $P$ with $m_P>1$. Since $K=\bQ$ we have that $\textnormal{N}(\fp_\nu)=p_\nu$ and $f(\nu)=1$ for all finite places $\nu$. Fix $P$ with $m_P>1$. Towards a contradiction suppose that $(t\cdot P)_{\nu_0}\neq 0 \mod m_P$ for some $\nu_0\notin S$.  Then $e_{\nu_0,P}(t)=0$. Notice that $H_{\fX,P}(t)=1$ means that $\lambda(P,t)$ is an $m_P$-power. Since if $\nu\neq \nu^\prime$ we have that $p_\nu\neq p_{\nu^{\prime}}$ we have by unique factorization of integers that
	
	\[\lambda(P,t)=\prod_{\nu\notin S}p_\nu^{m_{P}^{e_{\nu,P}(t)}\cdot q_{\nu,P}(t)}=(\prod_{\nu\notin S}p_\nu^{z_\nu(t)})^{m_P}\]
	for some integers $z_\nu(t)$. In particular for $\nu_0$ we have 
	
	\[p_{\nu_0}^{m_{P}^{e_{\nu_0,P}(t)}\cdot q_{\nu_0,P}(t)}=p_{\nu_0}^{{q_{\nu_0,P}(t)}}=p_{\nu_0}^{z_{\nu_0}(t)m_P}\]
	
	Thus $z_{\nu_0}(t)m_P=q_{\nu_0,P}(t)$ which contradicts $q_{\nu_0,P}(t)$ being indivisible by $m_P$. Thus for all $m_P>1$ and $\nu\notin S$ we have $(t\cdot P)_\nu\equiv 0 \mod m_P$ and $t$ is an $S$-integral point of $\fX$ by definition.
\end{proof}

\section{Stacky curves with coarse space $\bP^1$}\label{subsec:mainex}

We focus on the situation when the base curve is $\bP^1$. Let $X=\bP^1_\bQ$ and $S=\{\nu_\infty\}$ and take $\L$ to be $\O_{\bP^1}(1)$ so the ample height is the usual one. We consider the $M$-curve

\[\fX =(\bP^1_\bQ:(P_1,m_1),...,(P_r,m_r))\]
In this situation the $\lambda(P,t)$ can be easily computed. 
\begin{proposition}\label{prop:linform}
	Let $t=[x:y]$ and suppose that $P_i=(a_i:b_i)$ where $a_i,b_i$ are coprime integers.  Then $\lambda(P_i,t)=\mid a_iy-b_ix\mid.$
\end{proposition}
\begin{proof}
	We have that $(t\cdot P_i)_p=\max_n\{[x:y]\equiv [a_i:b_i]\mod p^n\}$. Note that this means there is some $\lambda\neq 0\mod p^n$ with $(x,y)=\lambda(a_i,b_i)\mod p^n$. Since $a_i,b_i$ have been taken coprime we may assume that $p$ does not divide $a_i$ or $b_i$. Suppose that $p\nmid a_i$ (the other case is similar). Then $\lambda=\frac{x}{a_i} \mod p^n$ and therefore $b_ix-a_iy= 0 \mod p^n$. Thus $(t\cdot P_i)_p=\ord_p(b_ix-a_iy)$. Then we have that
	
	\[\lambda(P_i,t)=\prod_{p}p^{\ord_p(b_ix-a_iy)}=\mid a_iy-b_ix\mid\]
	as needed.  
\end{proof}

\begin{definition}[Euler characteristic of $M$-curves \cite{MCurves}]
	Let $\fX=(X:(P_1,m_1),...,(P_r,m_r))$ be an $M$-curve. The \emph{Euler characteristic} of $\fX$ is defined by the formula 
	
	\[\chi(\fX)=2-2g(X)-\sum_{i=1}^r(1-\dfrac{1}{m_{i}})=2-2g(X)-r+\sum_{i=1}^r\dfrac{1}{m_{i}}\]
	where $g(X)$ is the genus of the curve $X$. We define the \emph{genus} of an $M$-curve by the formula $\chi(\fX)=2-2g(\fX)$.
\end{definition}

We now begin assembling the necessary ingredients to compare our heights constructed in \ref{sec:MCurves} to those in \cite{ESZ-B}. We first work with $H_{\L}^{\textnormal{st}}$.  

\begin{proposition}\label{prop:stablecomparison}
Let $\fX=\fX(X:(P_1,...,P_r),(m_1,...m_r)))$ and let $\L$ be a line bundle on $\fX$ with
	
	\[\L\cong \pi_\fX^*M\otimes \prod_{i=1}^r\L_{P_i}^{\otimes d_i}\]
	
	where $0\leq d_i<m_i$ and $M$ is a line bundle on $X$ with $\pi_
	\fX$ being the coarse space map. Let $H_{\L}^{\textnormal{st},\textnormal{ESZ-B}}$ be the stable height as constructed in \cite{ESZ-B}. Then $H_{\L}^{\textnormal{st}}=H_{\L}^{\textnormal{st},\textnormal{ESZ-B}}$.     
\end{proposition}
\begin{proof}
In \cite{ESZ-B} a general definition of the stable height is given. Let $m=\prod_{i=1}^rm_i$. By the properties of the stable height described in \cite{ESZ-B} we have that

\[H_{\L}^{\textnormal{st},\textnormal{ESZ-B}}=(H_{\L^{\otimes m}}^{\textnormal{st},\textnormal{ESZ-B}})^{\frac{1}{m}}\]

On the other hand, 

\[\L^{\otimes m}=\pi^*_\fX M^{\otimes m}\otimes \prod \pi_\fX^*\O_X(\frac{d_i m}{m_i} P_i)\]

It is a fact that if $L$ is a line bundle on $X$ then $H^{\textnormal{st},\textnormal{ESZ-B}}_{\pi_\fX^* L}=H_L\circ \pi_\fX$. Therefore we have

\[H^{\textnormal{st},\textnormal{ESZ-B}}_{\L^{\otimes m}}=H_{\L^{\otimes m}}=H_{M}^m\cdot \prod _{i=1}^rH_{P_i}^{\frac{d_i m}{m_i}}\]

Taking $m^{\textnormal{th}}$ roots gives the desired inequality.

\end{proof}

\begin{corollary}\label{cor:P^1stableheight}
Let $\fX=\fX(\bP^1:(P_1,...,P_r),(m_1,...m_r)))$ and let $\L$ be a line bundle on $\fX$. Then

\begin{equation}
H_{\L}^{\textnormal{st}}=H_{\bP^1}^{\deg_\fX \L}.
\end{equation}	
\end{corollary}
\begin{proof}
Let $\L$ be a line bundle on $\fX$. Then we may write

\[\L\cong \pi_\fX^*\O_{\bP^1}(d)\otimes\prod_{i=1}^r\L_{P_i}^{\otimes d_i}. \]

On $\bP^1$ we have that $\O_{\bP^1}(P_i)\cong \O_{\bP^1}(1)$. So by definition the stable height is

\[H_{\bP^1}^d\cdot \prod H_{\bP^1}^{\frac{d_i}{m_i}}=H_{\bP^1}^{d+\sum_{i=1}^r\frac{d_i}{m_i}}=H_{\bP^1}^{\deg_\fX \L}\]

as needed. 
\end{proof}

We can now precisely define the heights on a stacky curve with coarse space $\bP^1_\bQ$.

\begin{definition}\label{def:tangentheight}
	Let $\fX=(\bP^1_\bQ,(P_1,m_1),...,(P_r,m_r))$ be an $M$-curve. Set $P_i=[a_i:b_i]$ with $a_i,b_i$ coprime integers and $\ell_i(t)=ax-by$ when $t=[x:y]$ for $x,y$ coprime integers. Let $\L=\pi_\fX^* M\otimes \prod_{i=1}^r\L_{P_i}^{\otimes d_i}$ with $0\leq d_i<m_i$.
	
	Then we have
	
	\[H_\L=\max\{\vert x\vert,\vert y\vert\}^{\deg_\fX \L}\cdot \prod_{i=1}^r\phi_{m_i}(\ell_i(x,y)^{d_i})^{\frac{1}{m_i}}.\]
	
	In particular we have that
	
	\[H_{-K_{\fX}}(t)=\max\{\vert x\vert,\vert y\vert\}^{\chi(\fX)}\cdot \prod_{i=1}^r\phi_{m_i}(\ell_i(x,y))^{\frac{1}{m_i}}\]
	
	and
	
	\[H_{K_\fX}(t)=\max\{\vert x\vert,\vert y\vert\}^{-\chi(\fX)}\cdot \prod_{i=1}^rr_{m_i}(\ell_i(x,y))^{\frac{1}{m_i}}\]
	
	
\end{definition}

When $P_1=[1,0],P_2=[1,-1],P_3=[1:0]$ and $m_1=m_2=m_3=2$ we obtain the square root of the height function Question (\ref{eq:HeightDef}). We now have enough to prove our main comparison theorem.

\subsection{Proof of \ref{thm:maincomarison}}\label{pf:proofofmaincomparison}
We follow the argument given in \cite[Page 45]{ESZ-B}. Write $\L\cong \pi_\fX^*M\otimes \prod_{i=1}^r\L_{P_i}^{\otimes d_i}$. From \cite[Section 2.3]{ESZ-B} we have a decomposition

\begin{equation}\label{stackyeq1}\frac{H_{\L}^{\textnormal{ESZ-B}}}{H_{\L}^{\textnormal{st},\textnormal{ESZ-B}}}=\prod_p\delta_{\L,p}
\end{equation}
where the $\delta_{\L,p}$'s are the local discrepancies associated to $H_{\L}^{\textnormal{ESZ-B}}$. On the other hand, we have that 

\begin{equation}\label{stackyeq2}
\frac{H_\L(P)}{H^{\textnormal{st}}_\L(P)}=\prod_{i=1}^r\phi_{m_i}(\ell_i(P)^{d_i})^{\frac{1}{m_i}}.
\end{equation}

By Proposition \ref{prop:stablecomparison} we have that 
\begin{equation}\label{stackyeq3}
H_{\L}^{\textnormal{st},\textnormal{ESZ-B}}=H_{\L}^{\textnormal{st}}    
\end{equation} up to some positive constant. Therefore, it suffices to show that

\[\prod_p\delta_{\L,p}(P)=\prod_{i=1}^r\phi_{m_i}(\ell_i(P)^{d_i})^{\frac{1}{m_i}}.\]

Let $x\colon \textnormal{spec }\bQ\ra \fX$ be a rational point whose image is not any of the stacky points $P_i$. Then in \cite{ESZ-B} there is a 1-dimensional stack $\C$ called the \emph{tuning stack} and a diagram

\[\xymatrix{\textnormal{spec }\bQ\ar@/^2pc/[rr]^{x}\ar[r]\ar[dr] & \C\ar[r]^{\bar{x}}\ar[d]_\pi & \fX\ar[dl]\\ & \textnormal{spec }\bZ &}\]

Moreover, the local discrepancies can be computed at a prime $p$ by

\[p^{-\deg \pi_*\bar{x}^* \L^{-1}-(-\deg \bar{x}^*\L^{-1})}.\]

These degrees can be computed locally on $\fX$ in terms of the stacky points $P_i$. In other words,

\begin{align*}
\deg \bar{x}^* \L^{-1}&=\sum_{i=1}^r\deg_{P_i}\bar{x}^*\L^{-1},\\
\deg\pi_*\bar{x}^*\L^{-1}&=\sum_{i=1}^r\deg_{P_i}\pi_*\bar{x}^*\L^{-1}.
\end{align*}

The local degree of $\L$ at $P_i$ at $\bQ/\bZ$ is $\frac{d_i}{m_i}$. Following \cite{ESZ-B} we have the local degrees

\begin{align*}
\deg_{P_i} \bar{x}^*\L^{-1}&=-\frac{d_i\ord_p(\ell_i(x))}{m_i},\\
\deg_{P_i} \pi_*\bar{x}^*\L^{-1}&=\left\lfloor-\frac{d_i\ord_p(\ell_i(x))}{m_i}\right\rfloor.
\end{align*}

We obtain that the contribution at $P_i$ to the local discrepancy at $p$ can be written as 

\[p^{\left\lceil \frac{d_i\ell_i(x)}{m_i}\right\rceil-\frac{d_i\ell_i(x)}{m_i}}.\]

Now write $d_i\ell_i(x)=q_im_i+r_i$ with $0\leq r_i\leq m_i-1$. First suppose that $r_i=0$. Then

\[\left\lceil \frac{d_i\ell_i(x)}{m_i}\right\rceil-\frac{d_i\ell_i(x)}{m_i}=\left\lceil\frac{q_im_i}{m_i}
\right\rceil-\frac{q_im_i}{m_i}=q_i-q_i=0.\]

Now suppose that $r_i\neq 0$. Then we have that 

\begin{align*}
\left\lceil \frac{d_i\ell_i(x)}{m_i}\right\rceil-\frac{d_i\ell_i(x)}{m_i}&=\left\lceil \frac{q_im_i+r_i}{m_i}\right\rceil-\frac{q_im_i+r_i}{m_i}\\
&=q_i+\left\lceil\frac{r_i}{m_i}\right\rceil-\frac{r_i}{m_i}-q_i\\
&=1-\frac{r_i}{m_i}=\frac{m_i-r_i}{m_i}.
\end{align*}

In other words we have shown that

\begin{align*}
 \left\lceil \frac{d_i\ell_i(x)}{m_i}\right\rceil-\frac{d_i\ell_i(x)}{m_i}&=\begin{cases} 0 & d_i\ell_i(x)=0\mod m_i \\ -d_i\ell_i(x) \mod m_i & d_i\ell_i(x)\neq 0 \mod m_i \end{cases}\\
 &= \frac{-d_i\ell_i(x) \mod m_i}{m_i}=\frac{N_{m_i,d_i}(\ord_p(\ell_i(x)))}{m_i}.
\end{align*}

Thus the local discrepancies are given by

\begin{align*}
\prod_{p}\delta_{\L,p}(x)&=\prod_{i=1}^r\cdot \left(\prod_p p^{\frac{N_{m_i,d_i}(\ord_p(\ell_i(x)))}{m_i}}\right)\\
&=\prod_{i=1}^r\phi_{m_i}(\ell_i(x)^{d_i})^{\frac{1}{m_i}}.
\end{align*}

Combining equations (\ref{stackyeq1}), (\ref{stackyeq2}), and (\ref{stackyeq3}) gives the desired conclusion.






\begin{question}\label{ques:eszbques}
Consider a stacky curve $\fX=(X:(P_1,m_1),...,(P_r,m_r))$. A line bundle $M$ on $X$ with a chosen height function $H_M$ and integers $0\leq d_i<m_i$. Then does

\[H_\L(t)=H_M(t)\cdot \prod_{i=1}^rH_{P_i}^{\frac{d_i}{m_i}}(t)\cdot \prod_{i=1}^r \phi_{m_i}(\lambda(P_i,t)^{d_i})^{\frac{1}{m_i}}\]
agree with the height constructed in \cite{ESZ-B} up to a bounded constant when $X\neq \bP^1_\bQ$. In other words, does our stacky height machine recover the heights in \cite{ESZ-B} for all stacky curves over all number fields? If not, can one define \emph{different} size functions $\tilde{N}_{m_i,d_i}\colon \bZ_{\geq 0}\ra \bZ_{\geq 0}$ so that the ESZ-B height associated to $\L$ is of the form

\[H_\L(t)= \left(H_M(t)\cdot\prod_{i=1}^rH_{P_i}^{\frac{d_i}{m_i}}(t)\right)\cdot\left(\prod_{i=1}^r\left(\prod_p p^{\tilde{N}_{m_i,d_i}(t\cdot P_i)_p}\right)^{\frac{1}{m_i}}\right)\]

This result will follow provided one can show that the local degree of $\bar{x}^*\L$ with respect to $P_i$ over a prime $p$ is 

\[\frac{d_i\ord_p(\lambda(P_i,t))}{m_i}=\frac{d_i(t\cdot P_i)_p}{m_i}.\]

In this case, the argument given in \ref{thm:maincomarison} would give the desired result. One might further ask if these methods could be extended to compute the height functions of line bundles certain higher dimensional analogues of stacky curves.

\end{question}

\subsection{Morphisms of stacky curves}

We will require some results on morphisms between stacky curves.

\begin{definition}[Darmon, \cite{MCurves}]
	Let $\fX_1=(X_1,(P,m_P)),\fX_2=(X_2,(Q,m_Q))$ be $M$ curves defined over a number field $K$. A morphism of $M$-curves over $K$ is a morphism of algebraic curves $\pi\colon X_1\rightarrow X_2$ defined over $K$ such that for all $P\in X_1(K)$ with $\pi(P)=Q$ we have that \[m_Q\mid e_\pi(P)m_P\]
	where $e_\pi(P)$ is the ramification index of $\pi$ at $P$. We also define $e_{\underline{\pi}}(P)=\frac{e_{\pi}(P)m_P}{m_{Q}}$ the ramification index of $\underline{\pi}$ at $P$.
\end{definition}

Now let $\fX=(X:\bQ;(P_1,m_1),...,(P_r,m_r))$ be an $M$ curve. For any $s<r$ choose positive divisors $d_i$ of $m_i$ for $i=1,...,s$. Then there is a multiplicity lowering morphism \[\underline{\pi}(d_1,...,d_s)\colon \fX(X:(P_1,m_1),...,(P_r,m_r))\rightarrow \fX(X:(P_1,d_1),...,(P_s,d_s))\]  defined by the identity morphism on $X$.   The usefulness of this notion can be seen by the following.

\begin{proposition}[Darmon,\cite{MCurves}]
	Let $\underline{\pi}\colon \fX_1\rightarrow \fX_2$ be a morphism of $M$-curves defined over $K$. Then \[\underline{\pi}(\fX_1(\O_{K,S}))\subseteq \fX_2(\O_{K,S}).\]
	In otherwords, a morphism of $M$-curves preserves the notion of $S$-integral points.   	
\end{proposition}

\begin{lemma}\label{prop:Lprop}
Let $K$ be a field and 
\[U=\begin{bmatrix} 0 & 1\\
	 -1 &0
\end{bmatrix}\]

Define a bilinear form $L(\Bv,\Bx)=\Bv^TU^T\Bx$.
Let $T$ be a non-singular matrix with entries in $K$. Then

\[L(T\Bv,T\Bx)=(\det T)L(\Bv,\Bx).\]

\end{lemma}
\begin{proof}
We have that

\[L(T\Bv,T\Bx)=(T\Bv)^TU^TT\Bx=\Bv T^TU^TT\Bx=\Bv(UT)^TT\Bx.\]
Direct computation shows that $(UT)^TT=(\det T)U^T$.
We then have

\[L(T\Bv,T\Bx)=\Bv(UT)^TT\Bx=\Bv^T(\det T)U^T\Bx=(\det T)L(\Bv,\Bx).\]

\end{proof}
Note that if $P=[a:b]$ and $t=[x,y]$ where $a,b$ and $x,y$ are coprime integers then $\lambda(P,t)=ay-bx$ by 
(\ref{prop:linform}). In other words $\lambda(P,t)=L(P,t)$ as defined in (\ref{prop:Lprop}). 

\begin{lemma}\label{lem:autbound}
Let $P=[a:b],t=[x,y]\in \bP^1_\bQ$ with $a,b$ and $x,y$ coprime integers. Fix an integer $m>1$.Let $\alpha\colon \bP^1_\bQ\rightarrow \bP^1_\bQ$ be an automorphism. Let $\det \alpha$ be the smallest possible non-negative determinant of an integral representation of $\alpha$ and similarly for $\det \alpha^{-1}$.  Then we have

\begin{equation}
	(\rad_m(\det \alpha^{-1}))^{-(m-1)}\phi_m(\lambda(P,t))\leq \phi_m(\lambda(\alpha(P),\alpha(t))\leq\rad_m(\det\alpha)^{m-1}\phi_m(\lambda(P,t)) 
\end{equation}

\end{lemma}
\begin{proof}
Let $L$ be as (\ref{prop:linform}). Then we have that

\begin{equation}\label{eq:Lformeq}
\lambda(P,t)=L\left(\dfrac{\alpha(P)}{d_1},\dfrac{\alpha(t)}{d_2}\right)=\dfrac{(\det \alpha)}{d_1d_2}L(P,t)=\dfrac{(\det \alpha)}{d_1d_2}\lambda(P,t)
\end{equation}
for some integers $d_1,d_2$ which account for common factors of $\alpha(P)$ and $\alpha(t)$. Let $n_1=q_1m+r_1,n_2=q_2m+r_2$ be integers with $0\leq r_i<m$. 

\begin{align}
	\phi_m(p^{n_1+n_2})=
	\begin{cases} p^{m-r_1-r_2}=\dfrac{\phi_m(p^{n_1})\phi_m(p^{n_2})}{p^m} & \textnormal{ if }r_1+r_2<m\\
		p^{m-(r_1+r_2-m)}=p^{2m-r_1-r_2}=\phi_m(p^{n_1})\phi_m(p^{n_2}) &\textnormal{ if } r_1+r_2\geq m
	\end{cases}
\end{align}
Since $\phi_m$ is multiplicative we have that $\phi_m(zw)\leq \phi_m(z)\phi_m(w)\leq \rad_m(z)^{m-1}\phi_m(w)$. Therefore using (\ref{eq:Lformeq})

\[\phi_m(\lambda(\alpha(P),\alpha(t)))\leq \rad_m(\det \alpha)^{m-1}\phi_m(\lambda(P,t))\]
Applying the same reasoning using $\alpha^{-1}$ we have that

\[\phi_m(\lambda(P,t))\leq \rad_m(\det \alpha^{-1})^{m-1}\lambda(\phi(P),\phi(t))\]
Therefore we have 

\begin{equation}
	(\rad_m(\det \alpha^{-1}))^{-(m-1)}\phi_m(\lambda(P,t))\leq \phi_m(\lambda(\alpha(P),\alpha(t))\leq\rad_m(\det\alpha)^{m-1}\phi_m(\lambda(P,t)) 
\end{equation}

as required.

\end{proof}
\begin{corollary}\label{cor:3trans}
Let $\fX=\fX(\bP^1;(P_1,m_1),...,(P_r,m_r))$ Let $\alpha\colon \bP^1_\bQ\rightarrow \bP^1_\bQ$ be an automorphism of $\bP^1_\bQ$. Let $Q_i=\alpha(P_i)$ and  $\mathfrak{Y}=\mathfrak{Y}(\bP^1;(Q_1,m_1),...,(Q_r,m_r))$. Then $\alpha$ induces an isomorphism $\underline{\alpha}\colon \fX\rightarrow \mathfrak{Y}$. Let $\det\alpha$ and $\det \alpha^{-1}$ be as in (\ref{lem:autbound}). Let $D(\alpha,\fX)=\prod_{i=1}^r\rad_{m_i}(\det \alpha)^{1-1/m_i}$ and similarly define $D(\alpha^{-1},\fX)$. Suppose that $C_\phi$ is a constant (such a constant always exists) such that 

\begin{equation}
C_\alpha^{-1}H_{\bP^1}(t)^{\chi(\fX)}\leq H_{\bP^1}(\phi(t))^{\chi(\fX)}\leq C_\alpha H_{\bP^1}(t)^{\chi(\fX)}.
\end{equation}
Then 
\begin{equation}
D(\alpha^{-1},\fX)^{-1}C_\alpha^{-1} H_{\T_\fX}(t)\leq H_{\bT_{\mathfrak{Y}}}(\phi(t)) \leq D(\alpha,\fX)C_\alpha H_{\T_\fX}(t).
\end{equation}
\end{corollary}
\begin{proof}

By assumption
For each $i=1,...,r$ we have from (\ref{lem:autbound}) and our assumption that 

\begin{align*}
H_{\T_{\mathfrak{Y}}}(\phi(t))&=H_{\bP^1}(\phi(t))^{\chi(
	\mathfrak{Y})}\prod_{i=1}^r\phi_{m_i}(\lambda(Q_i,\phi(t)))^{\frac{1}{m_i}}\\ &\leq C_\alpha H_{\bP^1}(t)^{\chi(\fX)}\prod_{i=1}^r\rad_{m_i}(\det(\alpha))^{1-1/m_i}\phi_{m_i}(\lambda(P_i,t))^{\frac{1}{m_i}}\\&=C_\alpha D(\alpha,\fX)H_{\T_\fX}(t).
\end{align*}

Similarly we have

\begin{align*}
	H_{\T_{\mathfrak{Y}}}(\phi(t))&=H_{\bP^1}(\phi(t))^{\chi(
		\mathfrak{Y})}\prod_{i=1}^r\phi_{m_i}(\lambda(Q_i,\phi(t)))^{\frac{1}{m_i}}\\&\geq C_\alpha^{-1} H_{\bP^1}(t)^{\chi(\fX)}\prod_{i=1}^r\rad_{m_i}(\det(\alpha^{-1}))^{-(1-1/m_i)}\phi_{m_i}(\lambda(P_i,t))^{\frac{1}{m_i}}\\&=C_\alpha^{-1} D(\alpha^{-1},\fX)H_{\T_\fX}(t).
\end{align*}

\end{proof}

We can get slightly worse, but more understandable bounds as follows. We always have $\rad_m(x)\leq \rad(x)$. Note that we have that $\sum_{i=1}^r(1-1/m_i)=2-\chi(\fX)=2g(\fX)$. Thus we in fact have 

\begin{equation}
	C_\alpha^{-1}\rad(\det\alpha^{-1})^{-2g(\fX)}H_{\T_\fX}(t)\leq H_{\T_{\mathfrak{Y}}}(\phi(t))\leq C_\alpha\rad(\det\alpha)^{2g(\fX)}H_{\T_\fX}(t).
\end{equation}
	
Of particular note is that we see that when studying the Northcott property, we may change the height by an automorphism. Thus the Northcott property is stable under isomorphism as expected. 	

\subsection{Northcott property of the canonical height on stacky curves} We now investigate the properties of the canonical height on stacky curves, given by 
\begin{equation} \label{canheight} H_{(\Ba, \Bm)}(x,y) = \prod_{i=1}^n r_{m_i} (\ell_i(x,y))^{\frac{1}{m_i}} \cdot \max\{|x|, |y|\}^{-\delta(\Bm)}.
\end{equation}
When $\delta(\Bm) = 0$ we see that the canonical height exhibits a clear duality with the anti-canonical height, and so the same argument shows that Northcott's property will fail. When $\delta(\Bm) < 0$ we then see at once that $H_{(\Ba, \Bm)}$ will have Northcott's property as a consequence that the Weil height having Northcott's property. It remains to consider Northcott's property when $\delta(\Bm) > 0$. In this case we have $\Bm = (2, m_2, m_3)$ with
\[\frac{1}{m_2} + \frac{1}{m_3} > \frac{1}{2}.\]
It suffices to show that for any such pair $(m_2,m_3)$ there exist integers $a,b,c$ such that the curve
\[ax^2 + by^{m_2}  +c z^{m_3} = 0, \frac{1}{2} + \frac{1}{m_2} + \frac{1}{m_3} > 1\]
has infinitely many primitive integral solutions. This is the content of Beukers' paper \cite{Beu}, and we are done. 


\section{On the Northcott property of canonical and anti-canonical heights on stacky curves}
\label{NorthcottProp}

In this section we prove Theorem  \ref{intro:Northfail}, starting with Theorem \ref{Northfail}. We start with a reduction procedure of a curve $\fX(\bP^1 : (\Ba, \Bm))$ which we describe colloquially. By convention, we shall take our weight vectors $\Bm$ to have the property that
\[m_1 \leq m_2 \leq \cdots \leq m_n.\]

\begin{definition}
Consider a $M$-curve $\fX(X : (\Ba, \Bm))$ with $\Ba=(P_1,...,P_r)$. Let $\Bi=i_1,...,i_k$ be a sub-sequence of $1,2,...,r$. Then there is a morphism $\pi_{\Bi}\colon\fX(X : (\Ba, \Bm))\ra \fX(X : (\Ba^\prime, \Bm^\prime))$ where $\Ba^\prime=(P_{i_1},...,P_{i_k})$ and $\Bm^\prime=(m_{i_1},...,m_{i_k})$. The map is defined by taking the identity morphism on the coarse space $X$. We call such a morphism a \emph{totally ramified canonical covering}. 
\end{definition}
The above construction defines a morphism by the definition of a morphism of $M$- curves. It is totally ramified in the sense that if $i$ is some index that does not appear in $\Bi$ then $\pi_\Bi$ has ramification index $m_i$ at $P_i$. We use the term canonical this type of construction can be used for any $M$-curve. In particular, by taking $\Bi$ to be the empty set, we obtain the coarse space morphism.  We will show that if Theorem \ref{Northfail} holds for a totally ramified canonical covering of the shape
$\fX(\bP^1 : (\Ba^\prime, \Bm^\prime))$ where $\Ba^\prime, \Bm^\prime$ is obtained from $\Ba, \Bm$ respectively by removing a subset of indices, then it also holds for $\fX(\bP^1 : (\Ba, \Bm))$: see Theorem \ref{addpts} below. 

\begin{theorem} \label{addpts} Let $\fX(\bP^1 : (a_1, m_1), \cdots, (a_n, m_n))$ be a stacky curve. If the Northcott property fails for the height (\ref{Hdefgen}) for some totally ramified canonical cover of $\fX$, then it will also fail for $\fX$. 
\end{theorem}

\begin{proof} Let $\fX$ be given as in the statement of Theorem \ref{addpts}. We may assume, after reindexing the points if necessary, that the Northcott property for the ESZ-B height fails for the totally ramified canonical cover given by
\[\fX^{(1)} = \fX(\bP^1 : (a_1, m_1), \cdots, (a_k, n_k))\]
for some $k \leq n$. This implies that, for some positive number $C_k$ depending at most on $a_1, \cdots, a_k$ and $m_1, \cdots, m_k$, there are infinitely many integers $x,y$ such that
\[H_{(\Ba^{(k)}, \Bm^{(k)})} = \prod_{i=1}^k \phi_{m_i} (\ell_i(x,y))^{1/m_i} \max\{|x|, |y|\}^{\delta(\Bm^{(k)})} < C_k.\]

Next note that the quotient
\[\Q  = \frac{H_{(\Ba, \Bm)}(x,y)}{H_{(\Ba^{(k)}, \Bm^{(k)})}(x,y)} = \prod_{i={k+1}}^n \phi_{m_i}(\ell_i(x,y))^{1/m_i} \cdot \max\{|x|, |y|\}^{\sum_{i=k+1}^n (-1 + 1/m_i)}. \]
Observe that $\phi_m(s) \leq |s|^{m-1}$ for any integer $s$, with equality if and only if $s$ is square-free. It follows that
\[\Q \leq \prod_{i=k+1}^n |\ell_i(x,y)|^{1 - 1/m_i} \cdot \max\{|x|, |y|\}^{\sum_{i=k+1}^n (-1 + 1/m_i)},\]
and from here we immediately see from the triangle inequality that 
\[\Q \ll_{\Ba} 1.\]
Thus, by replacing $C_k$ with a larger positive number if necessary, we see that the Northcott property also fails for $H_{(\Ba, \Bm)}$ on $\fX$. \end{proof}

Now, given Theorem \ref{addpts}, it remains to consider certain \emph{minimal} choices of $\Bm$. We say that $\delta(\Bm)$ is \emph{minimally non-negative} if there is no subsequence $\Bm^\prime$ of $\Bm$ such that $\delta(\Bm^\prime) \leq 0$. We have the following lemma characterizing the minimally non-negative tuples:

\begin{lemma} \label{minneg} Suppose $\Bm = (m_1, \cdots, m_n)$ with $2 \leq m_1 \leq \cdots \leq m_n$ is minimally non-negative. Then $n \leq 4$. 
\end{lemma}

\begin{proof} Suppose $n \geq 5$. If there exist $m_i, m_j, m_k \geq 3$ then the sub-sequence $(m_i, m_j, m_k)$ satisfies $\delta((m_i, m_j, m_k)) \leq 0$, so $\Bm$ is not minimally non-negative. If $m_3 \geq 3$ then such a choice of $i,j,k$ exists, since $n \geq 5$. Therefore we may assume that $m_1 = m_2 = m_3 = 2$. But then $(2,2,2,m_4)$ satisfies $\delta((2,2,2,m_4)) \leq 0$, so $\Bm$ is not minimally non-negative.  
\end{proof}

It remains to deal with minimally non-negative tuples with $n = 3,4$. Before we proceed we will require Theorem \ref{g1 thm}, which we prove now. 

\begin{proof}[Proof of Theorem \ref{g1 thm}] As we remarked earlier, the proof given here is provided to us by A.~Shnidman in \cite{ShnidMO}. \\

For given non-singular binary quartic form $F \in \bZ[x,y]$ given by
\[F(x,y) = a_4 x^4 + a_3 x^3 y + a_2 x^2 y^2 + a_1 xy^3 + a_0y^4\]
we write $C_F$ for the curve defined by:
\[C_{F} : z^2 = F(u,v). \]
The Jacobian of the genus one curve $C_{a,b}$ is the elliptic curve $E_{a,b}$ given by
\[E_F : y^2 = x^3 - \frac{I(F)}{3}x - \frac{J(F)}{27}, \]
where $I,J$ are the basic invariants given by
\[I(F) = 12 a_4 a_0 - 3 a_3 a_1 + a_2^2, J(F) = 72 a_4 a_2 a_0 + 9 a_3 a_2 a_1 - 27 a_0 a_3^2 - 27 a_4 a_1^2 - 2 a_2^3.\]
By 2-descent, we see that $C_{F}$ corresponds to a class $c$ in $H^1(\bQ , E_{a,b}[2])$. Note that for any integer $d$, the group $H^1(\bQ, E_{F}^{(d)}[2])$ is canonically isomorphic to $H^1(\bQ , E_{F}^{(d)}[2])$ such that $c$ is the class of $C_{F}^{(d)}$ in $H^1(\bQ, E_{F}^{(d)}[2])$. \\

We now consider two cases. First suppose that $c$ does not come from 2-torsion. In this case it is immediate that $E_{F}^{(d_0)}$ has positive rank, where
\[d_0 = F(u_0, v_0), u_0, v_0 \in \bZ.\]
This is because in this case $C_{F}^{(d_0)}(\bQ) \ne \emptyset$. \\

If $c$ comes from 2-torsion, then we note that $C_{F}^{(d)}(\bQ)$ is non-empty for all $d \in \bZ$. That is, for all $d \in \bZ$ we have $C_{F}^{(d)}$ is isomorphic to $E_{F}^{(d)}$ over $\bQ$. We can then choose a class $c^\prime$ in $H^1(\bQ, E_{F}[2])$, represented by a different binary quartic form $G$, and choose $d$ such that the twist of the genus one curve $\C : z^2 = G(u,v)$ given by $\C^{(d)}$ has a rational point. This implies that $E_{F}^{(d)}$ has positive rank. Then with this choice of $d$, we find that $C_{F}^{(d)}(\bQ) \ne \emptyset$ and $E_{F}^{(d)}$ has positive rank, which completes the proof. \end{proof}

With Theorem \ref{g1 thm}, we proceed to handle minimally non-negative tuples, starting with the case $n = 4$.

\subsection{Minimally non-negative tuples with $n = 4$} 

We begin with the case $\Bm = (2,2,2,2)$, and we will need Theorem \ref{g1 thm}. By $3$-transivity of the action of $\PGL_2$ on $\bP^1$ and Lemma \ref{lem:autbound}, we may assume that three of the points are $0, 1, \infty$, corresponding to the linear forms $x,y, x+y$ in the variables $x,y$. We then write $\ell(x,y) = ax + by$ for the linear form representing the $4$-th half-point. \\ \\
We prove the following as a warm-up: 
\begin{lemma} There exist integers $a,b$ such that the stack $\fX(\bP_\bQ^1: (0, 2), (1, 2), (\infty, 2), (a/b, 2))$ has infinitely many rational points of E-S-ZB height equal to one.
\end{lemma}

\begin{proof} In this case the height is given by
\[H(x,y) = \sqf(x) \sqf(y) \sqf(x+y) \sqf(ax + by),\]
 so this is equal to one if and only if each of $x, y, x+y, ax+by$ is a square. To wit, we set
 \[x = x_1^2, y= x_2^2, x+y = x_3^2.\]
 This induces the equation
 \[x_1^2 + x_2^2 = x_3^2,\]
 which is solvable and whose (primitive) integral solutions are parametrized by
 \[x_1 = 2uv, x_2 = u^2 - v^2, x_3 = u^2 + v^2.\]
 Inserting this into $ax + by$ gives
 \[a(2uv)^2 + b(u^2 - v^2)^2 = F_{a,b}(u,v).\]
 We then fix $u = u_0, v = v_0$ so that $2u_0 v_0$ and $u_0^2 - v_0^2$ are co-prime, then solve the linear diophantine equation
 \[a (2u_0 v_0)^2 + b (u_0^2 - v_0^2) = 1.\]
 Given a solution $(a,b)$ to this diophantine equation, one obtains a genus curve defined by
 \[w^2 = F_{a,b}(u,v)\]
 which is isomorphic to an elliptic curve, since it has a rational point given by $w_0^2 = F_{a,b}(u_0, v_0)$. In particular, it must be isomorphic to its jacobian. A simple calculation shows that the jacobian of this genus one curve is given by the equation\\ 
 \begin{equation} E_{a,b} :  y^2 = x^3 - \frac{16(a^2 - ab + b^2)}{3} - \frac{64(a+b)(2a-b)(a-2b}{27} 
 \end{equation}
 \[= \left( x - \frac{4b - 8a}{3}\right)\left( x - \frac{4a - 8b}{3}\right)\left( x - \frac{4b + 4a}{3}\right),\]
 so it suffices to find $a,b$ such that $E_{a,b}$ has positive rank. We find that setting $u_0 = 1, v_0 = 5$ and $a = 17, b = -118$ that the curve $E_{a,b}$ has positive rank, and therefore $w^2 = F_{a,b}(u,v)$ will have infinitely many integral solutions $(u,v,w)$. This gives infinitely many pairs $u,v$ such that $F_{17,-118}(u,v)$ is a square. This implies our result, since 
 \[\sqf(x) \sqf(y) \sqf(x+y) \sqf(17x-118y) = \sqf((2uv)^2) \sqf((u^2 - v^2)^2) \sqf((u^2 + v^2)^2) \sqf(F_{17,-118}(u,v)) = 1.\]
\end{proof}

The general case will follow by applying the same ideas in tandem with Theorem \ref{g1 thm}. Indeed, Theorem \ref{g1 thm} gives that for any $a,b \in \bZ$ such that $F_{a,b}(u,v) = a(u^2 - v^2)^2 + 4bu^2 v^2$ is non-singular that there exists $d \in \bZ$ such that $C_F(\bQ) \ne \emptyset$ and $E_F$ has positive rank. Fixing such a $d$, we see that there are infinitely many co-prime integers $u,v,z$ such that
\[dz^2 = F(u,v).\]
Recall that in this set-up we have
\[x = (u^2 - v^2)^2, y = 4u^2 v^2, x + y = (u^2 + v^2)^2,\]
whence
\[H(x,y) = 1 \cdot 1 \cdot 1 \cdot \sqf(F(u,v)) \leq |d|.\]
This concludes the proof for the $\Bm = (2,2,2,2)$ case.  \\

Note that if $\Bm = (m_1, m_2, m_3, m_4)$ is minimally non-negative, then $m_1 = m_2 = 2$, since $\delta((3,3,3)) = 0$. Thus we may write
\[x = x_1 x_2^2, y = y_1 y_2^2\]
with $x_1, y_1$ square-free. We then write
\[\ell_3(x,y) = x + y = z_1 z_2^2 \cdots z_{m_3}^{m_3}, \ell_4(x,y) = w_1 w_2^2 \cdots w_{m_4}^{m_4}.\]
Again, we have $z_i, w_j$ are square-free for $1 \leq i \leq m_3 - 1$ and $1 \leq j \leq m_4 - 1$. \\

We now specialize to the points where $z_i = w_j = 1$ except for $i = 2$ and $j = 1,2$, as well as $x_1 = y_1 = 1$. Applying Theorem \ref{g1 thm} and using the same argument as in the $(2,2,2,2)$ case, we see that there is a choice of $w_1 = d$ such that there are infinitely many choices of $x_2, y_2, w_2, z_2$ satisfying
\[x = x_2^2, y = y_2^2, x + y = z_2^2, \ell_4(x,y) = dw_2^2.\]
The height of such a point is given by
\[\phi_2(x_2^2) \phi_2(y_2^2) \phi_{m_3} (z_2^2) \phi_{m_4} (d w_2^2) \max\{x_2^2, y_2^2\}^{1/m_3 + 1/m_4 - 1} \]
\[\ll |z_2|^{1 - 2/m_3} |w_2|^{1 - 2/m_4} \max\{|z_2|, |w_2|\}^{(2/m_3 - 1) + (2/m_4 - 1} \ll 1. \]
It follows that there are infinitely many points of bounded height, and so Northcott's property fails.

\subsection{Minimally non-negative tuples with $n = 3$} To complete the proof of Theorem \ref{Northfail}, it remains to handle the cases when $n = 3$ and $\chi(\fX) \leq 0$. We shall assume that $m_1 \leq m_2 \leq m_3$. We then note that $\delta(\Bm) \leq 0$ if and only if one of the following conditions is satisfied:
\begin{enumerate}
    \item $m_1 \geq 3$,
    \item $m_1 = 2, m_2 = 3, m_3 \geq 6$ 
    \item $m_1 = 2, m_2 \geq 4$.
\end{enumerate}

We deal with the first case. We then write
\begin{equation} \label{xyz} x = x_{1,1} x_{1,2}^2 \cdots x_{1,m_1}^{m_1}, y = x_{2,1} x_{2,2}^2 \cdots x_{2,m_2}^{m_2}, \end{equation}
\[x + y = x_{3,1} x_{3,2}^2 \cdots x_{3,m_3}^{m_3}.\]
Now set
\[x_{i,j} = 1 \text{ for } (i,j) \not\in \{(1,3), (2,3), (3,3), (3,1)\} \]
and
\[x_{1,3} = z_1, x_{2,3} = z_2, x_{3,3} = z_3, x_{3,1} = d.\]
Then the value of the height $H_{(\Ba, \Bm)}(x,y)$ in this case is given by
\begin{align*} H(z_1^3, z_2^3) & = \phi_{m_1}(z_1^3)^{1/m_1} \phi_{m_2}(z_2^3)^{1/m_2} \phi_{m_3}(d z_3^3)^{1/m_3} \max\{|z_1|^3, |z_2|^3\}^{1/m_1 + 1/m_2 + 1/m_3 - 1} \\
& \leq d^{m_3 - 1} |z_1|^{\frac{m_1 - 3}{m_1}} |z_2|^{\frac{m_2 - 3}{m_2}} |z_3|^{\frac{m_3 - 3}{m_3}} \max\{|z_1, |z_2|\}^{\left(\frac{3}{m_1} - 1 \right) + \left(\frac{3}{m_2} - 1 \right) + \left(\frac{3}{m_3} - 1 \right)}.
\end{align*} 
Observe that 
\[|z_i|^{\frac{m_i - 3}{3}} \max\{|z_1|, |z_2|\}^{-1 + \frac{3}{m_i}} \ll 1\]
for $i = 1,2,3$. \\

It remains to choose $d$ so that the plane cubic curve
\[z_1^3 + z_2^3 = dz_3^3\]
has infinitely many rational points. This is an easy consequence of the seminal work of Stewart and Top (Theorem 7, \cite{StewTop}). In particular, they showed that the number of cube-free integers $d$ with $|d| \leq X$ such that the equation
\[x^3 + y^3 = d\]
defines an elliptic curve with rank at least $2$ is asymptotically greater than $X^{1/3}$. We of course do not need such a strong statement, indeed we only need one such $d$. This completes the proof for the case $m_1 \geq 3$. \\

We proceed to handle the case $m_1 = 2, m_2 \geq 4$. Using the same notation as in (\ref{xyz}), we then set 
\[x_{i,j} = 1 \text{ for all } (i,j) \not\in \{(1,1), (1,2), (2,4), (3,4)\},\]
and set
\[x_{1,2} = z_1, x_{2,4} = z_2, x_{3,4} = z_3, x_{1,1} = d.\]
This gives a curve
\[dz_1^2 = z_3^4 - z_2^4.\]
We need to choose square-free $d$ so that this curve has infinitely many integral solutions, and such a $d$ exists by Theorem \ref{g1 thm}. The height $H_{(\Ba, Bm)}(x,y)$ is given by
\begin{align*} H_{(\Ba, \Bm)}(x,y) & = \phi_2(dz_1^2)^{1/2} \phi_{m_2}(z_2^3)^{1/m_i} \phi_{m_3}(d z_3^3)^{1/m_3} \max\{d|z_1|^2, |z_2|^4\}^{1/m_2 + 1/m_3 - 1/2} \\
& \leq d |z_2|^{1 - 4/m_2} |z_3|^{1 - 4/m_3} \max\{d|z_1|^2, |z_2|^4\}^{1/m_2 + 1/m_3 - 1/2}.
\end{align*}
Note that
\[|z_3|^4 \asymp \max\{dz_1^2, z_2^4\},\]
so we obtain the upper bound
\[\left(\frac{|z_2|}{\max\{|z_2|, |z_3|\}}\right)^{1 - \frac{4}{m_2}} \cdot \left(\frac{|z_3|}{\max\{|z_2|, |z_3|\}} \right)^{1 - \frac{4}{m_3}} \ll 1. \]
It follows that there are infinitely many integers $x,y$ such that $H_{(\Ba, \Bm)}(x,y)$ remains bounded. \\

Finally, we resolve the case $m_1 = 2, m_2 = 3, m_3 \geq 6$. In this case, we use the fact the there exist integers $a,b,c$ with $a$ square-free, $b$ cube-free, and $c$ 6-th power-free such that the equation
\begin{equation} \label{236c} ax^2 + by^3 + cz^6 = 0\end{equation}
has infinitely many primitive solutions (see \cite{DG}, Section 6.3). Thus, by fixing such a triple $(a,b,c)$ and setting
\[x_{1,1} = a, x_{2,1} x_{2,2}^2 = b, x_{3,1} \cdots x_{3,5} = c,\]
\[x_{1,2} = u_1, x_{2,3} = u_2, x_{3,6} = u_3\]
and
\[x_{5,j} = 1 \text{ for all } j \geq 7\]
we specialized a point on $\fX(\bP^1 : (0, 2), (\infty, 3), (-1, 6))$ to the curve given by (\ref{236c}). The height of such a point is then bounded in terms of $a,b,c$ only, and is thus absolutely bounded. This shows that the Northcott property fails in this case as well. \\

This concludes the proof of Theorem \ref{Northfail}. 

\subsection{Proof of Theorems \ref{Northsuccess}}  

We proceed to prove Theorem \ref{Northsuccess}. The claim when $\delta(\Bm) = 0$ is covered in Theorem \ref{Northfail}, so we will not discuss it again. When $\delta(\Bm) > 0$ we note that $n = 3$, and that in each such case there exist integers $a_\Bm, b_\Bm, c_\Bm$ such that the equation
\[a_\Bm x^{m_1} + b_{\Bm} y^{m_2} + c_{\Bm} z^{m_3} = 0\]
has infinitely many primitive integral solutions; see for example \cite{Beu}. This shows that the Northcott property fails for $H^0$. \\ 

We may now work with the case when $\delta(\Bm) < 0$. We see that the height $H^0$ is bounded below by 
\begin{equation} \label{htlb} \phi_{m_1}(x)^{1/m_1} \phi_{m_2}(y)^{1/m_2} \phi_{m_3}(x+y)^{1/m_3},\end{equation}
so it suffices to show that this quantity necessarily goes to infinity. We then use the notation from (\ref{xyz}), to obtain the equation
\begin{equation} \label{hgeq}x_{1,m_1}^{m_1} \prod_{j=1}^{m_1 - 1} x_{1,j}^j + x_{2,m_2}^{m_2} \prod_{j=1}^{m_2 - 1} x_{2,j}^{j} = x_{3,m_3}^{m_3} \prod_{j=1}^{m_3 - 1} x_{3,j}^j. \end{equation}
By convention, we have that $x_{i,j}$ is square-free for $1 \leq i \leq 3$ and $1 \leq j \leq m_i  - 1$. Thus (\ref{htlb}) is equal to 
\[\prod_{j=1}^{m_1 - 1} x_{1,j}^{\frac{m_1 - j}{m_1}} \prod_{j=1}^{m_2 - 1} x_{2,j}^{\frac{m_2 - j}{m_2}} \prod_{j=1}^{m_3 - 1} x_{3,j}^{\frac{m_3 - j}{m_3}}.\]
Viewing these products as coefficients in (\ref{hgeq}), we see that if $H^0$ is to be bounded, that these coefficients must be bounded. Therefore, it suffices to check that for a fixed triple of integers $a,b,c$ that the equation
\[ax^{m_1} + b y^{m_2} + c z^{m_3} = 0\]
has finitely many primitive integer solutions when $1/m_1 + 1/m_2 + 1/m_3 < 1$. But this is exactly the content of Darmon and Granville's paper \cite{DG}, so we are done. 

\section{Northcott property of perturbed anti-canonical heights and the $abc$-conjecture} 

In this section we prove Theorem \ref{MT4}, starting with Theorem \ref{nearnorth}. We consider the property of recovering Northcott's property on a modified ESZ-B anti-canonical height on the stacky curve
\[\fX = \fX(\bP_{\bQ}^1 : (\Ba, \Bm)).\]
Here the modified height takes the shape
\[H_{(\Ba, \Bm)}^\delta (x,y) = \prod_{i=1}^n \phi_{m_i} (\ell_i(x,y))^{1/m_i} \max\{|x|, |y|\}^{\delta}.\]
Since $\chi(\fX)=\delta(\Bm)=2-\sum_{i=1}^n(1-\frac{1}{m_i})\leq 0$, our goal is to show that 

\[\gamma(\fX)=\inf\{\delta\in \bR\colon H^\delta_{(\Ba,\Bm)}\textnormal{ has the Northcott property}\}=\chi(\fX)\]
assuming the $abc$-conjecture. Recall that we have shown that $ H^{\chi(\fX)}_{(\Ba,\Bm)}$ does not have the Northcott property unconditionally. Thus we must show that $ H^{\chi(\fX)+\kappa}_{(\Ba,\Bm)}$ has the Northcott property for all $\kappa>0$. First assume that $\chi(\fX)=0$. The Northcott property for the standard height implies that $H_{(\Ba, \Bm)}^\delta (x,y)$ has the Northcott property whenever $\delta>0$. So $\inf\{\delta\in \bR\colon H^\delta_{(\Ba,\Bm)}\textnormal{ has the Northcott property}\}=0=\chi(\fX)$ as needed. Now suppose that $\chi(\fX)<0$. Assume, without loss of generality, that $m_1 \leq  m_2 \leq  \cdots \leq  m_n$. We then write
\[\ell_i(x,y) = z_{i,1} z_{i,2}^2 \cdots z_{i,m_i - 1}^{m_i-1} z_{i, m_i}^{m_i}.\]
We have

	\[\phi_{m_i}(\ell_i(x,y))^{\frac{1}{m_i}}=\prod_{j=1}^{m_i-1}z_{i,j}^{\frac{m_i-j}{m_i}}.\]	
It follows that

\[H^{\chi(\fX)+\delta}_{(\Ba,\Bm)}=\max\{\vert x\vert,\vert y\vert\}^{\chi(\fX)+\delta}\prod_{i=1}^n\prod_{j=1}^{m_i-1}z_{i,j}^{\frac{m_i-j}{m_i}}\]

Suppose that the following inequality holds for any $\epsilon>0$,

\begin{equation}\label{eq:abcbackbone}
	\prod_{i=1}^n\prod_{j=1}^{m_i-1}z_{i,j}^{\frac{m_i-j}{m_i}}\gg_\ep \max\{\vert x\vert,\vert y\vert\}^{-\chi(\fX)-\ep }.
\end{equation}
Then by multiplying both sides of the equation by $ \max\{\vert x\vert,\vert y\vert\}^{\chi(\fX)+\kappa}$ we obtain

\begin{equation}
	\max\{\vert x\vert,\vert y\vert\}^{\chi(\fX)+\kappa}\prod_{i=1}^n\prod_{j=1}^{m_i-1}z_{i,j}^{\frac{m_i-j}{m_i}}\gg_\ep \max\{\vert x\vert,\vert y\vert\}^{\kappa-\epsilon }.
\end{equation}
Taking $\ep=\frac{\kappa}{2}$ we have that

\begin{equation}
	H_{(\Ba, \Bm)}^{\chi(\fX)+\kappa} (x,y)\gg_\ep \max\{\vert x\vert,\vert y\vert\}^{\frac{\kappa}{2} }.
\end{equation}
Thus $H_{(\Ba, \Bm)}^{\chi(\fX)+\kappa} (x,y)$ must have the Northcott property as it cannot remain bounded by the usual Northcott property for $\bP^1$. Therefore we are done if we can confirm inequality (\ref{eq:abcbackbone}). To do so we require the following proposition, due to Granville \cite{Gran}:

\begin{proposition}[Granville] \label{gran prop} Suppose that the $abc$-conjecture holds. Then for any binary form $F$ with non-zero discriminant and $\ep > 0$ we have
\[\operatorname{rad}(F(m,n)) = \prod_{p | F(m,n)} p \gg_{F, \ep} \max\{|m|,|n|\}^{\deg F - 2 - \ep}.\]
\end{proposition}
In other words, if the $abc$-conjecture holds then the radical of $F(m,n)$ will be quite large compared to the variables $m,n$ (provided that the degree is at least 3). \\

We will apply Proposition \ref{gran prop} to reduce the proof of Theorem \ref{nearnorth} to a \emph{linear programming problem}. 

\subsection{A linear program bound} Observe that for each $1 \leq i \leq n$ that
\[\prod_{j=1}^{m_i} z_{i, j} \geq \text{rad}\left( \prod_{j=1}^{m_i} z_{i,j} \right).\]
Applying Proposition \ref{gran prop} to the binary form 
\[Q_\Ba(x,y) = \prod_{i=1}^n \ell_i(x,y) \] 
in conjunction with the above observation, we obtain:
\begin{equation}\label{eq:Granvillecontraint}
	\prod_{i=1}^n \prod_{j=1}^{m_i} z_{i,j} \geq \rad \left(\prod_{i=1}^n \ell_i(x,y) \right) \gg_\ep \max\{|x|, |y|\}^{n - 2 - \ep}
\end{equation}		
Similarly for each $i$ we have the bound
\begin{equation}\label{eq:basicconstraint}
|\ell_i(x,y)| \ll \max\{|x|, |y|\}.
\end{equation}
Taking logarithms and writing $y_{i,j} = \log |z_{i,j}|$, we then have an optimization problem:
\begin{equation}\label{eq:objectivefunction}
\min \sum_{i=1}^n \frac{1}{m_i} \sum_{j=1}^{m_i - 1} (m_i -j) y_{i,j}
\end{equation}
subject to
\begin{equation}\label{eq:logconstraint1}
\sum_{i=1}^n \sum_{j=1}^{m_i} y_{i,j} \geq  (n-2 -\ep) \log B
	\end{equation}
and
\begin{equation}\label{basic:constraint1}
\sum_{j=1}^{m_i} j y_{i,j} \ll  \log B
\end{equation}	
where $B = \max\{|x|, |y|\}$. Further, we have $y_{i,j} \geq 0$ for all $i,j$. \\ 

We emphasize that, at this point, integrality no longer plays a role, and neither does the syzygies relating the $z_{i,j}$'s. Indeed, we only need to solve the above linear program allowing arbitrary real inputs. \\

Now put
\[c_{ij}=\dfrac{m_i-j}{m_i}\]  
for $1\leq i\leq n$ and $1\leq j\leq m_i-1$. Write $c_{i,m_i}=0$ and let $\Bc=(c_{i,j})$ to be the column vector with

\[\Bc^{T}=[c_{1,1},c_{1,2},\dots,c_{1,m_1-1},0,c_{2,1},c_{2,2},\dots,c_{n,1},\dots c_{n,m_n-1},0].\]
We have that $\Bc\in\bR^{N}$ where $N=\sum_{i=1}^n n m_i$. \\

 Let $A$ be the matrix with rows representing the constraints,

\begin{align}
C_0&:\sum_{i=1}^n \sum_{j=1}^{m_i} y_{i,j} \geq  (n-2 -\ep) \log B\\
C_i&:-\sum_{j=1}^{m_i} j y_{i,j} \gg  -\log B
\end{align}
If we have take $\Be_{i,j}$ to be a basis of $\bR^N$ then we have that the the rows of $A$ are given  by

\begin{align}
\Ba_0&=\sum_{i=1}^n\sum_{j=1}^{m_i}\Be_{i,j}\\
\Ba_i&=\sum_{j=1}^{m_i}-j\Be_{i,j}
\end{align}
Finally let $\Bb$ be the column vector with $n+1$ entries representing the the constraints given by (\ref{eq:logconstraint1}) and (\ref{basic:constraint1}). In other words we have

\begin{align}
b_0&=(n-2-\epsilon)\log B\\
b_i&=-\log B
\end{align}

Our linear programming problem is then the following: let $\By=(y_{ij})$ ordered as above.

\begin{align}\label{eq:primalLP}
	&\textnormal{Minimize: }\ \ \ \ \ \ \ \ \ \ \ \ \ \ \ \ \ \ \ \ \ \ \ \ \ \ \ \ \ \  \Bc^{T}\By\\
	&\textnormal{subject to: }\ \ \ \ \ \ \ \ \ \ \ \ \ \ \ \ \ \ \ \ \ \ \ \ \ \ \ \ \ \  A\By\geq \Bb\textnormal{ and }\By\geq 0. \notag
\end{align}

The dual linear program is

\begin{align}\label{eq:dualLP}
	&\textnormal{Maximize: }\ \ \ \ \ \ \ \ \ \ \ \ \ \ \ \ \ \ \ \ \ \ \ \ \ \ \ \ \ \  \Bb^{T}\Bx\\
	&\textnormal{subject to: }\ \ \ \ \ \ \ \ \ \ \ \ \ \ \ \ \ \ \ \ \ \ \ \ \ \ \ \ \ \  A^T\Bx\leq \Bc\textnormal{ and }\Bx\geq 0. \notag
\end{align}
where $\Bx=[x_0,x_1,...,x_n]$. We call a vector $\Bx$ dual feasible if $A^T\Bx\leq \Bc$ and and  vector $\By$ primal feasible if $A\By\geq \Bb$. We have the following well known \textbf{weak duality statement}.

\begin{lemma}[Weak duality]\label{lem:weak duality}
Let $A$ be an $m\times n$ matrix with real entries and $\Bc$ a $n\times 1$ real vector and $\Bb$ an $m\times 1$ real vector. Consider the primal linear program
	
\begin{align*}
	&\textnormal{Minimize: }\ \ \ \ \ \ \ \ \ \ \ \ \ \ \ \ \ \ \ \ \ \ \ \ \ \ \ \ \ \  \Bc^{T}\By\\
	&\textnormal{subject to: }\ \ \ \ \ \ \ \ \ \ \ \ \ \ \ \ \ \ \ \ \ \ \ \ \ \ \ \ \ \  A\By\geq \Bb\textnormal{ and }\By\geq 0.
\end{align*}

and the dual linear program

\begin{align*}
	&\textnormal{Maximize: }\ \ \ \ \ \ \ \ \ \ \ \ \ \ \ \ \ \ \ \ \ \ \ \ \ \ \ \ \ \  \Bb^{T}\Bx\\
	&\textnormal{subject to: }\ \ \ \ \ \ \ \ \ \ \ \ \ \ \ \ \ \ \ \ \ \ \ \ \ \ \ \ \ \  A^T\Bx\leq \Bc\textnormal{ and }\Bx\geq 0.
\end{align*}	
	
Let $\By$ be any primal feasible vector and $\Bx$ a dual feasible vector. Then

\[\Bc^{T}\By\geq \Bb^{T}\Bx.\] 	

\end{lemma}
\begin{proof}
Let $A=(a_{i,j})$. Because $\By$ is primal feasible we have $A\By\geq \Bb$. Therefore for all $1\leq i\leq m$ we have

\[\sum_{j=1}^na_{i,j}y_j\geq b_i.\]

Multiplying by $x_i$ and summing over all $i$ we have

\begin{equation}\label{eq:weakduality1}
	\sum_{i=1}^m\sum_{j=1}^na_{i,j}y_jx_i\geq \sum_{i=1}^n b_ix_i=\Bb^{T}\Bx. 
	\end{equation}

On the other hand because $\Bx$ is dual feasible we have that $A^T\Bx\leq \Bc$ so for each $1\leq j\leq n$ we have

\[\sum_{i=1}^ma_{j,i}x_i\leq c_j.\]

Multiplying by $y_j$ and summing over all $j$ gives

\begin{equation}\label{eq:weakduality2}
\sum_{j=1}^n\sum_{i=1}^ma_{i,j}x_iy_j\leq \sum_{j=1}^ny_jc_j=\Bc^{T}\By.
\end{equation}

Combining inequality (\ref{eq:weakduality1}) and inequality (\ref{eq:weakduality2}) gives

\[\Bc^{\ T}\By\geq \sum_{j=1}^n\sum_{i=1}^ma_{i,j}x_iy_j\geq \Bb^{T}\Bx.  \]

\end{proof}

Returning to our problem, the weak duality theorem tells us that it suffices to find a dual feasible solution $\Bx=[x_0,\dots x_n]$ such that $\Bb^{T}\Bx\geq -\chi(\fX)+\epsilon$. In other words we seek $\Bx=[x_0,\dots x_n]$ with

\begin{align*}
\Bb^{T}\Bx&=\log B\left((n-2-\epsilon)x_1-\sum_{j=1}^nx_i\right)\geq -\chi(\fX)+\epsilon \\
A^T\Bx&\leq \Bc\\
\Bx&\geq 0
\end{align*}

Take $\Bx=[1,\frac{1}{m_1},\frac{1}{m_2},\dots ,\frac{1}{m_n}]$. We first show that $\Bx$ is dual feasible. In this case $A$ is a $(n+1)\times \sum_{i=1}^n m_i$ matrix. So a row of $A^T$ is indexed by a pair $(i,j)$ with  $1\leq i\leq n$ and $1\leq j\leq m_i$. We have that the $(i,j)$ entry of $A^T\Bx=\Bx^{T}A$ can be computed as 

\[x_0-jx_i.\]

Therefore to show that $\Bx$ is dual feasible for an arbitrary $\bX$ we need that

\[x_0-jx_i\leq c_{i,j}=\dfrac{m_i-j}{m_i}=1-\dfrac{j}{m_i}.\]
 In our case the $(i,j)$ entry of $A^T\Bx$ is given by

\[1-\dfrac{j}{m_i}=c_{i,j},\]
so $\Bx$ is dual feasible. We then compute

\begin{align*}
\Bb^{T}\Bx&=\log B(n-2-\ep)x_0-\sum_{i=1}^nx_i\log B\\
&=\log B\left(n-2-\ep-\sum_{i=1}^n\dfrac{1}{m_i}\right)\\
&=\log B\left( -(2-\sum_{i=1}^n(1-\dfrac{1}{m_i})-\epsilon )\right)\\
&=\log B(-\chi(\fX)-\epsilon .)
\end{align*}

Therefore $\Bx=[1,\dfrac{1}{m_1},\dots ,\dfrac{1}{m_n}]$ is a dual feasible solution and

\[\Bb^{T}\Bx=\log B(-\chi(\fX)-\epsilon). \]
By the weak duality theorem we have that
\[\sum_{i=1}^n \frac{1}{m_i} \sum_{j=1}^{m_i - 1} (m_i -j) y_{i,j}\geq \log B(-\chi(\fX)-\epsilon). \]
Exponentiating gives

\[\prod_{i=1}^n\prod_{j=1}^{m_i-1}z_{i,j}^{\frac{m_i-j}{m_i}}\geq B^{(-\chi(\fX)-\epsilon)}. \]

As $B=\max\{\vert x\vert ,\vert y\vert\}$ we have verified inequality \ref{eq:abcbackbone} and consequently we have that conditional on the $abc$-conjecture that

\[\gamma(\fX)=\chi(\fX)\]
when $\chi(\fX)\leq 0$. 

\subsection{Proof of Theorem \ref{MT4}} 

One direction of the theorem is provided by Theorem \ref{nearnorth}, which we proved in the previous subsection. It suffices to prove the converse. \\ \\
Actually, for the convese we only need the assertion that for any $\kappa > 0$ and $m \geq 4$ that the function $H_{-K_{\fX_m}}(\Bx) \cdot H(\Bx)^\kappa$ has Northcott's property, where $\fX_m = \fX(\bP^1 : ((0, 1, \infty), (m, m, m))$. To see this, let us fix $\ep > 0$. Choose $0 < \kappa < \ep/3$ and choose $m \in \bN$ sufficiently large so that 
\[\frac{m-3}{m-1} + \frac{\kappa m}{2(m-1)} > \frac{1}{1 + \ep} .\]
By hypothesis, we have
\begin{equation} \label{abc1} H_{-K_{\fX_m}}(\Bx) = \phi_m(x)^{1/m} \phi_m(y)^{1/m} \phi_m(x+y)^{1/m}\max\{\vert x\vert,\vert y\vert\}^{\frac{3}{m}-1} \gg_\kappa \max\{|x|, |y|\}^{1 - \frac{3}{m}  - \frac{\kappa}{2} } \end{equation}
Trivially, we see that
\[\phi_m(u) \leq \rad(u)^{m-1}\]
for all $u \in \bZ$. Hence (\ref{abc1}) implies
\begin{equation} \rad(x)^{\frac{m-1}{m}} \rad(y)^{\frac{m-1}{m}} \rad(x+y)^{\frac{m-1}{m}} \gg_\kappa \max\{|x|, |y|\}^{1 - \frac{3}{m}  + \frac{\kappa}{2}}.
\end{equation}
Since $x,y,x+y$ are pairwise co-prime we have $\rad(x)\rad(y)\rad(x+y) = \rad(xy(x+y))$, hence
\[\rad(xy(x+y))^{\frac{m-1}{m}} \gg_\kappa \max\{|x|,|y|\}^{1 - \frac{3}{m}  + \frac{\kappa}{2}}.\]
Raising both sides to the $m/(m-1)$ power, we have
\[\rad(xy(x+y)) \gg_\kappa \max\{|x|,|y|\}^{\frac{m-3}{m-1} + \frac{\kappa m}{2(m-1)}} \gg_\kappa \max\{|x|, |y|\}^{\frac{1}{1+\ep}},\]
by out hypotheses on $m, \kappa$. It follows that
\[\rad(xy(x+y))^{1 + \ep} \gg_\ep \max\{|x|,|y|\}\]
which is plainly equivalent to the $abc$-conjecture, provided we adjust the implied constant.

\section{Quantitative arithmetic of stacky curves}
\label{count sec}

\subsection{Crude bound for $N_\Bm(T)$, with $\Bm = (2,2,m)$}

Here we deal with the case $\Bm = (2,2,m)$. The Euler characteristic is equal to 
\[\delta(\Bm) = 2 - \frac{1}{2} - \frac{1}{2} - 1 + \frac{1}{m} = \frac{1}{m}.\]
The height $H(x,y)$ is given by
\[H(x,y) = |x_1|^{1/2} |y_1|^{1/2} |z_1^{m-1} \cdots z_{m-1}|^{1/m} \max\{|x_1 x_2^2|, |y_1 y_2^2|\}^{1/m}, \]
where $x = x_1 x_2^2, y = y_1 y_2^2$ and
\begin{equation} \label{keyeq2} x_1 x_2^2 + y_1 y_2^2 = z_1 z_2^2 \cdots z_{m-1}^{m-1} z_m^m,\end{equation} 
with $x_1, y_1, z_1, \cdots, z_{m-1}$ square-free. We normalize the height by raising it to the $m$-th power, obtaining the bound
\begin{equation} \label{genhtbd}  |x_1 y_1|^{m/2} |z_1^{m-1} \cdots z_{m-1}| \max\{|x_1 x_2^2|, |y_1 y_2^2|\} \leq T.
\end{equation}
From here we see that
\[|z_1 z_2^2 \cdots z_m^m| \ll \max\{|x_1 x_2^2|, |y_1 y_2^2|\} \ll \frac{T}{|x_1 y_1|^{m/2} |z_1^{m-1} \cdots z_{m-1}|}, \]
whence we conclude that
\[|z_m| \ll \frac{T^{1/m}}{|x_1 y_1|^{1/2} |z_1 \cdots z_{m-1}|}.\]
This bound and $|z_m| \geq 1$ implies that 
\[|x_1 y_1|^{1/2} |z_1 \cdots z_{m-1}| \ll T^{1/m}.\]
From here we obtain a crude upper bound for $N_\Bm(T)$, which proves Theorem \ref{MT0}. Indeed, having chosen $x_1, y_1, z_1, \cdots, z_{m-1}$ there are then $O(T^{1/m}/(|x_1 y_1|^{1/2} |z_1 \cdots z_{m-1}|)$ possibilities for $z_m$. Having chose $z_m$ as well, there are then $O_\ep (T^\ep)$ possibilities for $x_2, y_2$, since $x_2, y_2$ are polynomially bounded so they will be determined by the norm-equation (\ref{keyeq2}) up to a log factor. Thus, there are
\begin{equation} \sum_{|x_1 y_1|^{1/2} |z_1 \cdots z_{m-1}| \leq T^{1/m}} O_\ep \left(\frac{T^{1/m + \ep}}{|x_1 y_1|^{1/2} |z_1 \cdots z_{m-1}|} \right)
\end{equation}
possible solutions to (\ref{keyeq2}) satisfying the height bound (\ref{genhtbd}).  We evaluate this as
\begin{align*} & \sum_{|x_1 y_1|^{1/2} |z_1 \cdots z_{m-1}| \leq T^{1/m}} O_\ep \left(\frac{T^{1/m + \ep}}{|x_1 y_1|^{1/2} |z_1 \cdots z_{m-1}|} \right) \\ 
& = \sum_{|x_1 y_1| \leq T^{2/m}} \frac{1}{|x_1 y_1|^{1/2}} \sum_{|z_1 \cdots z_{m-1}| \leq T^{1/m}/|x_1 y_1|^{1/2}} O_\ep \left(\frac{T^{1/m + \ep}}{|z_1 \cdots z_{m-1}|} \right) \\
& \ll_\ep \sum_{|x_1 y_1| \leq T^{2/m}} \frac{T^{1/m + \ep}}{|x_1 y_1|^{1/2}} \ll_\ep T^{2/m + \ep}. 
\end{align*}

To give a lower bound, we choose square-free integers $a,b,c$ so that the curve 
\[ax^2 + by^2 = cz^m\]
has a primitive integral solution. Such a triple is guaranteed to exist; see \cite{Beu}. Then we can parametrize (some) of the solutions by a triple of integral binary forms $(F,G, h)$ where $\deg F = \deg G = m$ and $\deg h = 2$. By 
\[x = F(u,v), y = G(u,v), z = h(u,v).\]
The height is 
\[|a|^{m/2} |b|^{m/2} |c|^{m-1} \max\{|ax^2|, |by^2|\},\]
so if we treat $a,b,c$ as constants then 
\[\max\{|x|, |y|\} \ll_{a,b,c} T^{1/2}.\]
Therefore, we are looking for solutions to the Thue inequality
\[\max\{|F(u,v)|, |G(u,v)|\} \ll_{a,b,c} T^{1/2}.\]
If we restrict $u,v$ so that
\[\max\{|u|, |v|\} \ll_{a,b,c} T^{1/(2m)},\]
then we see that the above height bound is satisfied. Thus $N_m(T) \gg T^{1/m}$.

\subsection{Proof of Theorem \ref{MT}} 

In this section, we prove Theorem \ref{MT}. To do so we will show that $N_2(T) = O \left(T^{1/2} (\log T)^3 \right)$ and give a separate argument to show that $N_2(T) \gg T^{1/2} (\log T)^3$. The incompatibility of these two arguments represents the main obstacle as to why an asymptotic formula for $N_2(T)$ remains elusive. \\ 

We count rational points of bounded height on the curve $\fX(\bP^1_\bQ;(0,2),(-1,2),(\infty,2))$ with the height on $\bP^1$ given by (\ref{eq:HeightDef}). On writing
\[a = x_1 y_1^2, b = x_2 y_2^2, x_1, x_2 \text{ square-free}\]
(note that this differs from the notation used elsewhere in the paper) we then have
\[H(a,b) = |x_1 x_2| \sqf(x_1 y_1^2 + x_2 y_2^2) \max\{|x_1 y_1^2|,  |x_2 y_2^2|\}.\]
and the max on the right hand side is dependent only on the relative size of $|a|,|b|$. If we write
\begin{equation} \label{variety} x_1 y_1^2 + x_2 y_2^2 = x_3 y_3^2,\end{equation} 
then we further obtain the expression
\[H(a,b) = \max\{|x_2 x_3 (x_1 y_1)^2|, |x_1 x_3 (x_2 y_2)^2|\}.\]
We may assume without loss of generality that $|x_1 y_1^2| \geq |x_2 y_2|^2$ and $x_1 > 0$, so that 
\[H(a,b) = |x_2 x_3 (x_1 y_1)^2|.\]

We consider the problem of counting integral points on the variety defined by (\ref{variety}), subject to the constraint
\begin{equation} \label{ht bd} 0 < |x_2 x_3 (x_1 y_1)^2| \leq T, |x_1 y_1^2| \geq |x_2 y_2^2|.  \end{equation} 

To obtain the upper bound we must dissect (\ref{ht bd}) into suitable ranges. When $|x_1 x_2 x_3| \leq T^{1/2}$ we fix $x_1, x_2, x_3$ and treat (\ref{variety}) as a diagonal ternary quadratic form, say $Q_\Bx$. It is then the case that 
\begin{equation} \label{yi bd} |y_i| \leq \frac{T}{|x_1 x_2 x_3| \cdot |x_i|}\end{equation} 
for $i = 1,2,3$, and by Corollary 2 of \cite{Brow-HB} we then have the estimate
\[O \left(d(x_1 x_2 x_3) \left(\frac{T^{1/2}}{|x_1 x_2 x_3|} + O(1) \right) \right) \]
for the number of $\By \in \bZ_{\ne 0}^3$ satisfying (\ref{ht bd}) and (\ref{variety}) provided that the quadratic form $Q_\Bx$ has a rational zero. Otherwise it is clear that there will be no contribution. Thus we must estimate 
\[\sum_{\substack{1 \leq |x_1 x_2 x_3| \leq T^{1/2} \\ Q_\Bx \text{ has a rational zero}}} d(x_1 x_2 x_3).  \]
This is similar to the work of Guo in \cite{Guo}, except he counted with respect to the height $\lVert \Bx \rVert_\infty$. Nevertheless the techniques are similar, and again this may be of independent interest. \\

Next we must deal with the case when $|x_1 x_2 x_3| \geq T^{1/2}$. For this it suffices to observe from (\ref{yi bd}) that $|x_1 x_2 x_3| \geq T^{1/2}$ implies 
\[|y_1 y_2 y_3| \leq \frac{T^{3/2}}{(x_1 x_2 x_3)^2} \leq T^{1/2}.\]
We then treat (\ref{variety}) as a linear form $L_\By$ in $\Bx$. We use this to show that the contribution for each $\By$ is $O\left(T^{1/2} |y_1 y_2 y_3|^{-1} + 1 \right)$, which gives an acceptable contribution upon summing over $\By$. \\

For the lower bound, we first restrict $y_1, y_2, y_3 \in \bZ_{\ne 0}$ satisfying 
\[|y_1 y_2 y_3| \leq T^\delta \]
for some explicit $\delta > 0$ to be specified later. We note that to obtain the correct order of magnitude it is permissible to choose any $\delta > 0$. \\

Having fixed $\By = (y_1, y_2, y_3)$, we consider the simultaneous conditions (\ref{variety}) and (\ref{ht bd}). This gives rise to a binary form inequality of the shape 
\begin{equation} \label{bin form in} |x_1^2 x_2 (y_1^2 x_1 + y_2^2 x_2)| \leq T y_3^2 y_1^{-2}. 
\end{equation}
Because $|y_1 y_2 y_3|$ is small, we can count the number of solutions $\Bx$ to this inequality with reasonable precision. However, even with $|y_1 y_2 y_3|$ counting the number of solutions $\Bx$ with enough uniformity appears to still be a challenging task, because the binary form in (\ref{bin form in}) is singular. This difficulty is exacerbated by the fact that we will need to apply a square-free sieve eventually to produce triples $\Bx$ with each coordinate square-free. \\

To get around this issue, we simply count solutions to (\ref{bin form in}) with $x_1, x_2$ satisfying the inequalities 
\[|x_i y_i^2| \leq c_i T^{1/4} |y_1 y_2 y_3|^{1/2}, i = 1,2\]
for some positive numbers $c_1, c_2$. This has the effect that the long cusps inherent in (\ref{bin form in}) are removed, and reduces the problem to a more straightforward geometry of numbers question. 

\subsubsection{Upper bounds}
\label{crude}

To obtain upper bounds, it is crucial to view (\ref{variety}) as a plane in $x_1, x_2, x_3$  when $|y_1 y_2 y_3| \leq T^{1/2}$ and viewing (\ref{variety}) as a conic in $y_1, y_2, y_3)$ when $|x_1 x_2 x_3| \leq T^{1/2}$. We call the former the \emph{linear case} and the latter the \emph{quadratic case}. We proceed to deal with the linear case below. \\

We shall first suppose that $|y_1 y_2 y_3| \leq T^{1/2}$ is fixed, and count the triples $(x_1, x_2, x_3)$ and $(y_1, y_2, y_3)$ for which (\ref{variety}) holds. \\

The key is the following lemma on counting points in sublattices of $\bZ^2$: 

\begin{lemma} \label{lat lem} Let $\Lambda \subset \bZ^2$ be a lattice. Then for all positive real numbers $R_1, R_2$ the number of primitive integral points $\Bx \in \Lambda$ satisfying $|x_i| \leq R_i, i = 1,2$ is at most $O \left(R_1 R_2/\det(\Lambda) + 1 \right)$. 
\end{lemma}

\begin{proof} If the rectangle $[-R_1, R_2] \times [-R_2, R_2]$ contains at least two primitive vectors in $\Lambda$, say $\Bx_1, \Bx_2$, then since this rectangle is convex it contains the parallelogram with end points $\pm \Bx_1, \pm \Bx_2$. The area of this parallelogram is at least as large as $\det \Lambda$, since the lattice spanned by $\Bx_1, \Bx_2$ is a sublattice of $\Lambda$. It thus follows that
\[R_1 R_2 \gg \det \Lambda.\]
Otherwise, the rectangle $[-R_1, R_1] \times [-R_2, R_2]$ contains at most one primitive vector in $\Lambda$. This completes the proof. 
\end{proof}
The strength of this lemma is that it gives a strong upper bound even in lopsided boxes. \\

Given (\ref{variety}), it follows that there is at least one $i \in \{2,3\}$ such that
\[|x_i y_i^2|/2  \leq x_1 y_1^2 \leq 2|x_i y_i^2|, \]
whence
\[\frac{x_1 y_1^2}{2 y_i^{-2}} \leq  |x_i| \leq \frac{2 x_1 y_1^2}{ y_i^{2}}.\]
Without loss of generality, we assume that this holds for $i = 2$. Suppose that $M_1 \leq x_1 < 2M_1$. By (\ref{ht bd}), we have
\[|x_3| \leq \frac{T}{|x_2 x_1^2 y_1^2|},\]
whence
\begin{align*} |x_3| & \leq T \cdot \frac{2 y_2^2}{(x_1 y_1^2)(x_1^2 y_1^2)}  \\
& \leq \frac{2y_2^2 T}{M_1^3 y_1^4}
\end{align*}
Applying Lemma \ref{lat lem} to the lattice defined by the congruence $y_1^2 x_1 - y_3^2 x_3 \equiv 0 \pmod{y_2^2}$ which has determinant equal to $y_2^2$, there are 
\[O \left(M_1 \cdot \frac{T y_2^2}{M_1^3 y_1^4} \cdot \frac{1}{y_2^2} + 1 \right) = O \left(\frac{T}{M_1^2 y_1^4} + 1 \right)\]
possibilities for $x_1, x_3$, which then determines $x_2 = (y_1^2 x_1 - y_3^2 x_3)/y_2^2$. Similarly, applying Lemma \ref{lat lem} to the lattice defined by $y_1^2 x_1 + y_2^2 x_2 \equiv 0 \pmod{y_3^2}$, with determinant equal to $y_3^2$, gives the estimate
\[O \left(M_1 \cdot \frac{y_1^2 M_1}{y_2^2} \frac{1}{y_3^2} + 1\right) = O \left(\frac{M_1^2 y_1^2}{y_2^2 y_3^2} + 1 \right) \]
for the number of $x_1, x_2$, which then also determine $x_3$. The two bounds coincide when 
\[M_1 = \frac{T^{1/4} |y_2 y_3|^{1/2}}{|y_1|^{3/2}}, \]
and we get the bound
\[O \left(\frac{T^{1/2} |y_2 y_3| y_1^2}{y_2^2 y_3^2 |y_1|^3} + 1 \right) = O \left(\frac{T^{1/2}}{|y_1 y_2 y_3|} + 1\right) \]
for the number of $x_1, x_2, x_3$ given $y_1, y_2, y_3$. Thus, we obtain an acceptable estimate whenever $|y_1 y_2 y_3| \ll T^{1/2}$, since
\begin{align*} \sum_{1 \leq |y_1 y_2 y_3| \leq T^{1/2}} \frac{T^{1/2}}{|y_1 y_2 y_3|} + 1 & \ll T^{1/2} \sum_{n \leq T^{1/2}} \frac{d_3(n)}{n} + \sum_{n \leq T^{1/2}} d_3(n)
\end{align*}
It is well-known that
\[\sum_{n \leq Z} d_3(n) = Z (\log Z)^2 + O(Z \log Z). \]
By partial summation, we have
\begin{align*} \sum_{n \leq Z} \frac{d_3(n)}{n} & = Z^{-1} \sum_{n \leq Z} d_3(n) + \int_1^Z \left(\sum_{n \leq t} d_3(n) \right) \frac{dt}{t^2} \\
& \ll (\log Z)^2 + \int_1^Z \frac{(\log t)^2 dt}{t} \\
& \ll (\log Z)^3
\end{align*}
It follows that 
\[T^{1/2} \sum_{n \leq T^{1/2}} \frac{d_3(n)}{n} + \sum_{n \leq T^{1/2}} d_3(n) \ll T^{1/2} (\log T)^3.\]

It remains to deal with the case when $|y_1 y_2 y_3| \gg T^{1/2}$, where we instead fibre over $\Bx$ and consider zeroes of the corresponding diagonal quadratic forms $Q_\Bx$. Since
\[|x_i y_i^2| \ll x_1 y_1^2 \]
for $i = 1,2$ by assumption, it follows that 
\[|x_1 x_2 x_3 y_1^2 y_2^2 y_3^2| \leq x_1^3 y_1^6,\]
hence
\[|y_1^2 y_2^2 y_3^2| \ll \frac{x_1^3 y_1^6}{x_1 |x_2 x_3|}.\]
If $|x_1 x_2 x_3| \gg T^{1/2}$, then 
\[x_1^3 y_1^6 \gg T^{3/2} \Leftrightarrow x_1 y_1^2 \gg T^{1/2}.\]
This implies that
\[|x_1 x_2 x_3| \cdot x_1 y_1^2 \gg T,\]
which violates (\ref{ht bd}) if the implied constants are sufficiently large. It thus follows that we must have $|x_1 x_2 x_3| \ll T^{1/2}$ in this case. \\

We now fix $x_1, x_2, x_3$ and consider (\ref{variety}) as a ternary quadratic form in $y_1, y_2, y_3$. We shall require the following version of Corollary 2 in \cite{Brow-HB}, which is an analogue of Lemma \ref{lat lem}: 

\begin{lemma} \label{quad lem} Let $x_1, x_2, x_3$ be pairwise co-prime square-free integers. Let $R_1, R_2, R_3$ be positive real numbers. Then the number of primitive solutions $y_1, y_2, y_3$ to the equation 
\[x_1 y_1^2 + x_2 y_2^2 = x_3 y_3^2\]
with $|y_i| \leq R_i$ is bounded by 
\[O \left( d(x_1 x_2 x_3) \left( \left(\frac{R_1 R_2 R_3}{|x_1 x_2 x_3|} \right)^{1/3} + 1 \right) \right).\]
\end{lemma}

Since $|x_i y_i^2| \ll x_1 y_1^2$ for $i = 1,2$, it follows that
\[|x_1 x_2 x_3 (x_i y_i^2)| \ll |x_1 x_2 x_3 (x_1 y_1^2)| \leq T\]
for $i = 1,2$, whence
\[|(x_1 y_1)^2 x_2 x_3|, |(x_2 y_2)^2 x_1 x_3|, |(x_3 y_3)^2 x_1 x_2| \ll T.\]
This implies that
\[(y_1 y_2 y_3)^2 (x_1 x_2 x_3)^4 \ll T^3,\]
hence
\[|y_1 y_2 y_3| \ll \frac{T^{3/2}}{(x_1 x_2 x_3)^2}.\]
Lemma \ref{quad lem} then implies that for fixed $x_1, x_2, x_3$ the number of primitive $\By = (y_1, y_2, y_3)$ satisfying (\ref{variety}) is 
\[O \left(d(x_1 x_2 x_3) \left(\frac{T^{1/2}}{|x_1 x_2 x_3|} + 1 \right) \right). \]
We now sum over primitive $\Bx \in \bZ^3$ satisfying $|x_1 x_2 x_3| \ll T^{1/2}$, with the property that the quadratic form $Q_\Bx$ given by (\ref{variety}) has a rational zero. By the Hasse-Minkowski theorem, this is tantamount to the form $Q_\Bx(\By) = x_1 y_1^2 + x_2 y_2^2 - x_3 y_3^2$ being everywhere locally soluble. The estimation of this is interesting on its own right and will be handled in a separate subsection.

\subsubsection{Counting soluble ternary quadratic forms} 

In this section, we consider the set
\[\S = \{(x_1, x_2, x_3) \in \bZ^3 : x_1, x_2, x_3 > 0, \gcd(x_1, x_2) = \gcd(x_1, x_3) = \gcd(x_2, x_3) = 1, \]
\[x_i \text{ square-free for } i = 1,2,3, x_1 y_1^2 + x_2 y_2^2 - x_3 y_3^2 \text{ is everywhere locally soluble}\}.\]
By a well-known theorem of Legendre (see \cite{Guo}) the indicator function for $\S$ is given by
\begin{equation} f_\S(x_1, x_2, x_3) = \left(2^{-\omega(x_1)} \sum_{a_1 | x_1} \left(\frac{x_2 x_3}{a_1} \right) \right) \left(2^{-\omega(x_2)} \sum_{a_2 | x_2} \left(\frac{x_1 x_3}{a_2} \right) \right) \left(2^{-\omega(x_3)} \sum_{a_3 | x_3} \left(\frac{-x_1 x_2}{a_3} \right) \right).
\end{equation}
We will now combine the ideas given in \cite{Guo} and those in \cite{FK}. \\

Put 
\begin{align*} \S(X) & = \sum_{1 \leq x_1 x_2 x_3 \leq X}  \sum_{(x_1, x_2, x_3) \in \S} \frac{d(x_1 x_2 x_3)}{x_1 x_2 x_3} \\ & = \sum_{1 \leq |x_1 x_2 x_3| \leq X} \frac{d(x_1 x_2 x_3)}{x_1 x_2 x_3}  f_\S(x_1, x_2, x_3). \end{align*}
Since $x_1, x_2, x_3$ are pairwise coprime and square-free, it follows that
\[d(x_1 x_2 x_3) = 2^{\omega(x_1 x_2 x_3)} = 2^{\omega(x_1)} \cdot 2^{\omega(x_2)} \cdot 2^{\omega(x_3)}, \]
where $\omega(n)$ is the number of distinct prime factors of $n$. It follows that
\begin{equation} \S(X) = \sum_{1 \leq x_1 x_2 x_3 \leq X} \frac{2^{\omega(x_1 x_2 x_3)}}{x_1 x_2 x_3} f_\S(x_1, x_2, x_3) 
\end{equation}
\[ = \sum_{1 \leq x_1 x_2 x_3 \leq X} \frac{1}{x_1 x_2 x_3} \left(1 + \left(\frac{x_2 x_3}{x_1} \right) \left(\frac{x_1 x_3}{x_2} \right) \left(\frac{-x_1 x_2}{x_3} \right)  + \sum_g g(x_1, x_2, x_3)\right),\]
where $g$ expresses a product of Jacobi symbols. The sum 
\begin{equation} \label{S1X} \S_1(X) = \sum_{1 \leq |x_1 x_2 x_3| \leq X} \frac{1}{x_1 x_2 x_3} \left(1 + \left(\frac{x_2 x_3}{x_1} \right) \left(\frac{x_1 x_3}{x_2} \right) \left(\frac{-x_1 x_2}{x_3} \right)  \right) \end{equation}
is expected to contribute the main term while the sum 
\begin{equation} \label{char sum} \S_2(x) = \sum_{1 \leq x_1 x_2 x_3 \leq X} \frac{1}{x_1 x_2 x_3} \sum_g g(x_1, x_2, x_3) \end{equation}
is expected to be negligible, due to the cancellation of characters. \\

By partial summation, we obtain:
\begin{equation} \label{par sum} \S_i(X) = \frac{1}{X} \Sigma_i(X) + \int_1^X \Sigma_i(t) \frac{t}{t^2},
\end{equation}
where
\[\Sigma_1(X) = \sum_{\substack{1 \leq |x_1 x_2 x_3| \leq X \\ x_1 x_2 x_3 \text{ square-free}  \\ Q_{(x_1, x_2, x_3)} \text{ is soluble} }} \left(1 + \left(\frac{x_2 x_3}{x_1} \right) \left(\frac{x_1 x_3}{x_2} \right) \left(\frac{x_1 x_2}{x_3} \right) \right) \]
and
\[\Sigma_2(X) = \sum_{1 \leq x_1 x_2 x_3 \leq X} \sum_g g(x_1, x_2, x_3). \]

Our situation differs from that of Guo in \cite{Guo} since we are counting over triples with $|x_1 x_2 x_3| \leq X$ rather than $\max\{|x_1|, |x_2|, |x_3|\} \leq X$, which introduces some difficulties. However, this is exactly analogous to the situation encountered by Fouvry and Kluners in \cite{FK}. \\

Our key proposition will be:  

\begin{proposition} \label{quad form prop} We have the asymptotic upper bound
\[\S(X) = O\left((\log X)^3 \right).\]
\end{proposition} 
In fact we can refine Proposition \ref{quad form prop} to give an asymptotic formula, but this is unnecessary for our purposes. \\

We 

We proceed to prove Proposition \ref{quad form prop} in the remainder of the section. We begin by showing that triples $(x_1, x_2, x_3)$ with $\mu^2(x_1 x_2 x_3) = 1$ and $\omega(x_1 x_2 x_3)$ large contribute negligibly. To wit, put
\[\S_2^{(r)}(X) = \sum_{\substack{1 \leq x_1 x_2 x_3 \leq X \\ \omega(x_1 x_2 x_3) = r}} \frac{1}{x_1 x_2 x_3} \sum_g g(x_1, x_2, x_3).\]
By the triangle inequality, it is clear that 
\[\left \lvert \S_2^{(r)}(X) \right \rvert \ll \sum_{\substack{n \leq X \\ \mu^2(n) = 1, \omega(n) = r}} \frac{d_3(n)}{n}. \]
By partial summation, we have
\[\sum_{\substack{n \leq X \\ \mu^2(n) = 1, \omega(n) = r}} \frac{d_3(n)}{n} = X^{-1} \sum_{\substack{n \leq X \\ \mu^2(n) = 1, \omega(n) = r}} d_3(n) + \int_1^X \left(\sum_{\substack{n \leq t \\ \mu^2(n) = 1, \omega(n) = r}} d_3(n) \right)\frac{dt}{t^2}. \]
To estimate the latter sum, we will need the following result, which is Lemma 11 in \cite{FK}:

\begin{lemma} There exists an absolute constant $B_0 \geq 1$ such that for every $r \geq 0$, we have
\[|\{n \leq X : \omega(n) = r, \mu^2(n) = 1\}| \leq B_0 \cdot \frac{X}{\log X} \cdot \frac{(\log \log X + B_0)^r}{r!} \]
\end{lemma}

Applying the lemma, we have for $\Omega = 30 (\log \log X + B_0)$
\begin{align*} \sum_{\substack{n \leq X \\ \mu^2(n) = 1, \omega(n) \geq \Omega}} d_3(n) & \ll \frac{X}{\log X} \sum_{r \geq \Omega} 3^r \cdot \frac{(\log \log X + B_0)^r}{r!} \\
& \ll \frac{X}{\log X} \sum_{r \geq \Omega} \left(\frac{3e (\log \log X + B_0)}{r} \right)^r \\
& \ll \frac{X}{\log X} \sum_{r \geq \Omega} \left(\frac{3e}{10} \right)^r,
\end{align*}
the final sum a convergent geometric series. Hence
\[\sum_{r \geq \Omega} \left(\frac{3e}{10} \right)^r \ll \left(\frac{3e}{10} \right)^{\Omega} \ll \frac{1}{\log X}. \]
We thus conclude that 
\begin{align} \sum_{r \geq \Omega} \left \lvert S_2^{(r)}(X) \right \rvert & \ll 1 + (\log X)^{-2} + \int_1^X \frac{dt}{t (\log t)^2}  \\
& = O(1) \notag 
\end{align}
and is thus negligible. \\

Note that $x_1, x_2, -x_3$ cannot all be the same sign, otherwise (\ref{variety}) will only have a trivial real solution. Hence the signs of $(x_1, x_2, x_3)$ must be $(+, +, +)$, or $(+, -, +)$, since we assumed $x_1 > 0$ and $x_1 y_1^2 \geq |x_2 y_2^2|$. By rearranging, we must thus assume $x_1, x_2, x_3 > 0$. \\

We then expand (\ref{char sum}) by writing $x_i = x_{i1}x_{i2}$ for $i = 1,2,3$, and
\[\sum_{\substack{1 \leq x_1 x_2 x_3 \leq X \\ \mu^2(x_1 x_2 x_3) = 1}} \sum_g g(x_1, x_2, x_3) = \sum_{\substack{(x_{11}x_{12})(x_{21}x_{22})(x_{31}x_{32})\leq X \\ 1 < x_{i1} < x_i \text{ for } 1 \leq i \leq 3}} \left(\frac{x_{21}x_{22} x_{31}x_{32} }{x_{11}} \right) \left( \frac{x_{11} x_{12} x_{31} x_{32}}{x_{21}} \right) \left(\frac{x_{11} x_{12} x_{21} x_{22}}{x_{31}} \right). \]
We now follow the strategy outlined in \cite{FK} and break up the set 
\[\{(x_{11}, x_{12}, x_{21}, x_{22}, x_{31}, x_{32}) \in \bN^6 : x_{11} x_{12} x_{21} x_{22} x_{31} x_{32} \leq X \} \]
by restricting the $x_{ij}$'s to intervals of the form 
\[[A_{ij}, \Delta A_{ij}),\]
where 
\[\Delta = 1 + (\log X)^{-3}.\] 
For a given $\BA = (A_{11}, A_{12}, A_{21}, A_{22}, A_{31}, A_{32})$, put 
\[\S_2(X; \BA) = \sum_{\substack{x_{ij} \in [A_{ij}, \Delta A_{ij}) \\ \mu^2(x_{11} x_{12} x_{21} x_{22} x_{31} x_{32}) = 1 \\ \prod_{i,j} x_{ij} \leq X }} \left(\frac{x_{21}x_{22} x_{31}x_{32} }{x_{11}} \right) \left( \frac{x_{11} x_{12} x_{31} x_{32}}{x_{21}} \right) \left(\frac{x_{11} x_{12} x_{21} x_{22}}{x_{31}} \right). \]

We then have the following lemma:

\begin{lemma} We have the bound
\[\sum_{\prod A_{ij} \geq \Delta^{-6} X} \left \lvert \S_2(X; \BA) \right \rvert = O \left(X (\log X)^{-1} \right). \]
\end{lemma}

\begin{proof} We have
\begin{align*} \sum_{\prod A_{ij} \geq \Delta^{-6} X} \left \lvert \S_2(X; \BA) \right \rvert & \leq \sum_{\substack{\Delta^{-6} X \leq n \leq X \\ \mu^2(n) = 1}} d_3(n) \\
& \ll \sum_{\Delta^{-6} X \leq n \leq X} 3^{\omega(n)} \\
& \ll (1 - \Delta^{-6} ) X (\log X)^2. 
\end{align*}
By Taylor's theorem, we have
\[\Delta^{-6} = (1 + (\log X)^{-3})^{-6} = 1 - 6 (\log X)^{-3} + O \left((\log X)^{-6} \right).\]
The proof then follows.
\end{proof}

To proceed, we shall require the following well-known lemma regarding character sums: 

\begin{lemma}[Double Oscillation Lemma] \label{double} Let $\{\alpha_n\}, \{\beta_m\}$ be two sequences of complex numbers with each term having absolute value bounded by $1$. Let $M,N$ be positive real numbers. Then  we have
\[\sum_{m \leq M} \sum_{n \leq N} \alpha_m \beta_n \mu^2(2m) \mu^2(2n) \left(\frac{m}{n} \right) \] 
\[ \ll \min \left\{ \left(M^{-1/2} + (N/M)^{-1/2}  \right), \left(N^{-1/2} + (M/N)^{-1/2} \right) \right\}  \]
and for every $\ep > 0$, 
\[\sum_{m \leq M} \sum_{n \leq N} \alpha_m \beta_n \mu^2(2m) \mu^2(2n) \ll_\ep MN \left(M^{-1/2} + N^{-1/2} \right) (MN)^\ep \] 
\end{lemma}

We will also need the following variant of the Siegel-Walfisz theorem: 

\begin{lemma} \label{S-W} Let $\chi_q$ be a primitive character modulo $q \geq 2$. Then for every $A > 1$ we have
\[\sum_{Y \leq p \leq X} \chi_q(p) = O_A \left(\sqrt{q} \cdot X (\log X)^{-A} \right) \]
uniformly for $X \geq Y \geq 2$. 
\end{lemma}

We now consider, as in \cite{FK}, the quantities
\begin{equation} X^\dagger = (\log X)^9, X^\ddagger = \exp\left( (\log X)^{1/8} \right).
\end{equation}
We now consider those $\BA$ with the property that at most $2$ entries larger than $X^\ddagger$. We dissect the sum according to the number $r \leq 2$ of terms $A_{ij}$ greater than $X^\ddagger$. Let $n$ be the product of those $x_{ij}$ which are larger than $X^\ddagger$, and $m$ the product of the remaining ones. We sum over $\BA$ with such properties to obtain
\begin{align*} \sideset{}{^{(2)}} \sum_{\BA}  |\S_2(X; \BA)| & \leq \sum_{r \leq 2} \sum_{m \leq (X^\ddagger)^{6-r}} \mu^2(m) d_{6-r}(m) \sum_{n \leq X/m} \mu^2(n) d_r(n) \\
& \ll \sum_{r \leq 2} \sum_{m \leq (X^\ddagger)^{6-r}} \mu^2(m) d_{6-r}(m) \left(\frac{X}{m} \right) (\log X)^{r-1}  \\ 
& \ll X \left(\sum_{r \leq 2} (\log X)^{r-1} \right) \left(\sum_{m \leq (X^\ddagger)^{6}} \frac{d_{6}(m)}{m} \right) \\
& \ll X (\log X)\left(\log \exp \left((\log X)^{1/8} \right) \right)^7 \\
& \ll X (\log X)^{15/8}. 
\end{align*}
This is sufficiently small for our purposes. \\ 

We may now assume that $A_{ij} \geq X^\ddagger$ for at least three pairs $i,j$ with $1 \leq i \leq 3, 1 \leq j \leq 2$. We now suppose that there exist $a \ne b$ such that 
\[A_{a,2}, A_{b,1} \geq X^\dagger. \]
The sum over $\BA$ satisfying these properties can be bounded by
\begin{align*} \sum_{\BA} \left \lvert \S_2(X; \BA) \right \rvert & \leq \sum_{x_{ij}, (i,j) \ne (a,2), (b,1)} \prod_{(i,j) \ne (a,2), (b,1)} \left \lvert \sum_{x_{a,2}} \sum_{x_{b,1}} \alpha_{(a,2)} \beta_{(b,1)} \left(\frac{x_{a,2}}{x_{b,1}} \right) \right \rvert,
\end{align*}
where $\alpha, \beta$ have modulus at most one. Lemma \ref{double} then applies, and since our variables $x_{a,2}, x_{b,1}$ range over intervals exceeding $X^\dagger$ in length, it follows that 
\[|\S_2(X;\BA)| \ll \left(\prod_{(i,j) \ne (a,2), (b,1)} A_{ij} \left(A_{a,2} A_{b,1} \left(A_{a,2}^{-1/3} + A_{b,1}^{-1/3} \right) \right) \right) \ll X (X^\dagger)^{-1/3} = O \left(X (\log X)^{-3} \right), \]
which is again enough. \\

Next consider the family where the two previous conditions do not hold, and in addition there exist $a \ne b$ such that $2 \leq A_{b,1} \leq X^\dagger$ and $A_{a,2} > X^\ddagger$. Under these conditions, we see that
\[|\S_2(X; \BA)| \ll \sum_{x_{ij}, (i,j) \ne (a,2), (b,1)} \sum_{x_{a,2}} \left \lvert \sum_{x_{b,1}} \mu^2 \left(\prod_{(i,j) \ne (a,2), (b,1)} x_{ij} \right) \left(\frac{x_{a,2}}{x_{b,1}} \right) \right \rvert,   \]
where $A_{ij} \leq x_{ij} \leq \Delta A_{ij}$ and $\omega(x_{ij}) \leq \Omega$ for $1 \leq i \leq 3, 1 \leq j \leq 2$. Now put $\ell = \omega(x_{a,2})$, writing 
\[x_{a,2} = p_1 \cdots p_\ell \]
with $p_1 < p_2 < \cdots < p_\ell$ we obtain
\[|\S_2(X; \BA)| \ll \sum_{\substack{x_{ij} \\ (i,j) \ne (a,2), (b,1)}} \sum_{x_{b,1}} \sum_{0 \leq \ell \leq \Omega} \left \lvert \sum_{\omega(x_{a,2}) = \ell} \mu^2\left(\prod_{i,j} x_{ij} \right) \left(\frac{x_{a,2}}{x_{b,1}} \right) \right \rvert,\]
the inner sum being bounded by 
\[\sum_{p_1 \cdots p_{\ell-1}} \left \lvert \sum_{p_\ell} \left(\frac{p_\ell}{x_{b,1}} \right) \right \rvert \] 
and $p_1, \cdots, p_\ell$ satisfy $A_{a,2} \leq p_1 \cdots p_{\ell} \leq \Delta A_{a,2}$. Note that
\[p_\ell \geq A_{a,2}^{1/\ell} \geq \exp \left((\log X)^{1/9} \right).  \]
We may now apply Lemma \ref{S-W} to obtain the bound 
\[\left \lvert \sum_{p_\ell} \left(\frac{p_\ell}{x_{b,1}} \right) \right \rvert \ll_A A_{b,1}^{1/2} \frac{A_{a,2}}{p_1 \cdots p_{\ell-1}} (\log X)^{-A/9} + \Omega, \]
with $A$ arbitrarily large. Note that $p_1 \cdots p_{\ell-1} \leq X$, hence
\[\sum_{p_1 \cdots p_{\ell-1} \leq X} (p_1 \cdots p_{\ell-1})^{-1} \ll \sum_{n \leq X} \frac{1}{n} \ll \log X. \]
Hence 
\[\sideset{}{^{(3)}} \sum_\BA |\S_2(X; \BA)| \ll A_{b,1}^{1/2}\prod_{i,j} A_{ij} (\log X)^{-A/9 + 1} \ll X (\log X)^{-A/9 + 11/2}. \]
Choosing $A$ large shows that this contribution is negligible. \\

The remaining case can be summarized by the following properties:
\begin{enumerate}
    \item $\prod_{i,j} A_{ij} \leq \Delta^{-6} X$; 
    \item $A_{ij} \geq X^\ddagger$ for at least three pairs of indices $(i,j)$; 
    \item If $A_{ij}, A_{k\ell} \geq X^\dagger$ then $j = \ell$; 
    \item If $A_{ij} \leq A_{k \ell}$ with $j \ne \ell$, then either $A_{ij} = 1$ or $2 \leq A_{ij} \leq X^\dagger$ and $A_{k \ell} < X^\ddagger$. 
\end{enumerate}

We now show that the second option in (4) cannot happen. This will imply that we have accounted for all possibilities for (\ref{char sum}), and hence reduced our problem to estimating $\S_1(X)$. \\

Suppose, without loss of generality, that $2 \leq A_{11} \leq X^\dagger$ and $A_{22} < X^\ddagger$. Since $A_{ij} \geq X^\ddagger$ for at least three pairs of indices $(i,j)$, one of $A_{12}$ or $A_{32}$ must exceed $X^\ddagger$. We then have $A_{11} \leq X^\dagger$ and $A_{32}$, say, exceeds $X^\ddagger$, which means that our earlier estimation covers this case. \\

The upshot now is that
\begin{equation} \Sigma_2(X) \ll_A X (\log X)^{15/8}
\end{equation}
for some $\kappa(A) > 0$. It follows from (\ref{par sum}) that
\begin{align*} \S_2(X) & = X^{-1} \Sigma_2(X) + \int_1^X \Sigma_2(t) \frac{dt}{t^2} \\
& \ll (\log X)^{15/8} + \int_1^X \frac{(\log t)^{15/8} dt}{t} \\
& = (\log X)^{23/8},
\end{align*}
which is sufficiently small for our purposes. \\

Finally, we may evaluate the main term, which is given by (\ref{S1X}). By the triangle inequality, we have
\[\S_1(X) \ll \sum_{x_1 x_2 x_3 \leq X} \frac{1}{x_1 x_2 x_3} = \sum_{n \leq X} \frac{d_3(n)}{n}\]
which is $O((\log X)^3)$. This completes the proof of the Proposition. 

\subsubsection{Lower bounds}

For the lower bound, we shall assume 
\[1\leq |y_1 y_2 y_3| \leq T^\delta\]
where $\delta$ is some explicit positive number which we shall specify later. We then consider $x_1, x_2$ satisfying 

\begin{equation} \label{box bd} |x_i y_i^2| \leq c_i T^{1/4} |y_1 y_2 y_3|^{1/2}, i = 1,2\end{equation} 
where $c_1, c_2$ are two small positive numbers. Note that 
\[|x_3 y_3^2| = |x_1 y_1^2 + x_2 y_2^2| \leq |x_1 y_1^2| + |x_2 y_2^2| \leq (c_1 + c_2) T^{1/4} |y_1 y_2 y_3|^{1/2}, \]
whence
\[|x_1 x_2 x_3| (y_1 y_2 y_3)^2 \leq c_3 T^{3/4} |y_1 y_2 y_3|^{3/2} \]
where $c_3 = c_1 c_2 (c_1 + c_2)$. Thus 
\[|(x_1 y_1^2) x_1 x_2 x_3| \leq (c_1 T^{1/4} |y_1 y_2 y_3|^{1/2}) (c_3 T^{3/4} |y_1 y_2 y_3|^{-1/2}) \]
which is less than $T$ provided that $c_1 c_3 \leq 1$. Therefore every pair $(x_1, x_2)$ satisfying (\ref{box bd}) with $x_1, x_2$ both square-free and $x_3 = (y_1^2 x_1 + y_2^2 x_2) y_3^{-2} \in \bZ$ square-free will contribute to $N(T)$. \\

We now count pairs $(x_1, x_2)$ such that
\begin{enumerate}
    \item $(x_1, x_2)$ satisfies (\ref{box bd}); 
    \item $\gcd(x_1, x_2) = 1$; 
    \item $x_1, x_2$ are square-free; and
    \item $y_1^2 x_1 + y_2^2 x_2 \equiv 0 \pmod{y_3^2}$, $(y_1^2 x_1 + y_2^2 x_2)y_3^{-2}$ is square-free. 
\end{enumerate}
For each prime $p$, we interpret conditions (2) to (4) modulo $p^2$. Condition (2) is the assertion that $p | x_1 \Rightarrow p \nmid x_2$, Condition (3) is the assertion that for all primes $p$ we have $p^2 \nmid x_1, x_2$, and Condition (4) is stating $y_3^2 | y_1^2 x_1 + y_2^2 x_2$, and if $p^{s} || y_3$, then $p^{2s + 2} \nmid y_1^2 x_1 + y_2^2 x_2$. Let 
\[\rho_\By(m) = \# \{(x_1, x_2) \pmod{m} : (2) \text{ to } (4) \text{ holds for all } p | m\}. \]
It is apparent that $\rho_\By(\cdot)$ is multiplicative. Put 
\[N^\ast(\By; T) = \# \{(x_1, x_2) \in \bZ^2 : (1) \text{ to } (4) \text{ hold} \} \]
and 
\[N_b^\ast(\By; T) = \#\{(x_1, x_2) \in \bZ^2 : (\ref{box bd}) \text{ holds, }(2) \text{ to } (4) \text{ holds mod } b \}\]
By standard arguments, we have 
\[N^\ast(\By; T) = \prod_{p \leq Y} \left(1 - \frac{\rho_\By(p^{2k})}{p^{2k}} \right) \frac{T^{1/2} }{|y_1 y_2 y_3|} + O \left(\sum_{Y < p < T^{1/8} |y_1 y_2 y_3|^{1/4} \max\{|y_1|^{-1},|y_2|^{-1}\}} \left(\frac{T^{1/2}}{p^2 |y_1 y_2 y_3|} + 1 \right) \right), \]
the error term being bounded by 
\[O \left(\frac{T^{1/2}}{Y |y_1 y_2 y_3|} + \frac{T^{1/8} |y_1 y_2 y_3|^{1/2}}{\min\{|y_1|, |y_2|\}} \right). \]
Since $|y_1 y_2 y_3| \leq T^\delta$, we obtain an acceptable error term provided that $\delta < 1/4$. This shows that 
\[N(T) \gg \sum_{1 \leq |y_1 y_2 y_3| \leq T^\delta} N^\ast(\By; T) \gg \sum_{1 \leq |y_1 y_2 y_3| \leq T^\delta} \frac{T^{1/2}}{|y_1 y_2 y_3|}.\]
Since
\[\sum_{1 \leq |y_1 y_2 y_3| \leq Z} |y_1 y_2 y_3|^{-1} \gg \sum_{n \leq Z} d_3(n) n^{-1} \gg  (\log Z)^3, \]
this confirms the lower bound. 

\subsection{Counting points with respect to the canonical height when $\chi(\fX) < 0$}

In this section, we first prove that the number of \emph{quadratic points} on a hyperelliptic curve given by the model
\[C_F: z^2 = F(x,y)\]
where $F$ is an integral, non-singular binary form having degree $2g+2$ with $g \geq 2$, is dominated by the ``obvious" points given by triples $(x,y, \sqrt{F(x,y)})$. To show that the \emph{proper} quadratic points are negligible, we note that when $g = 2$ the proper quadratic points, which come in conjugate pairs, are in bijection with the rational points of the Jacobian $\Jac(C_F)$ via the correspondence $[P] \mapsto [P_1 + P_2] - K_{C_F}$, where $K_{C_F}$ is the canonical divisor. Thus in this case the proper quadratic points of bounded height are given by the rational points of bounded height in $\Jac(C_F)(\bQ)$, for which there are $O_F((\log T)^{r_F})$ many, where $r_F$ is the Mordell-Weil rank of $\Jac(C_F)$. For $g \geq 3$ the proper quadratic points are finite by Faltings' theorem. Thus, the number of quadratic points on $C_F$ is asymptotically equal to the number of rational points in $\bP_\bQ^1$ of bounded height. \\ \\
To the contrary, for $\fX = \fX(\bP^1 : (\Ba, \Bm))$ with $\chi(\fX) < 0$, we get a much less reasonable result. This is because we have little control over the set of integers $x,y$ such that $\ell_i(x,y)$ is divisible by a large square for $i = 1, \cdots, n$. Even with the $abc$-conjecture there is only so much that can be shown. In the case when $\Bm = (2, \cdots, 2)$ we have the following:

\begin{theorem} \label{abc thm}  Let $\fX = \fX(\bP^1 : (\Ba, \Bm))$ be a stacky curve with $\Bm = \underbrace{(2, \cdots, 2)}_n$ with $n \geq 5$. Let $N_{\Ba, n}(T)$  be the number of rational points on $\fX$ satisfying $H_{(\Ba, \Bm)}(x,y) \leq T$. Assume that the $abc$-conjecture holds. Then for any $\ep > 0$ we have
\[N_{(\Ba, \Bm)}(T) \ll_\ep T^{\frac{1}{n-3} + \ep}.\]
\end{theorem} 

\begin{proof} This is similar to the proof of Theorem \ref{nearnorth}. We conclude from that proof that
\[\prod_{i=1}^n |x_i y_i| \geq \rad \left(\prod_{i=1}^n \ell_i(x,y) \right) \gg_\ep \max\{|x|, |y|\}^{n-2 - \ep},\]
and
\[\prod_{i=1}^n |x_i y_i^2| \ll \max\{|x|, |y|\}^n\]
by the triangle inequality. Comparing, we conclude that
\[\prod_{i=1}^n |y_i| \ll_\ep \max\{|x|, |y|\}^{2 + \ep}\]
and in turn 
\[\prod_{i=1}^n |x_i| \gg_\ep \max\{|x|, |y|\}^{n-4 - \ep}.\]
It follows that 
\[H_{(\Ba, \Bm)}(x,y) = \max\{|x|, |y|\}^{n-4} \prod_{i=1}^n |x_i| \gg_\ep \max\{|x|, |y|\}^{2n-8 - \ep}. \]
Hence $N_{(\Ba, \Bm)}(T)$ is bounded by the number of rational points in $\bP_\bQ^1$ have height at most $O_\ep \left(T^{\frac{1}{2n-8 - \ep}}\right)$, which is $O_\ep \left(T^{\frac{1}{n-4-\ep}} \right)$. By adjusting $\ep$. we see that
\[N_{(\Ba, \Bm)}(T) \ll_\ep T^{\frac{1}{n-4} + \ep}.\]
\end{proof}

\begin{remark} We do not expect the upper bound given in Theorem \ref{abc thm} to be sharp. Indeed, the bound we obtain essentially comes from the scenario that for almost all integers $m \ll_\ep T^{2+\ep}$ that there exist $x,y$ with $\max\{|x|, |y|\} \ll_\ep T^{2 + \ep}$ with $Q(x,y)$ divisible by $m^2$. We expect that this should not be the case. 
\end{remark}

\subsection{Hasse principle for integral points when $\Bm = (2, 2, 2)$}

We now consider the question of whether the Hasse principle holds for integral points on stacky curves of the shape $\fX = \fX(\bP_\bQ^1 : (a_1, a_2, a_3), (2,2,2))$. By Theorem \ref{height MT} it suffices to consider when it is possible for the stacky part of the height to be equal to one. This is tantamount to requiring the existence of co-prime integers $x,y$ and integers $y_1, y_2, y_3$ for which
\[|\ell_i(x,y)| = y_i^2 \text{ for } i = 1,2,3.\]
Here, as we recall, $\ell_i(x,y) = \alpha_i x - \beta_i y$, with $a_i = [\alpha_i : \beta_i]$. For $i = 1,2$ we obtain a system of linear equations
\[\begin{bmatrix} \alpha_1 & - \beta_1 \\ \alpha_2 & -\beta_2 \end{bmatrix} \begin{bmatrix} x \\ y \end{bmatrix} = \begin{bmatrix} y_1^2 \\ y_2^2 \end{bmatrix}.\]
Inverting, we find that
\[\begin{bmatrix} x \\ y \end{bmatrix} = \frac{1}{\alpha_1 \beta_2 - \alpha_2 \beta_1} \begin{bmatrix}  \beta_2 & - \beta_1 \\ \alpha_2 & -\alpha_1 \end{bmatrix} \begin{bmatrix} y_1^2 \\ y_2^2 \end{bmatrix}.\]
It follows that 
\[(\alpha_1 \beta_2 - \alpha_2 \beta_1) y_3^2 = \alpha_3 \left(\beta_2 y_1^2 - \beta_1 y_2^2 \right) - \beta_3 \left(\alpha_2 y_1^2 - \alpha_1 y_2^2 \right) \]
which we can write as 
\[(\alpha_2 \beta_3 - \alpha_3 \beta_2) y_1^2 - (\alpha_1 \beta_3 - \alpha_3 \beta_1) y_2^2 + (\alpha_1 \beta_2 - \alpha_2 \beta_1) y_3^2 = \begin{vmatrix} y_1^2 & y_2^2 & y_3^2 \\ \alpha_1 & \alpha_2 & \alpha_3 \\ \beta_1 & \beta_2 & \beta_3 \end{vmatrix} =  0.\]
Therefore, the existence of the integers $y_1, y_2, y_3$, and hence $x,y$, depends on whether this conic has a rational point.

\end{document}